\ifx\shlhetal\undefinedcontrolsequence\let\shlhetal\relax\fi


\input amssym.def
\input amssym.tex


\def\phi{\varphi}

\def\pmb#1{\setbox0=\hbox{#1}%
\leavevmode
\kern-.025em\copy0\kern-\wd0
\kern.05em\copy0\kern-\wd0
\kern-.025em\raise.0433em\box0 }


%


\def\calP{{\cal P}}





%

\def\fr#1{{\frak #1}}
        \def\fra{\fr{a}}
        
        \def\frc{\fr{c}}
        \def\frd{\fr{d}}

        \def\fri{\fr{i}}


%

\font\tenrm=cmr10
\font\ninerm=cmr9

\font\sevenrm=cmr7
\font\sixrm=cmr6
\font\fiverm=cmr5
%

\font\tensl=cmsl10
\font\ninesl=cmsl9

%

\font\tenit=cmti10
\font\nineit=cmti9

%

\font\tenbf=cmbx10
\font\ninebf=cmbx9

\font\sevenbf=cmbx7
\font\sixbf=cmbx6
\font\fivebf=cmbx5
%
\font\tensmc=cmcsc10
\font\ninesmc=cmcsc9

%

\font\tenmi=cmmi10
\font\ninemi=cmmi9

\font\sevenmi=cmmi7
\font\sixmi=cmmi6
\font\fivemi=cmmi5
%

\font\tensy=cmsy10
\font\ninesy=cmsy9

\font\sevensy=cmsy7
\font\sixsy=cmsy6
\font\fivesy=cmsy5
%
\font\tenex=cmex10
\font\nineex=cmex9
%

\font\tenmsb=msbm10
\font\ninemsb=msbm9

\font\sevenmsb=msbm7
\font\sixmsb=msbm6
\font\fivemsb=msbm5
%

%

\font\teneufm=eufm10
\font\nineeufm=eufm9

\font\seveneufm=eufm7
\font\sixeufm=eufm6
\font\fiveeufm=eufm5

\def\tenpoint{%
\lineskip=2pt\baselineskip=18pt\lineskiplimit=0pt%
\def\rm{\fam0\tenrm}%
\textfont0=\tenrm \scriptfont0=\sevenrm \scriptscriptfont0=\fiverm%
\textfont1=\tenmi \scriptfont1=\sevenmi \scriptscriptfont1=\fivemi%
\textfont2=\tensy \scriptfont2=\sevensy \scriptscriptfont2=\fivesy%
\textfont3=\tenex \scriptfont3=\nineex \scriptscriptfont3=\nineex%
\textfont\itfam=\tenit \def\it{\fam\itfam\tenit}%
\textfont\slfam=\tensl \def\sl{\fam\slfam\tensl}%
\textfont\bffam=\tenbf \scriptfont\bffam=\sevenbf%
\scriptscriptfont\bffam=\fivebf \def\bf{\fam\bffam\tenbf}%
\textfont\msbfam=\tenmsb \scriptfont\msbfam=\sevenmsb%
\scriptscriptfont\msbfam=\fivemsb \def\Bbb{\fam\msbfam\tenmsb}%
\textfont\eufmfam=\teneufm \scriptfont\eufmfam=\seveneufm%
\scriptscriptfont\eufmfam=\fiveeufm \def\frak{\fam\eufmfam\teneufm}%
\def\smc{\tensmc}
\rm}

\def\ninepoint{%
\lineskip=1pt\baselineskip=16pt\lineskiplimit=0pt%
\def\rm{\fam0\ninerm}%
\textfont0=\ninerm \scriptfont0=\sixrm \scriptscriptfont0=\fiverm%
\textfont1=\ninemi \scriptfont1=\sixmi \scriptscriptfont1=\fivemi%
\textfont2=\ninesy \scriptfont2=\sixsy \scriptscriptfont2=\fivesy%
\textfont3=\nineex \scriptfont3=\nineex \scriptscriptfont3=\nineex%
\textfont\itfam=\nineit \def\it{\fam\itfam\nineit}%
\textfont\slfam=\ninesl \def\sl{\fam\slfam\ninesl}%
\textfont\bffam=\ninebf \scriptfont\bffam=\sixbf%
\scriptscriptfont\bffam=\fivebf \def\bf{\fam\bffam\ninebf}%
\textfont\msbfam=\ninemsb \scriptfont\msbfam=\sixmsb%
\scriptscriptfont\msbfam=\fivemsb \def\Bbb{\fam\msbfam\ninemsb}%
\textfont\eufmfam=\nineeufm \scriptfont\eufmfam=\sixeufm%
\scriptscriptfont\eufmfam=\fiveeufm \def\frak{\fam\eufmfam\nineeufm}%
\def\smc{\ninesmc}%
\rm}


\catcode`\@=11

\newbox\secbox
\newbox\keybox
\newbox\itembox

\newskip\sectionskipamount

\newif\ifextra \extratrue

\newdimen\itemindent \itemindent=25pt

\def\item#1{\oldindent=\parindent\parindent=\itemindent%
\par\hangindent\parindent\indent\llap{\rm #1\enspace}%
\parindent=\oldindent\ignorespaces}
\def\footnoteindent{20pt}
\newinsert\footins
\skip\footins=12pt plus 2pt minus 4pt
\dimen\footins=30pc
\count\footins=1000
\def\footnoterule{\kern-3pt\hrule width 2in\kern 2.6pt}
\def\footnote#1#2{\unskip\edef\@sf{\spacefactor\the\spacefactor}%
{$^{#1}$}\@sf
\insert\footins{\parindent=\footnoteindent\ninepoint
\itemindent=\footnoteindent
\interlinepenalty=100 \let\par=\endgraf
\leftskip=0pt \rightskip=0pt
\baselineskip=11pt 
\splittopskip=12pt plus 1pt minus 1pt \floatingpenalty=20000
\smallskip\item{#1}{#2}}}

\sectionskipamount=20 pt plus 4 pt minus 4 pt
\def\sectionbreak{\par\ifdim\lastskip<\sectionskipamount
\removelastskip\vskip\sectionskipamount \fi\penalty-200}
%
\def\head#1 #2\endhead\par{\sectionbreak
   \message{#1 #2}
   \setbox\secbox=\hbox{\bf#1\enspace}\ignorespaces
   \vbox{\hangindent=\wd\secbox\noindent
   \hangafter=1\box\secbox{\bf\ignorespaces #2}\smallskip}\ignorespaces
   \nobreak\vskip-\parskip\noindent}

\def\subhead#1\endsubhead{\proclaim@{\smc#1\unskip\ifextra.\fi%
\extratrue}{}{\quad\rm}}
\def\subsubhead#1\endsubsubhead{\proclaim@{\sl#1\unskip\ifextra.\fi%
\extratrue}{}{\quad\rm}}

\def\introduction\par{\sectionbreak\noindent{\bf Introduction}
\par\noindent}

%
\def\proclaim@#1#2#3{\ifdim\lastskip<\medskipamount\removelastskip
\medskip\fi\penalty-55\noindent#1#2#3\ignorespaces}

%
\def\proclaim #1{\proclaim@{\smc #1:\enspace}{}\sl}

%
\def\proclaiminfo #1#2{\proclaim@{\smc #1\enspace}{\rm(#2):}%
{\enspace\sl}}

\def\endproclaim{\par\rm
 \ifdim\lastskip<\medskipamount\removelastskip\penalty55\medskip\fi}

%
\def\demo#1{\removelastskip\medskip\par\noindent{\sl #1:}\enspace}

\def\bull{\vrule height 6pt width 4pt depth 0.47pt} 
\def\qed{\qquad\bull\vadjust{\medskip}\ifmmode{}\else\par\fi}

%
\def\rem#1{\proclaim@{\sl#1:}{}{\quad\rm}}

%
\def\reminfo#1#2{\proclaim@{\sl#1:}{\enspace\sl#2.}{\quad\rm}}


\def\subdemo#1{\proclaim@{\smc#1\unskip\ifextra:\fi\extratrue}{}%
{\quad\rm}}

%
\def\subdemoinfo#1#2{\proclaim@{\smc#1\ifextra:\fi\extratrue}%
{\enspace\sl #2.}{\quad\rm}}

%
\def\|{|\;}

\def\cf{\rm}
\def\colon{{:}\;}
\def\comment#1\endcomment{}
\def\condition(#1){%
\setbox\itembox=\hbox{(#1)\enspace}\itemindent=\wd\itembox%
\item{\rm(#1)}%
}

\def\dotcup{
 \def\dotcupD{\cup\hskip-7.5pt\cdot\hskip4.5pt}
 \def\dotcupS{\cup\hskip-4pt\raise0.6pt\hbox{$\cdot$}}
 \mathchoice{\dotcupD}{\dotcupD}{\dotcupS}{}}
\def\dotunion{
\def\dotunionD{\bigcup\kern-9pt\cdot\kern5pt}
\def\dotunionT{\bigcup\kern-7.5pt\cdot\kern3.5pt}
\mathop{\mathchoice{\dotunionD}{\dotunionT}{}{}}}

\def\firstcondition (#1){\hangindent\parindent{\rm(#1)}\enspace
 \ignorespaces}
\def\freeprod{
\def\freeprodD{{\prod\kern-10.9pt{*}\kern5pt}}
\def\freeprodT{\prod\kern-9pt{*}\kern2pt}
\mathop{\mathchoice{\freeprodD}{\freeprodT}
{not defined}{not defined}}
}

\def\hefresh{
   \def\hefreshD{\mathop{\raise1.5pt\hbox{${\smallsetminus}$}}}
   \def\hefreshS{\mathop{\raise0.85pt\hbox{$\scriptstyle\smallsetminus$}}}
   \mathchoice{\hefreshD}{\hefreshD}{\hefreshS}{\hefreshS}}

\def\leaderfill{\leaders\hbox to 1em{\hss.\hss}\hfill}

\def\newline{\hfill\break}
\def\newpage{\vfill\eject}

\def\qed{\qquad\bull\vadjust{\medskip}\ifmmode{}\else\par\fi}
\def\remark#1{}
\outer\def\shalom{\bigskip\rightline\today\par\vfill\supereject\end}

\def\shvor{\discretionary{}{}{}}
\def\sub#1{_{\lower2pt\hbox{$\scriptstyle #1$}}}
\def\sdp{\hbox to 7.56942pt{$\times$\hskip-0.25em{\raise 0.3ex
\hbox{$\scriptscriptstyle|$}}}}


\def\thanks{\proclaim@{\smc Acknowledgement:\quad}{}{\rm}}

%
\newdimen\keywidth
\newdimen\oldindent
\def\references#1{
\setbox\keybox=\hbox{[#1]\enspace}\keywidth=\wd\keybox
\oldindent=\parindent \parindent=\keywidth
\ninepoint\frenchspacing\rm
   \ifdim\lastskip<\bigskipamount\removelastskip\bigskip\fi
   \filbreak%
\centerline{\bf References}\medskip\nobreak}
\def\ref[#1]{\par\smallskip\hang\indent\llap{\hbox to\keywidth%
    {[#1]\hfil\enspace}}\ignorespaces}
\def\endreferences{\tenpoint\nonfrenchspacing\parindent=\oldindent}


\def\and{{\rm and}}

\def\mod{\;\mathop{\rm mod}\hskip0.2em\relax}


\def\bd{{\rm bd}}


\def\today{\number\day\space \ifcase\month\or January\or February\or
March\or April\or May\or June\or July\or August\or September%
\or October\or November\or December\fi, \space\number\year}

\def\endletterhead{\endgroup \par\bigskip\sl\dspace\parindent=25pt
  \footline={\ifnum\pageno=1\hfil \else\hss\tenrm\folio\hss \fi}}
\def\endletterheadgermanstyle{\endgroup
  \par\bigskip\sl\dspace\parindent=0pt\parskip=\medskipamount
  \footline={\ifnum\pageno=1\hfil \else\hss\tenrm\folio\hss \fi}}
\def\sincerely#1{\par\bigskip \hbox to \hsize {\null \hfil #1 \hfil}}


\def\dspace{\lineskip=2pt\baselineskip=18pt\lineskiplimit=0pt}

\hsize=16 true cm
\tenpoint
\parindent=25pt
\footline={\ifnum\pageno=0\hfil \else\hss\tenrm\folio\hss\fi}

\catcode`\@=\active





\def\calP{{\cal P}}






\def\fr#1{{\frak #1}}
        \def\fra{\fr{a}}
        
        \def\frc{\fr{c}}
        \def\frd{\fr{d}}

        \def\fri{\fr{i}}



\def\case#1:#2.{\par\medskip \noindent {\smc Case #1:} {\sl #2.}\quad}
\def\cf{}
\def\claim#1:#2.{\par\medskip \noindent {\smc Claim #1:} {\sl #2.}\quad}
\def\Claim:#1.{\par\medskip\noindent{\smc Claim:} {\sl #1.}\quad}
\def\colon{{:}\;}
\def\comment#1\endcomment{}
\def\condition(#1){\item{\rm#1}}

\def\dotcup{
 \def\dotcupD{\cup\hskip-7.5pt\cdot\hskip4.5pt}
 \def\dotcupS{\cup\hskip-4pt\raise0.6pt\hbox{$\cdot$}}
 \mathchoice{\dotcupD}{\dotcupD}{\dotcupS}{}}
\def\dotunion{
\def\dotunionD{\bigcup\kern-9pt\cdot\kern5pt}
\def\dotunionT{\bigcup\kern-7.5pt\cdot\kern3.5pt}
\mathop{\mathchoice{\dotunionD}{\dotunionT}{}{}}}

\def\finishproclaim {\par\rm
 \ifdim\lastskip<\medskipamount\removelastskip\penalty55\medskip\fi}
\def\firstcondition (#1){\hangindent\parindent{\rm(#1)}\enspace
 \ignorespaces}
\def\freeprod{\mathop{\prod\kern-10.9pt{*}\kern5pt}}
\def\freeprodT{\prod\kern-9pt{*}\kern2pt}

\def\setminus{
   \def\hefreshD{\mathop{\raise1.5pt\hbox{${\smallsetminus}$}}}
   \def\hefreshS{\mathop{\raise0.85pt\hbox{$\scriptstyle\smallsetminus$}}}
   \mathchoice{\hefreshD}{\hefreshD}{\hefreshS}{\hefreshS}}

\def\leaderfill{\leaders\hbox to 1em{\hss.\hss}\hfill}

\def\longmapright#1#2{\smash{\mathop
                {\hbox to #1pt{\rightarrowfill}}\limits^{#2}}}

\def\newline{\hfill\break}
\def\newpage{\vfill\eject}

\def\part#1:#2.{\par\medskip \noindent {\smc Part #1:} {\sl #2.}\quad}
\def\proclaimauthor #1 (#2):{\medbreak\noindent{\smc#1}
{\rm(#2):}\enspace\sl\ignorespaces}
\def\proofof#1:{\par\medskip\noindent{\sl Proof of {\rm #1:}}}

\def\qed{\qquad\bull\vadjust{\medskip}\ifmmode{}\else\par\fi}
\def\remark#1{}
\def\rmapdown#1
   {\Big\downarrow\rlap{$\vcenter{\hbox{$\scriptstyle#1$}}$}}
\outer\def\shalom{\bigskip\rightline\today\par\vfill\supereject\end}

\def\shvor{\discretionary{}{}{}}
\def\sub#1{_{\lower2pt\hbox{$\scriptstyle #1$}}}
\def\sdp{\hbox to 7.56942pt{$\times$\hskip-0.25em{\raise 0.3ex
\hbox{$\scriptscriptstyle|$}}}}

\def\pmb#1{\setbox0=\hbox{#1}%
\leavevmode
\kern-.025em\copy0\kern-\wd0
\kern.05em\copy0\kern-\wd0
\kern-.025em\raise.0433em\box0 }


\def\and{{\rm and}}

\def\mod{\;{\rm mod}\hskip0.2em\relax}



\def\sqed#1{\qquad\bull\hbox{$_{#1}$}\vadjust{\medskip}\ifmmode{}%
\else\par\fi}
\def\pcf{{\rm pcf}}

\def\tcf{{\rm tcf}}
\def\mod{\mathop{\rm mod}}
\def\sup{\mathop{\rm sup}}
\def\otp{\mathop{\rm otp}}
\def\cf{\mathop{\rm cf}}
\def\Min{{\rm min}}
\def\Max{{\rm Max}}
\def\reg{{\rm reg}}
\def\hcf{{\rm hcf}}
\def\htcf{{\rm htcf}}
\def\atom{{\rm atom}}
\def\Depth{{\rm Depth}}

\def\Dom{{\rm Dom}}
\def\dom{{\rm dom}}
\def\cov{{\rm cov}}
\def\Ord{{\rm Ord}}
\def\eub{{\rm eub}}

\def\Rang{{\rm Rang}}
\def\DEPTH{{\rm DEPTH}}

\def\Car{{\rm Car}}
\def\Card{{\rm Card}}
\def\lub{{\rm lub}}

\def\wsat{{\rm wsat}}
\def\wlog{\mathop{\rm wlog}}

\def\Prod{{\prod}}
\def\trans#1#2{#1\;^\kappa\!#2}
\def\transa#1#2{#1\;^\theta\!#2}
\def\inv{{\rm inv}}
\def\newagei{{\fri}}
\def\Ga{{\fra}}
\def\lg{{\ell g}}


last revision - February 1995

\vskip3cm
\centerline{THE PCF THEOREM REVISITED}

\medskip

\centerline{by}
\centerline{Saharon Shelah\footnote*{Partially supported by the
Deutsche Forschungsgemeinschaft, grant Ko 490/7-1.
Publication no. 506.}}
\centerline{Institute of Mathematics, The Hebrew University,
Jerusalem, Israel}
\centerline{Department of Mathematics, Rutgers University, New
Brunswick, NJ, USA}
\bigskip
\centerline{\bf Dedicated to Paul Erd\H{o}s}
\bigskip

\rem{Abstract}
{\it The $\pcf$ theorem (of the possible cofinability theory) was
proved for reduced products $\prod_{i<\kappa} \lambda_i/I$, where
$\kappa< \min_{i<\kappa} \lambda_i$. Here we prove this theorem under
weaker assumptions such as $\wsat(I)< \min_{i<\kappa} \lambda_i$, where
$\wsat(I)$ is the minimal $\theta$ such that $\kappa$ cannot be
delivered to $\theta$ sets $\notin I$ (or even slightly weaker
condition). We also look at the existence of 
exact upper bounds relative to $<_I$ ($<_I-\eub$) 
as well as cardinalities
of reduced products and the cardinals $T_D(\lambda)$. Finally we apply this
to the problem of the depth of ultraproducts (and reduced products) of
Boolean algebras.}


\newpage

\head \S0 Introduction \endhead

The aim of the $\pcf$ theory is to answer the question,
what are the possible cofinalities ($\pcf$) of the partial orders
$\prod\limits_{i<\kappa} \lambda_i/I$, where $\cf(\lambda_i)=\lambda_i$, for
different ideals $I$ on $\kappa$. For a quick introduction to the $\pcf$
theory see [Sh400a], and for a detailed exposition, see [Sh-g] and
more history. 
In \S1 and \S2 we generalize the basic theorem of this theory by
weakening the assumption $\kappa< \min_{i<\kappa} \lambda_i$ to the
assumption that $I$ extends a fixed ideal $I^*$ with $\wsat(I^*)<
\min_{i<\kappa}\lambda_i$, where $\wsat(I^*)$ is the minimal $\theta$
such that $\kappa$ cannot be divided to $\theta$ sets $\notin I^*$ (not
just that the Boolean algebra $\calP(\kappa)/I^*$ has no $\theta$ pairwise
disjoint non zero elements).
So \S1, \S2 follow closely [Sh-g, Ch. I=Sh345a], [Sh-g, II 3.1],
[Sh-g, VIII \S1].
It is interesting to note that some of those proofs which look to be
superceded when by [Sh420, \S1] we know that for regular
$\theta<\lambda$, $\theta^+<\lambda \Rightarrow \exists$ stationary
$S\in I[\lambda]$, $S\subseteq \{\delta<\lambda: \cf(\delta)=\theta\}$,
give rise to proofs here which seem neccessary. Note $\wsat(I^*)\leq
|\Dom(I^*)|^+$ (and $\reg_*(I^*)\leq |\Dom(I^*)|^+$ so [Sh-g, I \S1,
\S2, II \S1, VII 2.1, 2.2, 2.6] are really a special case of the proofs here.

During the sixties the cardinalities of ultraproducts and
reduced products were much investigated (see Chang and Keisler
[CK]).
For this the notion ``regular filter'' (and $(\lambda,
\mu)$-regular filter) were introduced, as: if $\lambda_i\ge \aleph_0$,
$D$ a regular ultrafilter (or filter) on $\kappa$ then
$\prod\limits_{i<\kappa}\lambda_i /D=(\lim\inf_D\lambda_i)^\kappa$.
We reconsider these problems in \S3 (again continuing [Sh-g]).
We also draw a conclusion on the depth of the reduced product of
Boolen algebras partially answering a problem of Monk; and make it clear
that the truth of the full
excpected result is translated to a problem on pcf.
On those problems on Boolean algebras see Monk [M].
In this section we include known proofs for completeness (mainly 3.6).

Let us review the paper in more details. In 1.1, 1.2 we give basic
definition of cofinality, true cofinality, $\pcf(\bar \lambda)$ and
$J_{<\lambda}[\bar \lambda]$ where usually $\bar \lambda=\langle
\lambda_i: i< \kappa\rangle$ a sequence of regular cardinals, 
$I^*$ a fixed ideal
on $\kappa$ such that we consider only ideals extending it (and filter
disjoint to it). Let $\wsat (I^*)$ be the first $\theta$ such that we
cannot partition $\kappa$ to $\theta$ $I^*$-positive set (so they are
pairwise disjoint, not just disjoint modulo $I^*$). In 1.3, 1.4 we
give the basic properties. In lemma 1.5 we phrase the basic property
enabling us to do anything: (1.5$(*)$): $\lim\inf_{I^*} (\bar
\lambda)\geq \theta \geq \wsat(I^*)$, $\prod \bar \lambda/ I^*$ is
$\theta^+$-directed; we prove that $\prod \bar \lambda/
J_{<\lambda}[\bar \lambda]$ is $\lambda$-directed. In 1.6, 1.8 we deduce
more properties of $\langle J_{<\lambda}[\bar \lambda]: \lambda\in
\pcf(\bar \lambda)\rangle$ and in 1.7  deal with $<_{J_{<\lambda}[\bar
\lambda]}$-increasing sequence $\langle f_\alpha: \alpha<\lambda\rangle$
with no $<_{J_{<\lambda}[\bar \lambda]}$-bound in $\prod \bar \lambda$.
In 1.9 we prove $\pcf(\bar \lambda)$ has a last element and in 1.10,
1.11 deal with the connection between the true cofinality of
$\prod\limits_{i<\kappa} \lambda_i/ D^*$ and $\prod\limits_{i<\sigma}
\mu_i/E$ when $\mu_i=:\tcf(\prod\limits_{i<\kappa} \lambda_i / D_i)$ and
$D^*$ is the $E$-limit of the $D_i$'s.

In 2.1 we define normality of $\lambda$ for $\bar \lambda$: $J_{\leq
\lambda} [\bar \lambda] = J_{<\lambda} [\bar \lambda] + B_\lambda$ and
we define 
semi-normality: $J_{\leq} [\bar \lambda]= J_{<\lambda} [\bar \lambda] +
\{B_\alpha: \alpha< \lambda\}$ where $B_\alpha / J_{<\lambda}[\bar
\lambda]$ is increasing. We then (in 2.2) characterize semi normality
(there is a $<_{J_{<\lambda}[\bar \lambda]}$-increasing $\bar f=\langle
f_\alpha: \alpha<\lambda\rangle$ cofinal in $\prod \bar \lambda/D$ for
every ultrafilter $D$ (disjoint to $I^*$ of course) such that 
$\tcf(\prod \bar \lambda/D)=\lambda$)
and when semi normality implies normality (if some such $\bar f$ has a
$<_{J_{<\lambda}[\bar \lambda]} - \eub$).

We then deal with continuity system $\bar a$  and $<_{J_{<\lambda}[\bar
\lambda]}$- increasing sequence obeying $\bar a$, in a way adapted to
the basic assumption $(*)$ of 1.5.

Here as elsewhere if $\min(\bar \lambda)\geq \theta^+$ our life is
easier than when we just assume $\lim sup_{I^*} (\bar \lambda)\geq
\theta$, $\prod \bar \lambda /I^*$ is $\theta^+$-directed (where
$\theta\geq \wsat (I^*)$ of course). In 2.3 we give the definitions, in
2.4 we quote existence theorem, show existence of obedient sequences (in
2.5), essential uniqueness (in 2.7) and better consequence to 1.7 (in
the direction to normality). We define (2.9) generating sequence and
draw a conclusion (2.10(1)). Now we get some desirable properties: in
2.8 we prove semi normality, in 2.10(2) we compute $\cf(\prod\bar
\lambda/ I^*)$ as $\max\pcf(\bar \lambda)$. Next we relook at the whole
thing: define several variants of the $\pcf$-th (Definition 2.11).
Then (in 2.12)
we show that e.g. if $\min(\bar \lambda)> \theta^+$, we get the
strongest version (including normality using 2.6, i.e. obedience).
Lastly 
we try to map the implications between the various properties when we do
not use the basic assumption 1.5 $(*)$ (in fact there are considerable
dependence, see 2.13, 2.14).

In 3.1, 3.3 we present measures of regularity of filters, in 3.2 we
present measures of hereditary cofinality of $\prod \bar \lambda/D$:
allowing to decrease $\bar \lambda$ and/or increase the filter. 
In 3.4 - 3.8 we try
to estimate reduced products of cardinalities $\prod\limits_{i<\kappa}
\lambda_i/D$ and in 3.9 we give a reasonable upper bound by hereditary
cofinality ($\leq (\theta^\kappa /D+ \hcf_{D, \theta}
(\prod\limits_{i<\kappa}\lambda_I))^{<\theta}$ when $\theta \geq
\reg_{\otimes}(D)$).

In 3.10-3.11 we return to existence of $\eub$'s and obedience
(Saharon, new point over 2.6) and in 3.12 draw conclusion on ``downward
closure".

In 3.13 - 3.14 we estimate $T_D(\bar \lambda)$ and in 3.15 try to
translate it more fully to pcf problem (countable cofinality is somewhat
problematic (so we restrict ourselves to $T_D(\bar
\lambda)>\mu=\mu^{\aleph_0}$). We also mention $\aleph_1$-complete
filters; (3.16, 3.17) and see 
what can be done without relaying on $\pcf$ (3.20)).

Now we deal with depth: define it (3.18, see 3.19), give lower bound
(3.22), compute it for ultraproducts of interval Boolean algebras of
ordinals (3.24). Lastly we translate the problem ``does $\lambda_i<
\Depth^+ (B_i)$ for $i<\kappa$ implies $\mu<
\Depth^+(\prod\limits_{i<\kappa} B_i/D)$" at least when $\mu>
2^\kappa$ and $(\forall \alpha<\mu)[|\alpha|^{\aleph_0}<\mu]$,
to a $\pcf$ problem (in 3.26).

In the last section we phrase a reason 1.5$(*)$ works (see 4.1), analyze
the case we weaken to 1.5$(*)$ to $\lim\inf_{I^*}(\bar \lambda)\geq
\theta\geq \wsat (I^*)$ proving the pseudo $\pcf$-th (4.4).


\head \S1 Basic pcf\endhead

\rem{Notation 1.0}
$I, J$ denote ideals on a set $\Dom(I)$, $\Dom(J)$ resp., called its
domain (possibly $\bigcup_{A\in I}A\subset\Dom I$).
If not said otherwise the domain is an infinite cardinal
denoted by $\kappa$ and also the ideal is proper i.e. $\Dom(I)\notin I$.
Similarly $D$ denotes a filter on a set $\Dom D$; we do not
always distinguish strictly between an ideal on $\kappa$ and the
dual filter on $\kappa$.
Let $\bar\lambda$ denote a sequence of the form
$\langle\lambda_i\colon i<\kappa\rangle$.
We say $\bar\lambda$ is regular if every $\lambda_i$ is regular,
$\Min{\bar\lambda}=\Min\{\lambda_i\colon i<\kappa\}$ (of course
also in $\bar\lambda$ we can replace $\kappa$ by another set),
and let $\prod{\bar\lambda}=\prod\limits_{i<\kappa}\lambda_i$; usually we are
assuming $\bar \lambda$ is regular.
Let $I^*$ denote a fixed ideal on $\kappa$.
Let $I^+=\calP(\kappa)\setminus I$ (similarly $D^+=\{
A\subseteq\kappa\colon \kappa\setminus A\not\in D\}$), let
$$\lim\inf_I{\bar\lambda}=\min\{\mu\colon \{
i<\kappa\colon\lambda_i\le\mu\}\in I^+\}\quad \hbox{ and}$$
$$\lim\sup_I{\bar\lambda}=\Min\{\mu\colon\{ i<\kappa\colon
\lambda_i >\mu\}\in I\}\quad\hbox{ and}$$
$$\atom_I{\bar\lambda}=\{\mu\colon\{i\colon\lambda_i=\mu\}\in
I^+\}.$$
For a set $A$ of ordinals with no last element, $J_A^{\bd}=\{
B\subseteq A\colon\sup (B)<\sup (A)\}$, i.e. the ideal of bounded
subsets.
Generally, if $\inv(X)=\sup\{ |y|\colon \vDash\varphi(X, y)\}$ then
$\inv^+(X)=\sup\{ |y|^+\colon \vDash\varphi[X,y]\}$, and any $y$ such
that
$\vDash\varphi[X,y]$ is a wittness for $|y|\le\inv(X)$ (and
$|y|<\inv^+(X)$), and it exemplifies this.
Let $\bar A^*_\theta[{\bar\lambda}]=\langle
A^*_\alpha\colon\alpha<\theta\rangle=\langle
A^*_{\theta,\alpha}[{\bar\lambda}]\colon\alpha<\theta\rangle$ be
defined by: $A^*_\alpha=\{ i<\kappa\colon\lambda_i>\alpha\}$.
Let $\Ord$ be the class of ordinals.

\rem{Definition 1.1}
\item{(1)} For a partial order
\footnote\dag{actually we do not require
$p\leq q \leq p \Rightarrow p=q$ so we should say quasi partial order}
$P$:

\itemitem{(a)} $P$ is $\lambda$-directed if: for every
$A\subseteq P$, $|A|<\lambda$ there is $q\in P$ such that
$\bigwedge_{p\in A}p\le q$, and we say: $q$ is an upper bound of
$A$;

\itemitem{(b)} $P$ has true cofinality $\lambda$ if there is
$\langle p_\alpha\colon \alpha<\lambda\rangle$ cofinal in $P$, i.e.:
$\bigwedge_{\alpha<\beta}p_\alpha<p_\beta$ and
$\forall q\in P[\bigvee_{\alpha<\lambda}q\le p_\alpha]$ [and
one writes $\tcf(P)=\lambda$ for the minimal such $\lambda$]
(note: if $P$ is linearly ordered it always has a true cofinality but
e.g.\ $(\omega, <)\times(\omega_1, <)$ does not).

\itemitem{(c)} $P$ is called endless if $\forall p\in P\exists q\in
P[q>p]$ (so if $P$ is endless, in clauses (a), (b), (d) above we can
replace $\le$ by $<$).

\itemitem{(d)} $A\subseteq P$ is a cover if: $\forall p\in
P\exists q\in A[p\le q]$; we also say ``$A$ is cofinal in $P$''.

\itemitem{(e)} $\cf(P)=\min\{ |A|\colon A\subseteq P$ is a
cover$\}$.

\itemitem{(f)} We say that, in $P$, $p$ is a lub (least upper
bound) of $A\subseteq P$ if:

\hskip30pt ($\alpha$)\ $p$ is an upper bound of $A$ (see (a))

\hskip30pt ($\beta$)\ if $p'$ is an upper bound of $A$ then
$p\le p'$.

\item{(2)} If $D$ is a filter on $S$, $\alpha_s$ (for $s\in S$)
are ordinals, $f$, $g\in\prod\limits_{s\in S}\alpha_s$, then:
$f/D<g/D$, $f<_D g$ and $f<g\mod D$ all mean $\{ s\in S\colon
f(s)<g(s)\}\in D$. Also if $f$, $g$ are partial functions from $S$ to
ordinals, $D$ a filter on $S$ $\underline{\rm then}$ $f<g\ \mod\ D$
means $\{i\in \Dom(f): i\notin \Dom(f)\ \underline{or}\ f(i)< g(i)
\hbox{ (so both are defined)}\}$ belongs to $D$. We write $X=A\ \mod\ D$
if $\Dom(D)\setminus [(X\setminus A) \cup (A\setminus X)]$ belongs to $D$.
Similarly for $\le$, and we do not distinguish between a filter
and the dual ideal in such notions.
So if $J$ is an ideal on $\kappa$ and $f$,
$g\in\prod{\bar\lambda}$, then $f<g \mod J$ iff
$\{ i<\kappa\colon\lnot f(i)<g(i)\}\in J$.
Similarly if we replace the $\alpha_s$'s by partial orders.

\item{(3)} For $f$, $g\colon S\to$ Ordinals, $f<g$ means
$\bigwedge_{s\in S}f(s)<g(s)$; similarly $f\le g$.
So $(\prod{\bar\lambda}, \le)$ is a partial order, we denote it
usually by $\prod\bar\lambda$; similarly $\prod f$ or
$\prod\limits_{i<\kappa} f(i)$.

\item{(4)} If $I$ is an ideal on $\kappa$,
$F\trans \subseteq\Ord$, we call $g\trans \in\Ord$ an
$\le_I$-eub (exact upper bound) of $F$ if:
\itemitem{($\alpha$)} $g$ is an $\le_I$-upper bound of $F$ (in
$^\kappa\Ord$)

\itemitem{($\beta$)} if $h\trans \in\Ord$, $h<_I\Max\{g, 1\}$
$\underline{\rm then}$ for some $f\in F$, $h<\max\{f, 1\}\mod I$.

\itemitem{($\gamma$)} if $A\subseteq\kappa$, $A\not=\emptyset\mod
I$ and $[f\in F\Rightarrow f\upharpoonright A=_I 0_A$, i.e.\ $\{
i\in A\colon f(i)\not=0\}\in I]$ then $g\upharpoonright
A=_J 0_A$.

\itemitem{(5)\quad (a)}
We say the ideal $I$ (on $\kappa$) is $\theta$-weakly saturated
if $\kappa$ cannot be divided to $\theta$ pairwise disjoint
sets from $I^+$ (which is ${\cal P}(\kappa)\setminus I$)

\itemitem{(b)} $\wsat(I)=\Min\{\theta\colon I$ is
$\theta$-weakly saturated$\}$

\rem{Remark 1.1A}
\item{(1)} Concerning 1.1(4), note: $g'=\Max\{ g, 1\}$ means
$g'(i)=\Max\{g(i), 1\}$ for each $i<\kappa$;
if there is $f\in F$, $\{ i<\kappa\colon f(i)=0\}\in I$ we can
replace $\Max\{ g, 1\}$, $\Max\{ f,1\}$ by $g$, $f$
respectivally in clause ($\beta$) and omit clause ($\gamma$).

\item{(2)} Considering $\prod\limits_{i<\kappa} f(i)$, $<_I$ formally if
$(\exists i)f(i)=0$ then $\prod\limits_{i<\kappa} f(i)=\emptyset$; but
we usually ignore this, particularly when $\{i: f(i)=0\}\in I$.

\rem{Definition 1.2}
Below if $\Gamma$ is ``the ultrafilter disjoint to $I$", we write $I$ instead
$\Gamma$.
\item{(1)} For a property $\Gamma$ of ultrafilters (if $\Gamma$
is the empty condition, we omit it):
$$
\pcf_\Gamma(\bar\lambda)= \pcf(\bar\lambda, \Gamma)=
\{\tcf(\prod \bar\lambda/D)\colon D\;
\hbox{is an ultrafilter
on } \kappa\;\hbox{satisfying }\Gamma\}
$$
(so $\bar\lambda$ is a sequence of ordinals, usually of regular
cardinals, note: as $D$ is an
ultrafilter, $\Prod \bar\lambda/D$ is linearly ordered hence has true
cofinality).

\item{(1A)} More generally, for a property $\Gamma$ of ideals on
$\kappa$ we let $\pcf_\Gamma(\bar\lambda)=\{\tcf(\Prod \bar\lambda/J)
\colon J$ is an ideal on $\kappa$ satisfying $\Gamma$ such that
$\Prod \bar\lambda/J$ has true cofinality$\}$. Similarly below.

\item{(2)} $J_{<\lambda}[{\bar\lambda}, \Gamma]=\{ B\subseteq\kappa$:
for no ultrafilter $D$ on $\kappa$ satisfying $\Gamma$ to which $B$
belongs, is $\tcf(\Prod \bar\lambda/D)\ge\lambda\}$.

\item{(3)} $J_{\le\lambda}[\bar\lambda, \Gamma]=
J_{<\lambda^+}[\bar\lambda, \Gamma]$.

\item{(4)} $\pcf_\Gamma({\bar\lambda}, I)=\{\tcf(\Prod
\bar\lambda/D)\colon D$ a filter on $\kappa$ disjoint to $I$
satisfying $\Gamma\}$.

\item{(5)} If $B\in I^+$,
$\pcf_I({\bar\lambda}\upharpoonright B)=\pcf_{I+(\kappa\setminus B)}
({\bar\lambda})$ (so if $B\in I$ it is $\emptyset$), also
$J_{<\lambda}({\bar\lambda}\upharpoonright B, I)\subseteq\calP(B)$
is defined similarly.

\item{(6)}
If $I=I^*$ we may omit it, similarly in (2), (4).
\item{} If $\Gamma=\Gamma_{I^*}=\{ D\colon D$ an ultrafilter on $\kappa$
disjoint to $I^*\}$ we may omit it.

\rem{Remark}
We mostly use $\pcf(\bar\lambda)$, $J_{<\lambda}[\bar\lambda]$.

\rem{Claim 1.3}
\item{(0)} $(\Prod\bar\lambda, <_J)$ and $(\Prod\bar\lambda, \le_J)$
are endless (even when each $\lambda_i$ is just a limit ordinal);

\item{(1)} $\min(\pcf_I(\bar\lambda))\ge\lim\inf_I(\bar\lambda)$ for
$\bar \lambda$ regular; 

\itemitem{(2)\quad (i)} If $B_1\subseteq B_2$ are from
$I^+$ $\underline{\rm then}$ $\pcf_I(\bar\lambda\upharpoonright
B_1)\subseteq\pcf_I(\bar\lambda\upharpoonright B_2)$;

\itemitem{(ii)} if
$I\subseteq J$ $\underline{\rm then}$ $\pcf_J(\bar\lambda)\subseteq\pcf_I(\bar\lambda)$;
and
\itemitem{(iii)} for $B_1$, $B_2\subseteq\kappa$ we have
$\pcf_I(\bar\lambda\upharpoonright
(B_1\cup B_2))=\pcf_I(\bar\lambda\upharpoonright
B_1)\bigcup\pcf_I(\bar\lambda\upharpoonright B_2)$.
Also
\itemitem{(iv)} $A\in J_{<\lambda}
[\bar\lambda\upharpoonright (B_1\cup B_2)]
\Leftrightarrow  A\cap B_1\in
J_{<\lambda} [\bar\lambda\upharpoonright B_1]\ \&\ A\cap B_2\in
J_{<\lambda} [\bar\lambda\upharpoonright B_2]$
\itemitem{(v)} If $A_1$, $A_2\in I^+$, $A_1\cap A_2 = \emptyset$,
$A_1\cup A_2=\kappa$, and $\tcf(\prod\bar \lambda \restriction A_\ell,
<_I)=\lambda$ for $\ell=1, 2$ then $\tcf(\prod\bar \lambda,
<_I)=\lambda$; and if the sequence $\bar f=\langle f_\alpha:
\alpha<\lambda\rangle$ wittness both assumptions then it wittness the
conclusion. 

\itemitem{(3)\quad (i)} if $B_1\subseteq B_2\subseteq\kappa$, $B_1$
finite and $\bar\lambda$ regular $\underline{\rm then}$
$$
\pcf_I(\bar\lambda\upharpoonright
B_2)\setminus\Rang(\bar\lambda\upharpoonright B_1)
\subseteq\pcf_I(\bar\lambda\upharpoonright (B_2\setminus
B_1))\subseteq \pcf_I(\bar\lambda\upharpoonright B_2)
$$

\itemitem{(ii)} if in addition $i\in B_1\Rightarrow\lambda_i
<\Min (\Rang [\bar\lambda\upharpoonright(B_2\setminus B_1)])$,
\itemitem{} $\underline{\rm then}$ $\pcf_I(\bar\lambda\upharpoonright
B_2)\setminus\Rang(\bar\lambda\upharpoonright
B_1)=\pcf_I(\bar\lambda\upharpoonright(B_2\setminus B_1))$.

\item{(4)} Let $\bar\lambda$ be regular (i.e.\ each $\lambda_i$
is regular);
\itemitem{(i)} If $\theta=\lim\inf_I\bar\lambda$ $\underline{\rm then}$
$\prod\bar\lambda /I$ is $\theta$-directed

\itemitem{(ii)} If $\theta=\lim\inf_I\bar\lambda$ is singular
$\underline{\rm then}$ $\prod\bar\lambda /I$ is $\theta^+$-directed

\itemitem{(iii)} If $\theta=\lim\inf_I\bar\lambda$ is inaccessible
(i.e. a limit regular cardinal), the set $\{i<\kappa: \lambda_i=\theta\}$ is
in the ideal $I$ and for some club $E$ of $\theta$,
$\{i<\kappa\colon \lambda_i\in E\}\in I$
$\underline{\rm then}$
$\prod\bar\lambda /I$ is $\theta^+$-directed.
We can weaken the assumption to ``$I$ is not medium normal for
$\bar\lambda$'' (defined in the next sentence).
Let ``$I$ is not weakly normal for
$(\theta,\bar\lambda)$'' mean: for some $h\in\prod\bar\lambda$, for
no $j<\theta$ is $\{i<\kappa\colon
\lambda_i\leq\theta\Rightarrow h(i)<j\} =\kappa\mod
J$;
and let ``$I$ is not medium normal for $(\theta,\bar\lambda)$''
mean: for
some $h\in\prod\bar\lambda$, for no $\zeta<\lim\inf_I(\bar\lambda)=
\theta$, is $\{
i<\kappa\colon\lambda_i\leq\theta\Rightarrow h(i)<\zeta\}\in J^+$.

\itemitem{(iv)} If $\{i: \lambda_i=\theta\}=\kappa\ \mod\ I$ and $I$
is weakly normal then $(\prod\bar \lambda, <_I)$ has true cofinality
$\theta$.

\itemitem{(v)} If $\prod\bar\lambda /I$ is $\theta$-directed then
$\cf(\prod\bar\lambda /I)\ge\theta$ and $\min\pcf_I(\prod {\bar
\lambda})\geq \theta$.

\itemitem{(vi)} $\pcf_I(\bar \lambda)$ is non empty set of regular
cardinals. 
[see part (7)]

\item{(5)} Assume $\bar\lambda$ is regular 
and: if
$\theta'=:\limsup _I(\bar\lambda)$ is regular then $I$ is not
weakly normal for $(\theta',\bar\lambda)$.
$\underline{\rm Then}$
$\pcf_I(\bar\lambda)\not\subseteq(\lim\sup_I(\bar\lambda))^+$;
in fact for some ideal $J$ extending $I$, $\prod\bar\lambda/J$ is
$(\lim\sup_I(\bar\lambda))^+$-directed.

\item{(6)}
If $D$ is a filter on a set $S$ and for $s\in S$, $\alpha_s$ is a
limit ordinal $\underline{\rm then}$:
\itemitem{(i)} $\cf(\Prod_{s\in S}\alpha_s, <_D)=\cf(\Prod_{s\in
S}\cf(\alpha_s), <_D)= \cf(\Prod_{s\in S}(\alpha_s, <)/ D)$, and

\itemitem{(ii)} $\tcf(\Prod_{s\in S}\alpha_s, <_D)=\tcf(\Prod_{s\in
S}(\cf(\alpha_s), <_D))=\tcf(\Prod_{s\in S}(\alpha_s, <)/D)$.

In particular, if one of them is well defined, $\underline{\rm then}$
so are the others.
This is true even if we replace $\alpha_s$ by linear orders or
even partial orders with true confinality.

\item{(7)} If $D$ is an ultrafilter on a set $S$, $\lambda_s$ a
regular cardinal, $\underline{\rm then}$
$\theta =: \tcf(\Prod_{s\in S}\lambda_s,
<_D)$ is well defined and $\theta\in \pcf(\{\lambda_s\colon s\in S\})$.

\item{(8)} If $D$ is a filter on a set $S$, for $s\in S$,
$\lambda_s$ is a regular cardinal,
$S^*=\{\lambda_s\colon s\in S\}$ and
$$E=:\{B\colon\ \ \ B\subseteq S^*\ \hbox{ and }\ \{s\colon
\lambda_s\in B\}\in D\}$$
and $\lambda_s>|S|$ or at least
$\lambda_s>|\{t: \lambda_t=\lambda_s\}|$ for any $s\in S$ 
$\underline{\rm then}$:
\itemitem{(i)} $E$ is a filter on $S^*$, and if $D$ is an
ultrafilter on $S$ $\underline{\rm then}$ $E$ is an ultrafilter on $S^*$.

\itemitem{(ii)} $S^*$ is a set of regular cardinals and 
\itemitem{} if $s\in S\Rightarrow \lambda_s>|S|$ $\underline{\rm then}$
$(\forall \lambda\in S^*) \lambda>|S^*|$,

\itemitem{(iii)} $F=\{ f\in\Prod_{s\in S}\lambda_s\colon
s=t\Rightarrow f(s)=f(t)\}$ is a cover of $\Prod_{s\in
S}\lambda_s$,

\itemitem{(iv)} $\cf(\Prod_{s\in S}\lambda_s /D)=\cf(\Prod S^* /E)$
and $\tcf(\Prod_{s\in S}\lambda_s /D)=\tcf(\Prod S^* /E)$.

\item{(9)} Assume $I$ is an ideal on $\kappa$,
$F\trans \subseteq\Ord$ and $g\trans \in\Ord$.
If $g$ is a $\leq_I$-eub of $F$ $\underline{\rm then}$ $g$ is a
$\leq_I$-lub of $F$.

\item{(10)} $\sup\pcf_I(\bar\lambda)\le|\prod\bar\lambda /I|$

\item{(11)} If $I$ is an ideal on $S$ and $(\prod\limits_{s\in S}
\alpha_s, <_I)$ has true cofinality $\lambda$ as exemplified by $\bar
f=\langle f_\alpha: \alpha< \lambda\rangle$ $\underline{\rm then}$ the
function $\langle \alpha_s: s\in S\rangle$ is a $<_I-\eub$ (hence
$<_I-\lub$) of $\bar f$.

\item{(12)} The inverse of (11) holds: $\underline{\rm if}$ $I$ is an
ideal on $S$ and $f_\alpha\in {}^S\Ord$ for
$\alpha<\lambda=\cf(\lambda)$, $\langle f_\alpha: \alpha<\lambda\rangle$
is $<_I$-increasing with $<_I-\eub$ $f$ 
\item{} $\underline{\rm then}$
$\tcf(\prod\limits_{i} f(i), <_I)=\tcf(\prod \cf[f(i)], <_I)=\lambda$.

\item{(13)} If $I\subseteq J$ are ideals on $\kappa$ $\underline{\rm then}$

\itemitem{(a)} $\wsat (I)\geq \wsat(J)$

\itemitem{(b)} $\lim\inf_I(\bar \lambda)\leq \lim\inf_J(\bar \lambda)$

\itemitem{(c)} if $\lambda=\tcf(\prod\limits_{i<\kappa}\lambda_i,
<_I)$ then $\lambda=\tcf(\prod\limits_{i<\kappa}\lambda_i, <_J)$

\item{(14)} If $f_1$, $f_2$ are $<_I-\lub$ of $F$ $\underline{\rm
then}$ $f_1=_I f_2$.
\demo{Proof}
They are all very easy, e.g.

(0) We shall show $(\Prod\bar\lambda, <_J)$ is endless (assuming,
of course, that $J$ is a proper ideal on $\kappa$).
Let $f\in\Prod\bar\lambda$, then $g=:f+1$ (defined
$(f+1)(\gamma)=f(\gamma)+1)$ is in $\Prod\bar\lambda$ too as each
$\lambda_\alpha$ being an infinite cardinal is a limit ordinal
and $f<g\mod J$.

(5) Let $\theta' =: \lim\sup_I (\bar \lambda)$ and define
$$
J=:\{A\subseteq \kappa: \hbox{ for some }\theta<\theta', \{i<\kappa:
\lambda_i>\theta \hbox{ and } i\in A\} \hbox{ belongs to }I\}.
$$
Clearly $J$ is an ideal on $\kappa$ extending $I$ (and $\kappa\notin
J_0$) and $\lim\sup_J (\bar \lambda)= \lim \inf_J(\bar
\lambda)=\theta'$.

%

\noindent
$\underline{\rm Case~1}$: $\theta'$ is singular

By part (4) clause (ii), $\prod {\bar \lambda}/ J$ is
$(\theta')^+$-directed and by part (4) clause (iv) we get the
desired conclusion.

\noindent
$\underline{\rm Case~2}$: $\theta'$ is inaccessible ($>\aleph_0$)

Let $h\in \prod \bar \lambda$ wittness that ``$I$ is not weakly normal
for $(\theta', \bar \lambda)$"
and let
$$
J^*=\{ A\subseteq \kappa:
\hbox{ for some }j<\theta'\hbox{ we have }
\{i\in A: h(i)<j\}=A\
\mod\ I\}.
$$

Note that if $A\in J$ then for some $\theta<\theta'$, $A'=:\{i\in A:
\theta_i> \theta\}\in I$ hence
the choice $j=:\theta$ wittness $A\in J^*$. So $J\subseteq J^*$. Also
$J^*\subseteq {\cal P}(\kappa)$ by its definition. $J^*$ is closed under
subsets (trivial) and under union [why? assume $A_0$, $A_1\in J^*$,
$A=A_0\cup A_1$;
choose $j_0$,
$j_1<\theta'$ such that  $A'_\ell=: \{i\in A_\ell: h(i)<j_\ell\}=
A_\ell\ \mod\ I$, so $j=: \max\{j_0, j_1\}<\theta$ and $A'=\{i\in A:
h(i)<j\}= A\ \mod\ I$; so $A\in J^*$]. Also $\kappa\notin J^*$ [why? as
$h$ wittness that $I$ is not weakly normal for $(\theta', \bar \lambda)$].
So together
$J^*$ is an ideal on $\kappa$ extending $I$.
Now $J^*$ is not medium normal for $(\theta, \bar \lambda)$, as
wittnessed by $h$. Lastly $\prod {\bar
\lambda}/ J^*$ is $(\theta')^+$-directed (by part (4) clause (iii)), and
so $\pcf_J(\bar \lambda)$ is disjoint to $(\theta')^+$.

(9) Let us prove $g$ is a $\le_I-\lub$ of $F$ in
$(^\kappa\Ord, \le_I)$.
As we can deal separately with $I+A$, $I+(\kappa\setminus A)$
where $A=:\{i\colon g(i)=0\}$, and the later case is trivial we
can assume $A=\emptyset$.
So assume $g$ is not a $\leq_I-\lub$, so there is an upper bound $g'$ of
$F$, but not $g\le_I g'$.
Define $g''\trans \in\Ord$:
$$
g''(i)= \cases{ 0 & if $g(i)\le g'(i)$\cr
g'(i) & if $g'(i)<g(i)$\cr}.
$$
Clearly $g''<_I g$.
So, as $g$ in an $\leq_I-\eub$ of $F$ for $I$, there is $f\in F$ such that
$g''<_I\max\{ f,1\}$, but $B =:\{i\colon
g'(i)<g(i)\}\not=\emptyset\mod I$ (as not $g\le_I g'$) so
$g'\upharpoonright B=g''\upharpoonright B<_I\max \{
f,1\}\upharpoonright B$. But we know that $f\leq_I g'$ (as $g'$ is an
upper bound of $F$) hence $f\restriction B\leq_I g'\restriction B$, so
by the previous sentence neccessarily
$f\upharpoonright B=_I 0_B$ hence
$g'\upharpoonright B=_I 0_B$; as $g'$ is a $\le_I$-upper bound
of $F$ we know $[f'\in F\Rightarrow f'\upharpoonright B =_I 0_B]$,
hence by $(\gamma)$ of Definition 1.1(4) we have $g\upharpoonright B =_I
0_B$, a contradiction to  $B\notin I$ (see above).
\sqed{1.3}
\rem{Remark 1.3A}
In 1.3 we can also have the straight monotonicity properties of
$$
\pcf_I(\prod {\bar \lambda}, \Gamma).
$$
\proclaim{Claim 1.4}
\item{(1)} $J_{<\lambda}[\bar\lambda]$ is an ideal (of
$\calP(\kappa)$ i.e. on $\kappa$, but the ideal may not be proper).

\item{(2)} if $\lambda\le\mu$, then $J_{<\lambda}[\bar\lambda]\subseteq
J_{<\mu}[\bar\lambda]$

\item{(3)} if $\lambda$ is singular,
$J_{<\lambda}[\bar\lambda]=J_{<\lambda^+}[\bar\lambda]=J_{\le\lambda}
[\bar\lambda]$

\item{(4)} if $\lambda\not\in \pcf(\bar\lambda)$, then
$J_{<\lambda}[\bar\lambda]=J_{\le\lambda}[\bar\lambda]$.

\item{(5)} If $A\subseteq\kappa$, $A\not\in
J_{<\lambda}[\bar\lambda]$, and
$f_\alpha\in\prod\bar\lambda\upharpoonright A$, $\langle
f_\alpha\colon\alpha<\lambda\rangle$ is
$<_{J_{<\lambda}[\bar\lambda]}$-increasing cofinal in
$(\prod\bar\lambda\upharpoonright A)/J_{<\lambda}[\bar\lambda]$
$\underline{\rm then}$ $A\in J_{\le\lambda}[\bar\lambda]$.
Also this holds when we replace $J_{<\lambda}[\bar\lambda]$ by
any ideal $J$ on $\kappa$, $I^*\subseteq J\subseteq
J_{\le\lambda}[\bar\lambda]$.

\item{(6)} The earlier parts hold for $J_{<\lambda}[\bar \lambda,
\Gamma]$ too.
\endproclaim

\demo{Proof} Straight.

\proclaim{Lemma 1.5}
Assume
\item{($*$)} $\bar\lambda$ is regular and
\itemitem{($\alpha$)} $\Min\bar\lambda>\theta\ge \wsat (I^*)$
(see 1.1(5)(b)) or at least

\itemitem{($\beta$)} $\lim\inf_{I^*}(\bar\lambda)\ge\theta\ge
\wsat(I^*)$,
and $\prod\bar\lambda /I^*$ is $\theta^+$-directed\footnote\dag{note if
$\cf(\theta)<\theta$ then ``$\theta^+$-directed'' follows from
``$\theta$-directed'' which follows from \break
``$\lim\inf_{I^*}
(\bar\lambda)\ge \theta$, i.e. first part of clause $(\beta)$. Note also
that if clause $(\alpha)$ holds then $\prod\bar \lambda/ I^*$ is
$\theta^+$-directed (even $(\prod\bar \lambda, <)$ is
$\theta^+$-directed), so clause $(\beta)$ implies clause $(\alpha)$.}.

If
$\lambda$ is a cardinal $\ge\theta$, and $\kappa\not\in J_{<\lambda}
[\bar\lambda]$
$\underline{\rm then}$ $(\Prod\bar\lambda, <_{J_{<\lambda}[\bar\lambda]})$ is
$\lambda$-directed (remember: $J_{<\lambda}[\bar\lambda]=
J_{<\lambda}[\bar\lambda, I^*]$).
\endproclaim

\demo{Proof}
Note: if $f\in\Prod\bar\lambda$ then $f<f+1\in\Prod\bar\lambda$, (i.e.
$(\Prod\bar\lambda, <_{J_\lambda[\bar\lambda]})$ is endless)
where $f+1$ is defined by $(f+1)(i)=f(i) +1$).
Let $F\subseteq\Prod\bar\lambda$, $|F|<\lambda$, and we shall prove that
for some $g\in\Prod\bar\lambda$ we have $(\forall f\in F)(f\le g\mod
J_{<\lambda}[\bar\lambda])$, this suffices.
The proof is by induction on $|F|$.
If $|F|$ is finite, this is trivial.
Also if $|F|\le\theta$, when $(\alpha)$ of $(*)$ holds it is easy:
let $g\in\Prod\bar\lambda$ be $g(i)=\sup\{f(i)\colon
f\in F\}<\lambda_i$; when $(\beta)$ of $(*)$ holds use second clause of
$(\beta)$.
So assume $|F|=\mu$, $\theta<\mu<\lambda$ so let
$F=\{f_\alpha^0\colon \alpha<\mu\}$.
By the induction hypothesis we can choose by induction on
$\alpha<\mu$, $f^1_\alpha\in\Prod\bar\lambda$ such that:
\item{(a)} $f_\alpha^0\le f_\alpha^1\mod J_{<\lambda}[\bar\lambda]$

\item{(b)} for $\beta<\alpha$ we have $f_\beta^1\le f_\alpha^1\mod
J_{<\lambda}[\bar\lambda]$.

If $\mu$ is singular, there is $C\subseteq\mu$ unbounded,
$|C|=\cf(\mu)<\mu$, and by the induction hypothesis there is
$g\in\Prod\bar\lambda$ such that for $\alpha\in C$, $f_\alpha^1\le g\mod
J_{<\lambda}[\bar\lambda]$.
Now $g$ is as required: $f_\alpha^0\le f_\alpha^1\le f^1_{\min(C
\setminus \alpha)}\le g\mod J_{<\lambda}[\bar\lambda]$.
So without loss of generality $\mu$ is regular.
Let us define $A^*_\varepsilon =:\{ i<\kappa\colon\lambda_i>|\varepsilon|\}$
for $\varepsilon<\theta$, so $\varepsilon<\zeta<\theta\Rightarrow
A^*_\zeta\subseteq A^*_\varepsilon$ and $\varepsilon<\theta\Rightarrow
A^*_\varepsilon=\kappa\mod I^*$.
Now we try to define by induction on $\varepsilon<\theta$,
$g_\varepsilon$, $\alpha_\varepsilon=\alpha(\varepsilon)<\mu$ and
$\langle B_\alpha^\varepsilon\colon \alpha< \mu\rangle$ such that:
\item{(i)} $g_\varepsilon\in\Prod\bar\lambda$

\item{(ii)} for $\varepsilon<\zeta$ we have
$g_\varepsilon\upharpoonright
A^*_\zeta\le g_\zeta\upharpoonright A^*_\zeta$

\item{(iii)} for $\alpha<\mu$ let $B^\varepsilon_\alpha
=:\{i<\kappa\colon f^1_\alpha(i)>g_\varepsilon(i)\}$

\item{(iv)} for each $\varepsilon<\theta$, for every
$\alpha\in [\alpha_{\varepsilon+1}, \mu)$,
$B^\varepsilon_\alpha\not= B^{\varepsilon+1}_\alpha\mod
J_{<\lambda}[\bar\lambda]$.

We cannot carry this definition: as letting
$\alpha(*)=\sup\{\alpha_\varepsilon\colon \varepsilon<\theta\}$,
then $\alpha(*)<\mu$
since $\mu=\cf(\mu)>\theta$.
We know that $B^\varepsilon_{\alpha(*)}\cap A^*_{\varepsilon+1}\not=
B^{\varepsilon+1}_{\alpha(*)} \cap A^*_{\varepsilon+1}$ for
$\alpha<\theta$ (by (iv) and as $A^*_{\varepsilon+1}=\kappa
\mod I^*$ and $I^*\subseteq J_{<\lambda}[\bar\lambda]$)
and $B^\varepsilon_{\alpha(*)}\subseteq\kappa$
(by (iii)) and $[\varepsilon<\zeta\Rightarrow B^\zeta_{\alpha(*)}\cap
A^*_\zeta\subseteq
B^\varepsilon_{\alpha(*)}]$ (by (ii)), together $\langle
A^*_{\varepsilon+1}\cap (B^\varepsilon_{\alpha(*)}\setminus
B^{\varepsilon+1}_{\alpha(*)})\colon \varepsilon<\theta\rangle$ is a
sequence of
$\theta$ pairwise disjoint members of $(I^*)^+$, a
contradiction\footnote\dag{i.e we have noted that for no
$B_\varepsilon\subseteq \kappa$ ($\varepsilon<\theta$) do we have:
$B_\varepsilon \neq B_{\varepsilon+1}\ \mod\ I^*$ and
$\varepsilon<\zeta<\theta \Rightarrow B_\varepsilon \cap
A_\zeta\subseteq B_\zeta$ where $A_\zeta=\kappa\ \mod\ I^*$ (e.g.
$A_\zeta=A^*_\zeta$)}  to the
definition of $\theta=\wsat(I^*)$.

\noindent Now for $\varepsilon=0$ let $g_i$ be $f^1_0$ and
$\alpha_\varepsilon=0$.

\noindent
For $\varepsilon$ limit let
$g_\varepsilon(i)=\bigcup_{\zeta<\varepsilon}g_\zeta(i)$ for
$i\in A^*_\varepsilon$ and zero otherwise (note:
$g_\varepsilon\in\Prod\bar\lambda$ as $\varepsilon<\theta$,
$\lambda_i>\varepsilon$ for $i \in A^*_\varepsilon$ and
$\bar\lambda$ is a sequence of regular cardinals) and let
$\alpha_\varepsilon=0$.

\noindent
For $\varepsilon=\zeta+1$, suppose that $g_\zeta$ hence $\langle
B^\zeta_\alpha\colon \alpha<\mu\rangle$ are defined.
If $B^\zeta_\alpha\in J_{<\lambda}[\bar\lambda]$ for unboundedly many
$\alpha<\mu$ (hence for every $\alpha<\mu$) then $g_\zeta$ is an upper
bound for $F\ \mod\ J_{<\lambda}[\bar \lambda]$ and the proof
is complete.
So assume this fails, then there is a minimal $\alpha(\varepsilon)<\mu$
such that $B^\zeta_{\alpha(\varepsilon)}\not\in
J_{<\lambda}[\bar\lambda]$.
As $B^\zeta_{\alpha(\varepsilon)}\not\in J_{<\lambda}[\bar\lambda]$, by
Definition 1.2(2) for some ultrafilter $D$ on $\kappa$ disjoint to
$J_{<\lambda}[\bar\lambda]$ we have $B^\zeta_{\alpha(\varepsilon)}\in D$
and $\cf(\Prod\bar\lambda /D)\ge\lambda$.
Hence $\{ f^1_\alpha /D\colon \alpha<\mu\}$ has an upper bound 
$h_\varepsilon/D$ where $h_\varepsilon\in\Prod\bar\lambda$.
Let us define $g_\varepsilon\in\Prod\bar\lambda$:
$$
g_\varepsilon(i)=\Max\{ g_\zeta(i), h_\varepsilon(i)\}.
$$

Now (i), (ii) hold trivially and $B^\varepsilon_\alpha$ is defined by
(iii).
Why does (iv) hold (for $\zeta$) with $\alpha_{\zeta+1}=
\alpha_\varepsilon =: \alpha(\varepsilon)$?
Suppose $\alpha(\varepsilon)\le \alpha<\mu$.
As $f^1_{\alpha(\varepsilon)}\le f^1_\alpha\mod J_{<\lambda}
[\bar\lambda]$ clearly $B^\zeta_{\alpha(\varepsilon)}\subseteq
B^\zeta_\alpha \mod J_{<\lambda}[\bar\lambda]$.
Moreover $J_{<\lambda}[\bar\lambda]$ is disjoint to $D$ (by its
choice) so $B^\zeta_{\alpha(\varepsilon)}\in D$ implies
$B^\zeta_\alpha\in D$.

\noindent
On the other hand $B^\varepsilon_\alpha$ is $\{i<\kappa\colon
f^1_\alpha(i)>g_\varepsilon(i)\}$ which is equal to $\{i\in
\bar\lambda\colon f^1_\alpha(i)>g_\zeta(i),
h_\varepsilon\shvor(i)\}$
which does not belong to $D$ ($h_\varepsilon$ was chosen such that
$f_\alpha^1\le h_\varepsilon\mod D$).
We can conclude $B^\varepsilon_\alpha\not\in D$, whereas
$B^\zeta_\alpha\in D$;
so they are distinct $\mod J_{<\lambda}[\bar\lambda]$ as required in
clause (iv).

Now we have said that we cannot carry the definition for all
$\varepsilon<\theta$, so we are stuck at some $\varepsilon$;
by the above $\varepsilon$ is successor, say $\varepsilon=\zeta+1$, and
$g_\zeta$ is as required: an upper bound for $F$ modulo
$J_{<\lambda}[\bar \lambda]$.
\sqed{1.5}

\proclaim{Lemma 1.6}
If (*) of 1.5, $D$ is an ultrafilter on $\kappa$ disjoint to $I^*$
and $\lambda=\tcf(\Prod\bar\lambda, <_D)$, $\underline{\rm then}$
for some $B\in D$,
$(\Prod\bar\lambda\upharpoonright B, <_{J_{<\lambda}[\bar\lambda]})$ has
true cofinality $\lambda$.
(So $B\in J_{\le\lambda}[\bar\lambda]\setminus
J_{<\lambda}[\bar\lambda]$ by 1.4(5).)
\endproclaim

\demo{Proof}
By the definition of $J_{<\lambda}[\bar\lambda]$ clearly $D\cap
J_{<\lambda}[\bar\lambda]=\emptyset$.

\noindent
Let $\langle f_\alpha/ D\colon \alpha<\lambda\rangle$ be increasing
unbounded in $\Prod\bar\lambda /D$ (so $f_\alpha\in \Prod\bar\lambda$).
By 1.5 without loss of generality
$(\forall \beta<\alpha)(f_\beta<f_\alpha\mod
J_{<\lambda}[\bar\lambda])$.

\noindent
Now 1.6 follows from 1.7 below: its hypothesis clearly holds.
If $\bigwedge_{\alpha<\lambda}B_\alpha=\emptyset\mod D$, (see (A) of 1.7)
then (see (D) of 1.7) $J\cap D=\emptyset$ hence (see (D) of 1.7)
$g/ D$ contradicts the choice of $\langle f_\alpha/ D\colon
\alpha<\lambda\rangle$.
So for some $\alpha<\lambda$, $B_\alpha\in D$;
by (C) of 1.7 we get the desired conclusion.
\sqed{1.6}

\proclaim{Lemma 1.7}
Suppose (*) of 1.5, $\cf(\lambda)>\theta$,
$f_\alpha\in\Prod\bar\lambda$,
$f_\alpha<f_\beta\mod J_{<\lambda}[\bar\lambda]$ for
$\alpha<\beta<\lambda$, and there
is no $g\in\Prod\bar\lambda$ such that for every $\alpha<\lambda$,
$f_\alpha<g\mod J_{<\lambda}[\bar\lambda]$.

\noindent
$\underline{\rm Then}$ there are $B_\alpha$ (for $\alpha<\lambda$) such that:
\item{(A)} $B_\alpha\subseteq\kappa$ and for some
$\alpha(*)<\lambda\colon \alpha(*)\le\alpha<\lambda\Rightarrow
B_\alpha\not\in J_{<\lambda}[\bar\lambda]$

\item{(B)} $\alpha<\beta\Rightarrow B_\alpha\subseteq B_\beta\mod
J_{<\lambda}[\bar\lambda]$ (i.e. $ B_\alpha\setminus B_\beta\in
J_{<\lambda}[\bar\lambda]$)

\item{(C)} for each $\alpha$, $\langle f_\beta\upharpoonright
B_\alpha\colon
\beta<\lambda\rangle$ is cofinal in $(\Prod\bar\lambda\upharpoonright
B_\alpha,
<_{J_{<\lambda}[\bar\lambda]})$ (better restrict yourselves to
$\alpha\ge \alpha(*)$ (see (A))
so that necessarily $ B_\alpha\not\in J_{<\lambda}[\bar\lambda]$);.

\item{(D)} for some $g\in \Prod\bar\lambda$, $\bigwedge_{\alpha<\lambda}
f_\alpha\le g\mod J$ where\footnote\dag{Of course, if $B_\alpha
=\kappa\mod J_{<\lambda}[\bar\lambda]$, this becomes trivial.}
$J=J_{<\lambda}[\bar\lambda]+\{
B_\alpha\colon \alpha<\lambda\}$;

\noindent
in fact
\item{\quad (D)$^+$} for some $g\in\prod\bar\lambda$ for every $\alpha<\lambda$,
we have$^{\dag{}}$ $f_\alpha\le g\mod (J_{<\lambda}[\bar\lambda]+
B_\alpha)$, in fact $B_\alpha=\{ i<\kappa\colon f_\alpha(i)>g(i)\}$

\item{(E)} if $g\le g'\in\Prod\bar\lambda$, then for arbitrarily large
$\alpha<\lambda$:
$$
\{i<\kappa\colon [g(i)\ge f_\alpha(i)\Leftrightarrow
g'(i)\ge f_\alpha(i)]\}=\kappa\mod J_{<\lambda}[\bar\lambda]
$$
(hence for every large enough $\alpha<\lambda$ this holds)

\item{(F)} if $\delta$ is a limit ordinal $<\lambda$, $f_\delta$
is a $\le_{J_{<\lambda}[\bar\lambda]}-\lub$ of
$\{f_\alpha\colon\alpha <\delta\}$ $\underline{\rm then}$
$ B_\delta$ is a $\lub$
of $\{ B_\alpha\colon\alpha<\delta\}$ in $\calP(\kappa)/
J_{<\lambda}[\bar\lambda]$.
\endproclaim

\demo{Proof of 1.7}
Remember that for $\varepsilon<\theta$, $A_\varepsilon^*=\{
i<\kappa\colon\lambda_i>|\varepsilon|\}$ so $A^*_\varepsilon=\kappa\mod
I^*$
and $\varepsilon<\zeta\Rightarrow A^*_\zeta\subseteq A^*_\varepsilon$.
We now define by induction on $\varepsilon<\theta$, $g_\varepsilon$,
$\alpha(\varepsilon)<\lambda$, $\langle B_\alpha^\varepsilon\colon
\alpha<\lambda \rangle$ such that:
\item{(i)} $g_\varepsilon\in\Prod\bar\lambda$

\item{(ii)} for $\zeta<\varepsilon$, $g_\zeta\upharpoonright
A^*_\varepsilon\le g_\varepsilon\upharpoonright A^*_\varepsilon$

\item{(iii)} $ B_\alpha^\varepsilon =:\{i\in\kappa\colon
f_\alpha(i)>g_\varepsilon(i)\}$

\item{(iv)} if $\alpha(\varepsilon)\le \alpha<\lambda$ then
$ B_\alpha^\varepsilon\not=  B_\alpha^{\varepsilon+1}\mod J_{<\lambda}
[\bar\lambda]$

\noindent
For $\varepsilon=0$ let $g_\varepsilon=f_0$, and
$\alpha(\varepsilon)=0$.

\noindent
For $\varepsilon$ limit let
$g_\varepsilon(i)=\bigcup_{\zeta<\varepsilon}g_\zeta(i)$ if
$i\in A^*_\varepsilon$ and zero otherwise; now
$$
[\zeta<\varepsilon\Rightarrow g_\zeta\upharpoonright A^*_\varepsilon
\le g_\varepsilon\upharpoonright A^*_\varepsilon]
$$ 
holds trivially and
$g_\varepsilon\in\Prod\bar\lambda$ as each $\lambda_i$
is regular and $[i\in A^*_\varepsilon\Leftrightarrow \lambda_i>
\varepsilon]$), and let $\alpha(\varepsilon)=0$.

\noindent
For $\varepsilon=\zeta+1$, if $\{\alpha<\lambda\colon B_\alpha^\zeta\in
J_{<\lambda}[\bar\lambda]\}$ is unbounded in $\lambda$, then $g_\zeta$
is a bound for $\langle f_\alpha\colon \alpha<\lambda\rangle\mod
J_{<\lambda}[\bar\lambda]$, contradicting an assumption.
Clearly
$$
\alpha<\beta<\lambda\ \ \Rightarrow B_\alpha^\zeta\subseteq B_\beta^\zeta\ 
\mod J_{<\lambda}[\bar\lambda]
$$ 
hence $\{\alpha<\lambda\colon B_\alpha^\zeta\in J_{<\lambda}[\bar\lambda]\}$
is an initial segment of $\lambda$. So by the previous sentence there is
$\alpha(\varepsilon)<\lambda$ such that for every $\alpha\in
[\alpha(\varepsilon), \lambda)$, we have $B_\alpha^\zeta
\not\in J_{<\lambda}[\bar\lambda]$ (of course, we may increase
$\alpha(\varepsilon)$ later).
If $\langle B_\alpha^\zeta\colon \alpha<\lambda\rangle$ satisfies the
desired conclusion, with $\alpha(\varepsilon)$ for $\alpha(*)$ in (A)
and $g_\zeta$ for $g$ in (D), (D)$^+$ and (E), we are done.
Now among the conditions in the conclusion of 1.7, (A) holds by
the definition of $ B_\alpha^\zeta$ and of $\alpha(\varepsilon)$,
(B) holds by $ B_\alpha^\zeta$'s definition as $\alpha<\beta\Rightarrow
f_\alpha<f_\beta\mod J_{<\lambda}[\bar\lambda]$,
(D)$^+$ holds with $g=g_\zeta$ by the choice of $B_\alpha^\beta$ hence
also clause (D) follows.
Lastly if (E) fails, say for $g'$, then it can serve as
$g_\varepsilon$.
Now condition (F) follows immediately from (iii) (if (F) fails
for $\delta$, then there is $B\subseteq B_\delta^\zeta$ such that
$\bigwedge_{\alpha<\delta} B_\alpha^\zeta\subseteq B\mod
J_{<\lambda}[\bar\lambda]$ and
$ B_\delta^\zeta\setminus B\not\in J_{<\lambda}[\bar\lambda]$;
now
the function
$g^* =: (g_\zeta\upharpoonright(\kappa\setminus B))\cup
(f_\delta\upharpoonright B)$ contradicts ``$f_\delta$ is a
$\le_{J_{<\lambda}[\bar\lambda]}-\lub$ of $\{f_\alpha\colon
\alpha<\delta\}$'', because: $g^*\in \prod {\bar \lambda}$ (obvious),
$\neg (f_\delta \leq g^*\ \mod\ J_{<\lambda}[{\bar \lambda}])$ [why? as
$B^\zeta_\delta \setminus B\notin J_{<\lambda}[\bar \lambda]$ and
$g^*\restriction (B^\zeta_\delta\setminus B)= g_\zeta\restriction
(B^\zeta_\delta \setminus B)< f_\delta\restriction (B^\zeta_\delta
\setminus B)$ by the choice of $B^\zeta_\delta$], and for
$\alpha<\delta$ we have: 
$$f_\alpha \restriction B
\leq_{J_{<\lambda}[\bar \lambda]} f_\delta \restriction B= g^*
\restriction B\quad\hbox{ and}$$
$$f_\alpha\restriction (\kappa\setminus
B)\leq_{J_{<\lambda}[\bar \lambda]} g_\zeta \restriction
(\kappa\setminus B)= g^*\restriction (\kappa \setminus B)$$
(the $\leq_{J_{<\lambda}[\bar \lambda]}$ holds as $(\kappa\setminus B)\cap
B^\zeta_\alpha \in J_{<\lambda}[\bar \lambda]$ and the definition of
$B^\zeta_\alpha$).
So only clause (C) (of 1.7) may fail, without loss of generality for
$\alpha=\alpha(\varepsilon)$.
I.e. $\langle f_\beta\upharpoonright  B_{\alpha(\varepsilon)}^\zeta
\colon \beta<\lambda\rangle$ is not cofinal in
$(\Prod\bar\lambda\upharpoonright B^\zeta_{\alpha(\varepsilon)},
<_{J_{<\lambda}[\bar\lambda]})$.
As this sequence of functions is increasing w.r.t.
$<_{J_{<\lambda}[\bar\lambda]}$, there is
$h_\alpha\in\prod(\bar\lambda\upharpoonright
B^\zeta_{\alpha(\varepsilon)})$  such that for no
$\beta<\lambda$ do we have
$h_\alpha\le f_\beta\upharpoonright B_{\alpha(\varepsilon)}^j\mod
J_{<\lambda}[\bar\lambda]$.
Let $h'_\varepsilon=h_\varepsilon\cup
0_{(\kappa\setminus B_{\alpha(\varepsilon)}^\zeta)}$ and
$g_\varepsilon\in \Prod\bar\lambda$ be defined by $g_\varepsilon(i)=
\Max\{g_\zeta(i), h'_\varepsilon(i)\}$.
Now define $ B_\alpha^\varepsilon$ by (iii) so (i), (ii), (iii) hold
trivially, and we can check (iv).

So we can define $g_\varepsilon$, $\alpha(\varepsilon)$ for
$\varepsilon<\theta$, satisfying (i)--(iv).
As in the proof of 1.5, this is impossible: because
(remembering $\cf(\lambda)=\lambda>\theta$) letting
$\alpha(*)=:\bigcup_{\varepsilon<\theta}\alpha(\varepsilon)<\lambda$ we
have:
$\langle B^\varepsilon_{\alpha(*)}\cap A^*_\zeta\colon
\varepsilon<\zeta\rangle$ is
$\subseteq$-decreasing, for each $\zeta<\theta$, and
$A^*_\varepsilon=\kappa\mod I^*$ and
$B^{\varepsilon+1}_{\alpha(*)}\not= B^\varepsilon_{\alpha(*)}\mod
J_{<\lambda} [\bar\lambda]$ so
$\langle B^\varepsilon_{\alpha(*)}\cap A^*_{\varepsilon+1}\setminus
B^{\varepsilon+1}_{\alpha(*)}\colon \varepsilon<\theta\rangle$ is a
sequence of $\theta$ pairwise disjoint members of $(J_{<\lambda}[\bar
\lambda])^+$ hence of $(I^*)^+$ which
give the
contradiction to $(*)$ of 1.5; so the lemma cannot fail.
\sqed{1.7}

\proclaim{Lemma 1.8}
Suppose $(*)$ of 1.5.
\item{(1)} For every $ B\in J_{\le\lambda}[\bar\lambda]\setminus
J_{<\lambda}[\bar\lambda]$, we have:
$$
(\prod\bar\lambda\upharpoonright B, <_{J_{<\lambda}
[\bar\lambda]})\;
\hbox{has
true cofinality }\ \lambda\;\hbox{(hence $\lambda$ is regular)}.
$$

\item{(2)} If $D$ is an ultrafilter on $\kappa$, disjoint to $I^*$,
$\underline{\rm then}$ $\cf(\Prod\bar\lambda /D)$ is
$\min\{\lambda\colon D\cap
J_{\le\lambda}[\bar\lambda]\not=\emptyset\}$.

\itemitem{(3)\quad\ \ (i)} For $\lambda$ limit
$J_{<\lambda}[\bar\lambda]=\bigcup_{\mu<\lambda} J_{<\mu}[\bar\lambda]$
hence

\itemitem{(ii)} for every $\lambda$,
$J_{<\lambda}[\bar\lambda]=\bigcup_{\mu<\lambda} J_{\le\mu}
[\bar\lambda]$.

\item{(4)}
$J_{\le\lambda}[\bar\lambda]\not =J_{<\lambda}[\bar\lambda]$
$\underline{\rm iff}$
$J_{\le\lambda} [\bar\lambda]\setminus
J_{<\lambda}[\bar\lambda]\not=\emptyset$
$\underline{\rm iff}$
$\lambda\in\pcf(\bar\lambda)$.

\item{(5)} $J_{\le\lambda}[\bar\lambda] /
J_{<\lambda}[\bar\lambda]$ is $\lambda$-directed (i.e.\ if
$B_\gamma \in J_{\le\lambda}[\bar\lambda]$ for
$\gamma<\gamma^*$, $\gamma^*<\lambda$ $\underline{\rm then}$ for
some $B\in J_{\le\lambda}[\bar\lambda]$ we have
$B_\gamma\subseteq B\mod J_{<\lambda}[\bar\lambda]$ for every
$\gamma<\gamma^*$.)
\endproclaim

\demo{Proof}
(1) Let 

$$\eqalign{J=\{ B\subseteq\kappa\colon & B\in J_{<\lambda}[\bar\lambda]\ \hbox{
or }\ B\in J_{\le\lambda}[\bar\lambda]\setminus
J_{<\lambda}[\bar\lambda]\ \hbox{ and}\cr
\ \ & (\Prod\bar\lambda\upharpoonright B,
<_{J_{<\lambda}[\bar\lambda]})\ \hbox{ has true cofinality }\lambda\}.}$$

\noindent
By its definition clearly $J\subseteq J_{\le\lambda}[\bar\lambda]$;
it is quite easy to
check it is an ideal (use 1.3(2)(v)).
Assume $J\not= J_{\le\lambda}[\bar\lambda]$ and we shall get a
contradiction.
Choose $ B\in J_{\le\lambda}[\bar\lambda]\setminus J$;
as $J$ is an ideal, there is an ultrafilter $D$ on $\kappa$ such
that: $D\cap J=\emptyset$ and $ B\in D$.
Now if $\tcf(\Prod\bar\lambda /D)\ge \lambda^+$, then $ B\not\in
J_{\le\lambda}[\bar\lambda]$ (by the definition of
$J_{\le\lambda}[\bar\lambda]$); contradiction.
On the other hand if $F\subseteq\Prod\bar\lambda$, $|F|<\lambda$ then
there is $g\in\Prod\bar\lambda$ such that $(\forall f\in F)(f<g\mod
J_{<\lambda}[\bar\lambda])$ (by 1.5), so $(\forall f\in F)[f<g\mod D]$
(as $J_{<\lambda}[\bar\lambda]\subseteq J$, $D\cap J=\emptyset$), and
this implies $\cf(\Prod\bar\lambda /D)\ge\lambda$.
By the last two sentences we know that $\tcf(\Prod\bar\lambda/ D)$ is
$\lambda$.
Now by 1.6 for some $C\in D$, $\big(\Prod(\bar\lambda\upharpoonright C),
<_{J_{<\lambda}[\bar\lambda]}\big)$ has true cofinality $\lambda$, of
course $C\cap B \subseteq C$ and $C\cap B\in D$ hence $C\cap B\notin
J_{<\lambda} [\bar \lambda]$.
Clearly if $C'\subseteq C$, $C'\not\in
J_{<\lambda}[\bar\lambda]$ then also $(\Prod\bar\lambda\upharpoonright C',
<_{J_{<\lambda}[\bar\lambda]})$ has true cofinality $\lambda$, hence by
the last sentence
without loss of generality $C\subseteq B$; hence by 1.4(5) we know that
$C\in J_{\le\lambda}[\bar\lambda]$ hence by the definition of $J$
we have $C\in J$.
But this contradicts the choice of $D$ as disjoint from $J$.

We have to conclude that $J=J_{\le\lambda}[\bar\lambda]$ so we have
proved 1.8(1).

(2) Let $\lambda\in\pcf(\bar\lambda)$ be minimal such that $D\cap
J_{\le\lambda}[\bar\lambda]\not=\emptyset$ (it exists as by 1.3(10)
$J_{<(\prod \bar \lambda)^+}[\bar \lambda] = {\cal P}(\kappa)$)
and choose $ B\in D\cap
J_{\le\lambda}[\bar\lambda]$.
So $[\mu<\lambda\Rightarrow B\not\in J_{\le\mu}[\bar\lambda]]$
(by the choice of $\lambda$) hence by 1.8(3)(ii) below, we have
$ B\not\in J_{<\lambda}[\bar\lambda]$.
It similarly follows that
$D\cap J_{<\lambda}[\bar\lambda]=\emptyset$.
Now $(\Prod\bar\lambda\upharpoonright B,
<_{J_{<\lambda}[\bar\lambda]})$ has
true cofinality $\lambda$ by 1.8(1).
As we know that $B\in D\cap J_{\leq \lambda}[\bar \lambda]$, and
$J_{<\lambda}[\bar\lambda]\cap D=\emptyset$;
clearly we have finished the proof.

(3) (i) Let $J=:\bigcup_{\mu<\lambda}J_{<\mu}[\bar\lambda]$.
Now $J$ is an ideal by 1.4(2) and $(\Prod\bar\lambda, <_J)$ is
$\lambda$-directed;
i.e. if $\alpha^*<\lambda$ and
$\{f_\alpha\colon\alpha<\alpha^*\}\subseteq\Prod\bar\lambda$, then there
exists $f\in \Prod\bar\lambda$ such that
$$
(\forall \alpha<\alpha^*)(f_\alpha<f \mod J).
$$
[Why?
if $\alpha^*<\theta^+$ as
$(*)$ of 1.5 holds, this is obvious, suppose not;
$\lambda$ is a limit cardinal, hence there is $\mu^*$
such that $\alpha^*<\mu^*<\lambda$.
Without loss of generality $|\alpha^*|^+<\mu^*$.
By 1.5, there is $f\in\Prod\bar\lambda$ such that
$(\forall\alpha<\alpha^*)(f_\alpha<f\mod J_{<\mu^*}[\bar\lambda])$.
Since $J_{<\mu^*}[\bar\lambda]\subseteq J$, it is immediate that
$$
(\forall\alpha<\alpha^*)(f_\alpha<f \mod J).]
$$
Clearly $\bigcup_{\mu<\lambda} J_{<\mu}[\bar\lambda]\subseteq
J_{<\lambda}[\bar\lambda]$ by 1.4(2).
On the other hand, let us suppose that there is $ B\in
(J_{<\lambda}[\bar\lambda]\setminus\bigcup_{\mu<\lambda}
J_{<\mu}[\bar\lambda])$.
Choose an ultrafilter $D$ on $\kappa$ such that $ B\in D$ and
$D\cap J=\emptyset$.
Since $(\Prod\bar\lambda, <_J)$ is $\lambda$-directed and $D\cap
J=\emptyset$, one has $\tcf(\Prod\bar\lambda/ D)\ge\lambda$, but
$ B\in D\cap J_{<\lambda}[\bar\lambda]$, in contradiction to
Definition 1.2(2).

(3)(ii) If $\lambda$ limit --- by part (i) and 1.4(2); if
$\lambda$ successor --- by 1.4(2) and Definition 1.2(3).

(4) Easy.

(5) Let $\langle f_\alpha^\gamma\colon\alpha<\lambda\rangle$ be
$<_{J_{<\lambda}[\bar\lambda]+(\kappa\setminus
B_\gamma)}$-increasing and cofinal in $\Prod\bar\lambda$ (for
$\gamma<\gamma^*$).
Let us choose by induction on $\alpha<\lambda$ a function
$f_\alpha\in\prod\bar\lambda$, as a
$<_{J_{<\lambda}[\bar\lambda]}$-bound to $\{ f_\beta\colon
\beta<\alpha\}\cup\{ f_\alpha^\gamma\colon\gamma<\gamma^*\}$,
such $f_\alpha$ exists by 1.5 and apply 1.7 to $\langle
f_\alpha\colon\alpha<\lambda\rangle$, 
getting $\langle B'_\alpha\colon\alpha<\lambda\rangle$, now
$B'_\alpha$ for $\alpha$ large enough is as required.
\sqed{1.8}

\proclaim{Conclusion 1.9}
If $(*)$ of 1.5, then $\pcf(\bar\lambda)$ has a last element.
\endproclaim

\demo{Proof}
This is the minimal $\lambda$ such that $\kappa\in
J_{\le\lambda}[\bar\lambda]$.
[$\lambda$ exists, since $\lambda^* =: |\Prod\bar\lambda |\in
\{\lambda\colon
\kappa\in J_{\le\lambda}[\bar\lambda]\}\not=\emptyset]$ and by
1.8(2).
\sqed{1.9}

\proclaim{Claim 1.10}
Suppose (*) of 1.5 holds.
Assume for $j<\sigma$, $D_j$ is a
filter on $\kappa$ extending $\{\kappa\setminus A\colon A\in I^*\}$,
$E$ a filter on $\sigma$ and $D^*=\{ B\subseteq \kappa\colon\{j<\sigma
\colon B\in D_j\}\in E\}$ ( a filter on $\kappa$).
Let  $\mu_j=:\tcf (\Prod\bar\lambda, <_{D_j})$ be well defined for
$j<\sigma$, and assume further $\mu_j>\sigma+\theta$.

\noindent
Let
$$
\lambda=\tcf(\prod\bar\lambda, <_{D^*}), \mu
=\tcf(\prod\limits_{j<\sigma}\mu_j, <_E).
$$
$\underline{\rm Then}$ $\lambda=\mu$ (in particular, if one is well defined,
than so is the other).
\endproclaim

\demo{Proof}
Wlog $\sigma\ge\theta$ (otherwise we can add $\mu_j=:\mu_0$, $D_j=: D_0$
for $j\in\theta\setminus\sigma$,
and replace $\sigma$ by $\theta$ and $E$ by
$E'=\{ A\subseteq\theta\colon A\cap\sigma\in E\})$.
Let $\langle f^j_\alpha\colon\alpha<\mu_j\rangle$ be an
$<_{D_j}$-increasing cofinal sequence in $(\prod\bar\lambda, <_{D_j})$.

Now $\ell=0, 1$, for each $f\in\prod\bar\lambda$, define
$G_\ell (f)\in\prod\limits_{j<\sigma}\mu_j$ by
$G_\ell(f)(j)=\min\{\alpha<\mu_j\colon$ if $\ell=1$ then $f\le f^j_\alpha
\mod D_j$ $\underline{\rm and}$ if $\ell=0$ then: not
$f^j_\alpha\leq f\ \mod\ D_j\}$
(it is well defined for $f\in\Prod\bar\lambda$ by the choice of $\langle
f_\alpha^j\colon\alpha<\mu_j\rangle$).

\noindent
Note that for $f^1$, $f^2\in\Prod\bar\lambda$ and $\ell<2$ we have:
$$
\eqalign{
f^1 & \le f^2\mod D^*\Leftrightarrow B(f^1, f^2)=:\{i<\kappa\colon
f^1(i)\le f^2(i)\}\in D^* \cr
& \Leftrightarrow A(f^1, f^2) =:\{j<\sigma\colon B(f^1, f^2)\in
D_j\}\in E \cr
& \Leftrightarrow \hbox{ for some }A\in E, \hbox{ for every }i\in A \hbox{
we have } f^1\leq_{D_i} f^2\cr
& \Rightarrow \hbox{ for some } A\in E \hbox{ for every } i\in A
\hbox{ we have } \cr
& \qquad\qquad\qquad G_\ell(f^1)(i)\le G_\ell(f^2)(i) \cr
& \Leftrightarrow G_\ell(f^1)\le
G_\ell(f^2)\mod E. \cr}
$$
So
\item{$\otimes_1$} $G_\ell$ is a mapping from $(\prod\bar\lambda, \le_{D^*})$
into
$(\prod\limits_{j<\sigma}\mu_j, \le_E)$ preserving order.

Next we prove  that

\item{$\otimes_2$} for every $g\in\Prod_{j<\sigma}\mu_j$ for
some $f\in\Prod\bar\lambda$, we have $g \leq G_0(f)\mod E$.

\noindent
[Why? Note that $\min\{\mu_j\colon j<\sigma\}\ge\sigma^+\ge\theta^+$ and
$J_{\le\theta}[\bar\lambda]\subseteq
J_{\le\sigma}[\bar\lambda]$.
By 1.5 we know $(\Prod\bar\lambda, <_{J_{\le\sigma}[\bar\lambda]})$ is
$\sigma^+$-directed, hence for some $f\in\Prod\bar\lambda$:
\item{$(*)_1$} for $j<\sigma$ we have $f^j_{g(j)}<f\mod
J_{\le\sigma}[\bar\lambda]$.

\noindent
We here assumed $\sigma<\mu_j$, hence
$J_{\le\sigma}[\bar\lambda]\subseteq J_{<\mu_j}[\bar\lambda]$ (by
1.4(2)) but $J_{<\mu_j}[\bar \lambda]$ is disjoint to
$D_j$ by the definition of $J_{<\mu_j}[\bar \lambda]$ (by 1.8(2) + 1.3(13)(c))
so together with $(*)_1$:
\item{$(*)_2$} for $j<\sigma$, $f^j_{g(j)}<f\mod D_j$.

\noindent
So for every $j<\sigma$ we have $g(j)\leq G_0(f)(j)$ hence clearly
$g \leq G_0(f)$.] 
\item{$\otimes_3$} for $f\in \prod \bar \lambda$ we have $G_0(f) \leq G_1(f)$
[Why? read the definitions]

\item{$\otimes_4$} if $f_1$, $f_2 \in \prod \bar \lambda$ and $G_1(f_1)
<_E G_0 (f_2)$ then $f_1 <_{D^*} f_2$

[Why? as $G_1(f_1)<_E G_0(f_2)$ there is $B\in E$ such that: $j\in B
\Rightarrow G_1(f_1)(j) < G_0(f_2)(j)$ so for each $j\in A$ we have
$f_1\leq _{D_j} f^j_{G_1(f_1)(j)}$ by the definition of $G_1(f_1)$) and
$f^j_{G_1(f_1)(j)} <_{D_j} f_2$ (as $G_1(f_1)(j) < G_0(f_2)(j)$ and  the
definition of $G_0(f_2)(j)$) so together $f_1<_{D_j} f_2$. So $A(f_1,
f_2)=\{i<\kappa: f_1(i)< f_2(i)\}$ satisfies: $A(f_1, f_2)\in D_j$ for
every $j\in B$, hence $A(f_1, f_2)\in D^*$ (by the definition of $D^*$)
hence $f_1<_{D^*} f_2$ as required]

\noindent
Now first assume $\lambda= \tcf (\prod \bar \lambda, <_{D^*})$ is well
defined, so there is a sequence $\bar f=\langle f_\alpha:
\alpha<\lambda\rangle$ of members of $\prod \bar \lambda$,
$<_{D^*}$-increasing and cofinal. So $\langle G_0(f_\alpha):
\alpha<\lambda\rangle$ is $\leq_E$-increasing in $\prod\limits_{j<\sigma}
\mu_j$ (by $\otimes_1$), for every $g\in \prod\limits_{j<\sigma} \mu_j$ for
some $f\in \prod \bar \lambda$ we have $g\leq_E G_0(f)$ (why? by
$\otimes_2$), but by the choice of $\bar
f$ for some $\beta<\lambda$ we have $f <_{D^*} f_\beta$ hence by
$\otimes_1$ we have $g\leq_E
G_0(f) \leq_E G_0(f_\beta)$, so
$\langle G_0(f_\alpha): \alpha<\lambda\rangle$ is
cofinal in $(\prod\limits_{j<\sigma} \mu_j, <_E)$. Also for every
$\alpha<\lambda$, applying the previous sentence to $G(f_\alpha)+1$ ($\in
\prod\limits_{j<\sigma} \mu_j$) we can find $\beta<\lambda$ such that
$G(f_\alpha)+1 \leq_E G(f_\beta)$, so $G(f_\alpha)<_E G(f_\alpha)$, so for
some club $C$ of $\lambda$, $\langle G_0(f_\alpha): \alpha\in C\rangle$ is
$<_E$-increasing cofinal in $(\prod\limits_{j<\sigma} \mu_j, <_E)$. So
if $\lambda$ is well defined then $\mu=\tcf(\prod\limits_{j<\sigma}
\mu_j, <_E)$ is well defined and equall to $\lambda$.

Lastly assume that $\mu$ is well defined i.e. $\prod\limits_{j<\sigma}
\mu_j / E$ has true cofinality $\mu$, let $\bar g = \langle g_\alpha:
\alpha< \mu\rangle$ exemplifies it. Choose by induction on $\alpha<\mu$,
a function $f_\alpha$ and ordinals $\beta_\alpha$, $\gamma_\alpha$ such
that
\itemitem{(i)} $f_\alpha \in \prod \bar \lambda$

\itemitem{(ii)} $g_{\beta_\alpha} <_E G_0(f_\alpha) \leq_E
G_1(f_\alpha) <_E g_{\gamma_\alpha}$ (so $\beta_\alpha< \gamma_\alpha$)

\itemitem{(iii)} $\alpha_1 < \alpha_2 <\mu \Rightarrow
\gamma_{\alpha_1} < \beta_{\alpha_2}$ (so $\beta_\alpha \geq \alpha$)

\noindent
In stage $\alpha$, first choose $\beta_\alpha=\bigcup \{
\gamma_{\alpha_1}+1: \alpha_1 <\alpha\}$, then choose $f_\alpha \in \prod
\bar \lambda$ such that $g_{\beta_\alpha}+1 <_E G_0(f_\alpha)$ (possible
by $\otimes_2$) then choose $\gamma_\alpha$ such that $G_1(f_\alpha) <_E
g_{\gamma_\alpha}$. Now $G_0(f_\alpha) \leq_E G_1(f_\alpha)$ by
$\otimes_3$. By $\otimes_4$ we have $\alpha_1 < \alpha_2 \Rightarrow
f_{\alpha_1} <_{D^*} f_{\alpha_2}$. Also if $f\in \prod \bar \lambda$
then $G_1(f)\in \prod\limits_{j<\sigma} \mu_j$ hence by the choice of
$\bar g$, for some $\alpha< \mu$ we have $G_1(f) <_E g_\alpha$ but
$\alpha\leq \beta_\alpha$ so $G_1 (f)<_E g_\alpha \leq_E G_0(f_\alpha)$
hence by $\otimes_4$, $f<_{D^*} f_\alpha$. Altogether, $\langle
f_\alpha: \alpha< \mu\rangle$ exemplifies that $(\prod \bar \lambda,
<_{D^*})$ has true cofinality $\mu$, so $\lambda$ is well defined and
equal to $\mu$.
\sqed{1.11}

\proclaim{Conclusion 1.12}
If $(*)$ of 1.5 holds, and $\sigma$, $\bar\mu=\langle\mu_j\colon
j<\sigma\rangle$, $\langle D_j\colon j<\sigma\rangle$ are as in
1.10 and $\sigma+\theta<\min(\bar\mu)$, and $J$ is an ideal on $\sigma$
and $I$ an ideal on $\kappa$ such that $I^*\subseteq I\subseteq
\{A\subseteq \kappa:$ for some $B\in J$ for every $j\in \sigma\setminus
A$ we have $B\notin D_j\}$
(e.g. $I=I^*$)
$\underline{\rm then}$ $\pcf_J(\{\mu_j\colon j<\sigma\})\subseteq
\pcf_I(\bar\lambda)$.
\endproclaim

\demo{Proof} Let $E$ be an ultrafilter on $\sigma$ disjoint to $J$ then
we can define an ultrafilter $D^*$ on $\kappa$ as in 1.10, so clearly
$D^*$ is disjoint to $J$.
\sqed{1.11}


\head \S2 Normality of $\lambda\in\pcf(\bar\lambda)$ for
$\bar\lambda$ \endhead

Having found those ideals $J_{<\lambda}[\bar\lambda]$, we would like
to know more.
As $J_{<\lambda}[\bar\lambda]$ is increasing continuous in $\lambda$,
the question is how $J_<[\bar\lambda]$, $J_{<\lambda^+}[\bar\lambda]$
are related.

The simplest relation is
$J_{<\lambda^+}[\bar\lambda]=J_{<\lambda}[\bar\lambda]+ B$ for some
$B\subseteq\kappa$, and then we call $\lambda$ normal (for
$\bar\lambda$) and denote $ B= B_\lambda[\bar\lambda]$ though it is
unique only modulo $J_{<\lambda}[\bar\lambda]$.
We give a sufficient condition for exsitence of such $B$,
using this in 2.8; giving
the necessary definition in 2.3 and needed information in 2.4,
2.5, 2.6; lastly 2.7 is the essential uniqueness of cofinal sequences
in appropriate $\Prod\bar\lambda/I$.

\rem{Definition 2.1}
\item{(1)} We say $\lambda\in\pcf(\bar\lambda)$ is normal (for
$\bar\lambda$) if for some $ B\subseteq\kappa$,
$J_{\le\lambda}[\bar\lambda]=J_{<\lambda}[\bar\lambda]+ B$.

\item{(2)} We say $\lambda\in\pcf(\bar\lambda)$ is semi-normal (for
$\bar\lambda$) if there are $B_\alpha$ for $\alpha<\lambda$ such that:
\itemitem{(i)} $\alpha<\beta\Rightarrow B_\alpha\subseteq B_\beta\mod
J_{<\lambda}[\bar\lambda]$

\item{} and
\itemitem{(ii)} $J_{\le\lambda}[\bar\lambda]=J_{<\lambda}[\bar\lambda]+
\{ B_\alpha\colon \alpha<\lambda\}.$

\item{(3)} We say $\bar\lambda$ is normal if every
$\lambda\in\pcf(\bar\lambda)$ is normal for $\bar\lambda$.
Similarly for semi normal.

\item{(4)} In (1), (2), (3) instead $\bar\lambda$ we can say
$(\bar\lambda, I)$ or $\Prod\bar\lambda/I$ or
$(\Prod\bar\lambda, <_{I})$ if we replace $I^*$ by $I$ (an ideal
on $\Dom(\bar\lambda)$).

\rem{Fact 2.2}
Suppose (*) of 1.5 and $\lambda\in\pcf(\bar\lambda)$.
Now:
\item{(1)} $\lambda$ is semi-normal for $\bar\lambda$ iff for some
$F=\{f_\alpha\colon\alpha<\lambda\}\subseteq\Prod\bar\lambda$ we have:
$[\alpha<\beta\Rightarrow f_\alpha<f_\beta\mod
J_{<\lambda}[\bar\lambda]]$ and for every
ultrafilter $D$ over $\kappa$ disjoint to $J_{<\lambda}[\bar \lambda]$,
$F$ is unbounded
in $(\Prod\bar\lambda, <_D)$ whenever $\tcf(\Prod\bar\lambda, <_D)=
\lambda$.

\item{(2)} In 2.1(2), without loss of generality, we may assume
that 
\item{}$\quad\quad$ either: $B_\alpha= B_0\mod J_{<\lambda}[\bar\lambda]$
(so $\lambda$ is normal) 
\item{}$\quad\quad$ or: $ B_\alpha\not= B_\beta\mod
J_{\le\lambda}[\bar\lambda]$ for $\alpha<\beta<\lambda$.

\item{(3)} Assume $\lambda$ is semi normal for $\bar \lambda$.
Then $\lambda$ is normal for $\bar \lambda$ $\underline{\rm iff}$
for some $F$ as in part (1) (of 2.2) $F$ has
${\rm a}<_{J_{<\lambda}[\bar\lambda]}$-exact upper bound
$g\in\Prod_{i<\kappa}(\lambda_i+1)$ and then $B=:\{i<\kappa\colon
g(i)=\lambda_i\}$ generates $J_{\le\lambda}[\bar\lambda]$ over
$J_{<\lambda}[\bar\lambda]$.

\item{(4)} If $\lambda$ is semi regular for $\bar \lambda$ then for some
$\bar f=\langle f_\alpha: \alpha<\lambda\rangle$, $\bar B=\langle
B_\alpha: \alpha<\lambda\rangle$ we have: $\bar B$ is increasing modulo
$J_{<\lambda}[\bar \lambda]$, $J_{\leq \lambda}[\bar \lambda]= J_{<\lambda}[\Ga]+
\{B_\alpha: \alpha<\lambda\}$, and $\langle f_\alpha:
\alpha<\lambda\rangle$ is $<_{J_{<\lambda}[\Ga]}$-increasing and
cofinal, and $\bar f$, $\bar B$ as in 1.7.

\demo{Proof}
1) For the direction $\Rightarrow$, given $\langle B_\alpha: \alpha <
\lambda \rangle$ as in Definition 2.1(2), for each $\alpha<\lambda$, by
1.8(1) we have $(\prod \bar \lambda \restriction B_\alpha,
<_{J_{<\lambda} [\bar \lambda]})$ has true cofinality $\lambda$, and let
it be exemplified by $\langle f^\alpha_\beta: \beta< \lambda\rangle$. By
1.5 we can choose by induction on $\gamma< \lambda$ a function
$f_\gamma\in \prod \bar \lambda$ such that: $\beta$, $\gamma\leq \alpha
\Rightarrow f^\alpha_\beta \leq_{J_{<\lambda} [\bar \lambda]} f_\gamma$
and $\beta<\gamma \Rightarrow f_\beta<_{J_{<\lambda}[\bar \lambda]}
f_\gamma$.
\bigskip
\bigskip
\bigskip

Now $F=: \{f_\alpha: \alpha< \lambda\}$ is as required. [Why? First,
obviously $\alpha<\beta \Rightarrow f_\alpha< f_\beta\ \mod\
J_{<\lambda} [\bar \lambda]$. Second, if $D$ is an ultrafilter on
$\kappa$ disjoint to $I^*$ and $(\prod \bar \lambda, <_D)$ has true
cofinality $\lambda$, then by 1.6 for some $B \in J_{\leq \lambda} [\bar
\lambda] \setminus J_{<\lambda}[\bar \lambda]$ we have $B\in D$, so for
some $\alpha<\lambda$, $B\subseteq B_\alpha\ \mod\ J_{<\lambda} [\bar
\lambda]$ hence $B_\alpha\in D$. As $f^\alpha_\beta
\leq_{J_{<\lambda}[\bar \lambda]} f_\beta$ for $\beta\in [\alpha,
\lambda)$ clearly $F$ is cofinal in $(\prod \bar \lambda, <_D)$.]

The other direction, $\Leftarrow$ follows from 1.7 applied to
$F=\{f_\alpha: \alpha<\lambda\}$. [Why? we get there $\langle B_\alpha:
\alpha<\lambda\rangle$, $B_\alpha\in J_{\leq \lambda}[\bar \lambda]$
increasing modulo $J_{<\lambda}[\bar \lambda]$ so $J=: J_{<\lambda}[\bar
\lambda]+ \{B_\alpha: \alpha<\lambda\}\subseteq J_{\leq \lambda}[\bar
\lambda]$.

If equality does hold then for some ultrafilter $D$ over $\kappa$,
$D\cap J =\emptyset$ but $D\cap J_{\leq \lambda}[\bar \lambda]\not=
\emptyset$ so by clause (D) of 1.7, $F$ is bounded in $\prod \lambda/D$
whereas by 1.8(1),(2), $\tcf(\prod \bar \lambda, <_D)=\lambda$
contradicting the assumption on $F$.]

2) Because we can replace $\langle B_\alpha: \alpha< \lambda\rangle$ by
$\langle B_{\alpha_i}: i<\lambda\rangle$ whenever $\langle \alpha_i:
i<\lambda\rangle$ is non decreasing, non eventually constant.

3) If $\lambda$ is normal for $\bar \lambda$, let $B\subseteq \kappa$ be
such that $J_{\leq \lambda}[\bar \lambda]= J_{<\lambda}[\bar \lambda]
+B$. By 1.8(1) we know that $(\prod (\bar \lambda \restriction B),
<_{J_{<\lambda}[\bar \lambda]})$ has true cofinality $\lambda$, so let
it be exemplified by $\langle f^0_\alpha: \alpha<\lambda\rangle$. Let
$f_\alpha= f^0_\alpha \cup 0_{(\kappa\setminus B)}$ for
$\alpha<\lambda$. Now $\langle f_\alpha: \alpha<\lambda\rangle$ is as
required by 1.3(11).

Now suppose $\langle f_\alpha: \alpha<\lambda\rangle$ is as in part (1)
of 2.2 and $g$ is a $<_{J_{<\lambda}[\bar \lambda]}-\eub$ of $F$, $g\in
\prod\limits_{i<\kappa} (\lambda_i+1)$ and $B=\{i: g(i)=\lambda_i\}$.
Let $D$ be an ultrafilter on $\kappa$ disjoint to $J_{<\lambda}[\bar
\lambda]$. If $B\in D$ then for every $f\in \prod \bar \lambda$, let
$f'= (f\restriction B)\cup 0_{(\kappa \setminus B)}$, now necessarily
$f'<\max\{ g, 1\}$ (as $[i\in B \Rightarrow  f'(i)< \lambda_i=g(i)]$ and
$[i\in \kappa\setminus B \Rightarrow f'(i)=0\leq g(i)]$), hence (see
Definition 1.2(4)) for some $\alpha<\lambda$ we have $f'< \max\{
f_\alpha, 1\}\ \mod\ J_{<\lambda}[\bar \lambda]$ hence
for some $\alpha<\lambda$, $f'\leq f_\alpha\ \mod\
J_{<\lambda}[\bar \lambda]$ hence $f\leq f' \leq f_\alpha\ \mod\ D$; also
$\alpha<\beta \Rightarrow f_\alpha < f_\beta\ \mod\ D$, hence together
$\langle f_\alpha: \alpha<\lambda\rangle$ exemplifies $\tcf(\prod \bar
\lambda, <_D)=\lambda$. If $B\notin D$ then $\kappa\setminus B\in D$ so
$g'= g\restriction (\kappa \setminus B) \cup 0_B= g\ \mod\ D$ and
$\alpha< \lambda \Rightarrow f_\alpha <_D f_{\alpha+1} \leq_D g=_D g'$,
so $g'\in \prod \bar \lambda$ exemplifies $F$ is bounded in $(\prod \bar
\lambda, <_D)$ so as $F$ is as in 2.2(1), $\tcf(\prod \bar \lambda,
<_D)=\lambda$ is impossible. As $D$ is disjoint to $J_{<\lambda}[\bar
\lambda]$, necessarily $\tcf(\prod \bar \lambda, <_D)> \lambda$. The
last two arguments together give, by 1.8(2) that $J_{\leq \lambda}[\bar
\lambda] = J_{<\lambda}[\bar \lambda] + B$ as required in the definition
of normality.

4) Should be clear.
\sqed{2.2}

We shall give some sufficient conditions for normality.

\rem{Remark}
In the following definitions we slightly diviate from [Sh-g, Ch I
=Sh345a].
The ones here are perheps somewhat artificial but enable us to
deal also with case $(\beta)$ of 1.5(*). I.e. in Definition 2.3 below we
concentrate on the first $\theta$ elements of an $a_\alpha$ and for
``obey" we also have $\bar A^*=\langle A_\alpha: \alpha<\theta\rangle$
and we want to cover also the case $\theta$ is singular.

\rem{Definition 2.3}
Let there be given regular $\lambda$, $\theta<\mu<\lambda$,
$\mu$ possibly an ordinal, $S\subseteq\lambda$, $\sup (S)=\lambda$ and
for simplicity $S$ is a set of limit ordinals or at least have no two
successive members.
\item{(1)} We call $\bar a=\langle
a_\alpha\colon\alpha<\lambda\rangle$ a continuity condition for
$(S, \mu, \theta)$ (or is an $(S, \mu, \theta)$-continuity condition) if: $S$ is an unbounded subset of $\lambda$,
$a_\alpha\subseteq\alpha$,
$\otp(a_\alpha)<\mu$,
and $[\beta\in a_\alpha\Rightarrow a_\beta=a_\alpha\cap\beta]$
 and,
for every club $E$ of $\lambda$,
for some\footnote\dag{Note: if $\otp(a_\delta)=\theta$ and
$\delta=\sup(a_\delta)$ (holds if $\delta\in S$, $\mu=\theta+1$ and
$\bar a$ continuous in $S$ (see below)) $\underline{\rm then}$
$\delta\in E$.}
$\delta\in S$ we have $\theta=\otp
\{\alpha\in a_\delta\colon\otp(a_\alpha)<\theta$ and for no
$\beta\in a_\delta\cap\alpha$ is $(\beta, \alpha)\cap E=\emptyset\}$. We
say $\bar a$ is continuous in $S^*$ if $\alpha\in S^* \Rightarrow
\alpha=\sup (a_\alpha)$.

\item{(2)}
Assume $f_\alpha\trans{\in}{\Ord}$ for $\alpha<\lambda$ and $\bar
A^*=\langle A^*_\alpha\colon\alpha<\theta\rangle$ be a
decreasing sequence of subsetes of $\kappa$ such that
$\kappa\setminus A^*_\alpha\in I^*$.
We say $\bar f=\langle
f_\alpha\colon\alpha<\lambda\rangle$ obeys $\bar a=\langle
a_\alpha\colon\alpha<\lambda\rangle$ for $\bar A^*$ if:
\itemitem{(i)} for $\beta\in a_\alpha$,
if $\varepsilon=:\otp(a_\alpha)<\theta$ then we have $f_\beta\upharpoonright
A^*_\varepsilon\le f_\alpha\upharpoonright A^*_\varepsilon$
(note: $\bar A^*$ determine $\theta$).

\item{(2A)} Let $\kappa$, $\bar\lambda$, $I^*$ be as usual.
We say $\bar f$ obeys $\bar a$ for $\bar A^*$
continuously on $S^*$ if: $\bar a$ is continuous in $S^*$ and
$\bar f$ obeys $\bar a$ for $\bar A^*$
and in addition $S^*\subseteq S$ and for
$\alpha\in S^*$ (a limit ordinal) we have
$f_\alpha=f_{a_\alpha}$ from (2B), i.e.\
for every $i<\kappa$ we have
$f_\alpha(i)=\sup\{ f_\beta(i)\colon\beta\in a_\alpha\}$ when
$|a_\alpha|<\lambda_i$.

\item{(2B)} For given $\bar\lambda=\langle\lambda_i\colon
i<\kappa\rangle$, $\bar f=\langle
f_\alpha\colon\alpha<\lambda\rangle$
 where $f_\alpha\in\prod\bar\lambda$ and $a\subseteq\lambda$, and $\theta$
let $f_a\in\prod\bar\lambda$ be defined by:
$f_a(i)$ is
$0$ if $|a|\ge\lambda_i$ and $\cup\{ f_\alpha(i)\colon\alpha\in a\}$
if $|a|<\lambda_i$.

\item{(3)}
Let ($S, \theta$) stands for $(S, \theta+1, \theta)$; 
$(\lambda, \mu, \theta)$ stands for
``$(S, \mu, \theta)$ for some unbounded subset $S$ of $\lambda$'' and
$(\lambda, \theta)$ stands for $(\lambda, \theta+1, \theta)$.

If each $A^*_\alpha$ is $\kappa$ then we omit ``for $\bar
A^*$'' (but $\theta$ should be fixed or said).

\item{(4)} We add to ``continuity condition'' (in part (1)) the
adjective ``weak'' [``$\theta$-weak''] if
``$\beta\in a_\alpha\Rightarrow a_\beta=
a_\alpha\cap\beta$'' is replaced by ``$\alpha\in S\&\beta\in
a_\alpha\Rightarrow (\exists\gamma
<\alpha)[a_\alpha\cap\beta\subseteq a_\gamma\ \&\ \gamma<\min(a_\alpha
\setminus (\beta+1))\ \&\ [|a_\alpha \cap \beta|<\theta \Rightarrow
|a_\gamma \cap \beta| <\theta]]$'' [and we demand that $\gamma$ exists
only if  $\otp(a_\alpha\cap\beta)<\theta$].
(Of course a continuity condition is a weak continuity
condition which is a $\theta$-weak continuity condition).

\rem{Remark 2.3A} There are some obvious monotonicity implications, we
state below only 2.4(3).

\rem{Fact 2.4}
\item{(1)} Let $\theta_r=\cases{\theta & $\cf(\theta)=\theta$\cr
\theta^+ & $\cf(\theta)<\theta$\cr}$ and assume
$\lambda=\cf(\lambda)>\theta_r^+$.
$\underline{\rm Then}$ for some stationary
$S\subseteq\{\delta<\lambda\colon\cf(\delta)=
\theta_r\}$,
there is a continuity condition $\bar a$ for $(S, \theta_r)$;
moreover, it is continuous in $S$ and $\delta\in S\Rightarrow\otp
(a_\delta)=\theta_r$; so for 
every club $E$ of $\lambda$ for some $\delta\in S$,
$\forall\alpha,\beta [\alpha<\beta\ \&\ \alpha\in a_\delta\ \&\ \beta\in
a_\delta\to(\alpha, \beta)\cap E\not=\emptyset\}]$.

\item{(2)} Assume $\lambda=\theta^{++}$, then for some stationary
$S\subseteq \{\delta<\lambda: \cf(\delta)=\cf(\theta)\}$ there is a
continuity condition for $(S, \theta+1, \theta)$.

\item{(3)} If $\bar a$ is a $(\lambda, \mu, \theta_1)$-continuity
condition and $\theta_1\geq \theta$ $\underline{\rm then}$ there is a
$(\lambda, \theta+1, \theta)$-continuity condition.

\demo{Proof}
1) By [Sh420, \S1].

2) By [Sh351, 4.4(2)] and\footnote\dag{the definition of $B^\alpha_i$ in
the proof of [Sh-g, III 2.14(2)] should be changed as in [Sh351,
4.4(2)]} [Sh-g, III 2.14(2), clause (c), p.135-7].

3) Check.
\sqed{2.4}

\rem{Remark 2.4A} Of course also if $\lambda=\theta^+$ the conclusion of
2.4(2) may well hold. We suspect but do not know that the negation is
consistent with ZFC.

\rem{Fact 2.5}
Suppose (*) of 1.5,
$f_\alpha\in\Prod\bar\lambda$ for $\alpha<\lambda$,
$\lambda=\cf(\lambda)\leq \theta$ (of course
$\kappa=\dom(\bar\lambda))$ and $\bar A^*=\bar A^*[\bar\lambda]$ is
as in the proof of 1.5 (i.e.\ $A^*_\alpha=\{
i<\kappa\colon\lambda_i>\alpha\}$).
$\underline{\rm Then}$
\item{(1)}
Assume $\bar a$ is a $\theta$-weak continuity condition for $(S,
\theta)$, $\lambda=\sup (S)$, $\underline{\rm then}$
we can find $\bar f'=\langle f'_\alpha\colon
\alpha<\lambda\rangle$ such that:
\itemitem{(i)}$f'_\alpha\in\Prod\bar\lambda$,

\itemitem{(ii)} for $\alpha<\lambda$ we have
$f_\alpha\le f'_\alpha$

\itemitem{(iii)} for $\alpha<\beta<\lambda$ we have $f'_\alpha
<_{J_{<\lambda}[\bar\lambda]}f'_\beta$

\itemitem{(iv)} $\bar f'$ obeys $\bar a$ for $\bar A^*$

\item{(2)} If in addition $\min(\bar\lambda)>\mu$,
$S^*\subseteq S$ are stationary subsets of $\lambda$ and $\bar a$ is
a continuity condition for $(S, \mu, \theta)$
$\underline{\rm then}$ we can find $\bar f'=\langle
f'_\alpha\colon\alpha<\lambda\rangle$ such that
\itemitem{(i)} $f'_\alpha\in\prod\bar\lambda$

\itemitem{(ii)} for $\alpha\in\lambda\setminus S^*$ we have
$f_\alpha\le f'_\alpha$ and $\alpha=\beta+1\in \lambda\setminus S^*\ \&\
\beta\in S^* \Rightarrow f_\beta \leq f'_\alpha$

\itemitem{(iii)} for $\alpha<\beta<\lambda$ we have
$f'_\alpha<_{J_{<\lambda}[\bar\lambda]}f'_\beta$

\itemitem{(iv)} $\bar f'$ obeys $\bar a$ for $\bar A^*$
continuously on $S^*$

\item{(3)} Suppose $\langle
f'_\alpha\colon\alpha<\lambda\rangle$ obeys $\bar a$ continuously
on $S^*$ and satisfies 2.5(2)(ii) (and 2.5(2)'s assumption holds).
If $g_\alpha\in\Prod\bar\lambda$ and $\langle
g_\alpha\colon\alpha<\lambda\rangle$ obeys $\bar a$ continuously on
$S^*$ and $[\alpha\in\lambda\setminus S^*\Rightarrow
g_\alpha\le f_\alpha]$ $\underline{\rm then}$
$\bigwedge_\alpha g_\alpha\le f'_\alpha$.

\item{(4)} If $\zeta<\theta$, for $\varepsilon<\zeta$ we have $\bar
f^\varepsilon=\langle
f^\varepsilon_\alpha\colon\alpha<\lambda\rangle$,
where $f^\varepsilon_\alpha\in\prod\bar\lambda$, $\underline{\rm then}$
in 2.5(1) (and
2.5(2)) we can find $f'$ as there for all $\bar f^\varepsilon$
simultaneously. Only in clause (ii) we replace $f_\alpha\leq f'_\alpha$
by $f_\alpha \restriction A^*_\zeta \leq f'_\alpha \restriction
A^*_\zeta$ (and $f_\beta\leq f'_\alpha$ by $f_\beta \restriction
A^*_\zeta \leq f'_\alpha \restriction A^*_\zeta$.

\demo{Proof}
Easy (using 1.5 of course).

\proclaim{Claim 2.5A}
In 2.5 we can replace ``(*) from 1.5'' by ``$\prod\bar\lambda/J_{<\lambda}
[\bar\lambda]$ is $\lambda$-directed''.
\endproclaim

\proclaim{Claim 2.6}
Assume (*) of 1.5 and let $\bar A^*$ be as there,
\item{(1)} in 1.7, if $\langle
f_\alpha\colon \alpha<\lambda\rangle$ obeys some
$(S, \theta)$-continuity condition or just a $\theta$-weak one
for $\bar A^*$ (where $S\subseteq \lambda$ is unbounded) 
$\underline{\rm then}$ we can deduce also:\newline
$(G)$ the sequence $\langle B_\alpha/J_{<\lambda}[\bar\lambda]\colon
\alpha<\lambda\rangle$ is eventually constant.

\item{(2)} If $\theta^+ <\lambda$ $\underline{\rm then}$
$J_{\le\lambda}[\bar\lambda]/J_{<\lambda}[\bar\lambda]$ is
$\lambda^+$-directed (hence if $\lambda$ is semi normal for $\bar
\lambda$ then it is normal to $\bar \lambda$).

\endproclaim

\demo{Proof}
1) Assume not, so for some club $E$ of $\lambda$ we have
\item{$(*)$} $\alpha<\delta<\lambda\&\delta\in E\Rightarrow
B_\alpha\not= B_\delta\mod J_{<\lambda}[\bar\lambda]$.

As $\bar a$ is a $\theta$-weak $(S, \theta)$-continuity condition, there is
$\delta\in S$ such that $b=:\{\alpha\in
a_\delta\colon\otp(a_\delta\cap\alpha)<\theta$ and for no
$\beta\in a_\delta\cap\alpha$ is $(\beta,\alpha)\cap
E=\emptyset\}$ has order type $\theta$.
Let $\{\alpha_\varepsilon\colon\varepsilon<\theta\}$ list $b$
(increasing with $\varepsilon$). 
So for every $\varepsilon<\theta$ there is
$\gamma_\varepsilon\in(\alpha_\varepsilon,\alpha_{\varepsilon+1})\cap
E$, and let $\beta_\varepsilon<\alpha_{\varepsilon+1}$
be such that $a_\delta\cap\alpha_\varepsilon\subseteq
a_{\beta_\varepsilon}$ and $\otp(a_{\beta_\varepsilon}\cap
\alpha_\varepsilon)<\theta$; by shrinking and renaming
$\wlog\beta_\varepsilon<\gamma_\varepsilon$
and $\alpha_\varepsilon\in a_{\beta_\varepsilon}$. Let $\xi(\varepsilon)
=: \otp(a_{\beta_\varepsilon}\cap \alpha_\varepsilon)$.
Lastly let $B^0_\varepsilon =: \{ i<\kappa\colon
f_{\alpha_\varepsilon}(i)< f_{\beta_\varepsilon}(i)<
f_{\gamma_\varepsilon}(i)<f_{\alpha_{\varepsilon+1}}(i)\}$,
clearly it is $=\kappa\mod I^*$ and let (remember $(*)$ above)
$B_\varepsilon^*=: A^*_{\xi(\varepsilon)+1}\cap(B_{\gamma_\varepsilon}
\setminus B_{\beta_\varepsilon})\cap B_\varepsilon^0$, now
$B_{\alpha_\varepsilon}\subseteq
B_{\beta_\varepsilon}\subseteq B_{\gamma_{\varepsilon}}\mod
J_{<\lambda}[\bar\lambda]$ by clause (B) of 1.7, and
$B_{\gamma_\varepsilon}\not= B_{\beta_{\varepsilon}}$ by $(*)$
above hence $B_{\gamma_{\varepsilon}}\setminus
B_{\beta_\varepsilon}\not=\emptyset\mod
J_{<\lambda}[\bar\lambda]$. Now $B_\varepsilon^0$,
$A^*_{{\xi(\varepsilon)+1}}=\kappa\mod I^*$ by the previous sentence and
by 1.5$(*)$ which we
are assuming respectively and $I^*\subseteq J_{<\lambda}[\bar\lambda]$ by the
later's definition; so we have gotten
$B^*_\varepsilon\not=\emptyset\mod J_{<\lambda}[\bar\lambda]$.
But for $\varepsilon<\zeta<\theta$ we have $B^*_\varepsilon\cap
B^*_\zeta=\emptyset$, for suppose $i\in B^*_\varepsilon\cap
B^*_\zeta$, so $i\in A^*_{{\xi(\varepsilon)+1}}$ and also 
$f_{\gamma_\varepsilon} (i)< f_{\alpha_{\varepsilon+1}}(i)\le
f_{\beta_\zeta}(i)$
(as $i\in B_\varepsilon^0$ and as $\alpha_{\varepsilon+1}\in
a_{\beta_\zeta}\ \&\ i\in A^*_{\xi(\zeta)+1}$ respectively);
now $i\in B^*_\varepsilon$ hence $i\in B_{\gamma_{\varepsilon}}$
i.e.\ (where $g$ is from 1.7 clause (D)$^+$)
$f_{\gamma_\varepsilon}(i)>g(i)$ hence (by the above)
$f_{\beta_\zeta}(i)>g(i)$
hence $i\in B_{\beta_\zeta}$ hence $i\not\in B^*_\zeta$,
contradiction.
So $\langle B^*_\varepsilon\colon\varepsilon<\theta\rangle$ is a
sequence of $\theta$ pairwise disjoint members of
$(J_{<\lambda}[\bar\lambda])^+$, contradiction.

2) The proof is similar to the proof of 1.8(4), using 2.6(1) instead 1.7
(and $\bar a$ from 2.4(1) if $\lambda> \theta^+_r$ or 2.4(2) if $ \lambda=
\theta^{++}$). \sqed{2.6}

We note also (but shall not use):

\proclaim{Claim 2.7}
Suppose $(*)$ of 1.5 and
\item{(a)} $f_\alpha\in\Prod\bar\lambda$ for $\alpha<\lambda$,
$\lambda\in\pcf(\bar\lambda)$ and $\bar f=\langle
f_\alpha\colon\alpha<\lambda\rangle$ is
$<_{J_{<\lambda}[\bar\lambda]}$-increasing

\item{(b)} $\bar f$ obeys $\bar a$ continuously on $S^*$, where $\bar a$
is a continuity condition for $(S, \theta)$ and $\lambda=\sup(S)$ (hence
$\lambda>\theta$ by the last phrase of 2.3(1))

\item{(c)} $J$ is an ideal on $\kappa$ extending $J_{<\lambda}[\bar\lambda]$,
and $\langle f_\alpha/J\colon\alpha<\lambda\rangle$ is cofinal in
$(\Prod\bar\lambda,<_J)$ (e.g. $J=J_{<\lambda}[\bar\lambda]+(\kappa\setminus
B)$, $ B\in J_{\le\lambda}[\bar\lambda]\setminus J_{<\lambda}[\bar\lambda]$).

\item{(d)} $\langle f'_\alpha\colon\alpha<\lambda\rangle$
satisfies (a), (b) above.

\item{(e)} $f_\alpha\le f'_\alpha$ for $\alpha\in\lambda\setminus
S^*$
(alternatively: $\langle f'_\alpha\colon\alpha<\lambda\rangle$
satisfies (c)).

\item{(f)} if $\delta\in S^*$ $\underline{\rm then}$ $J$ is
$\cf(\delta)$-indecomposable (i.e.\ if $\langle
A_\varepsilon\colon\varepsilon<\cf(\delta)\rangle$ is a
$\subseteq$-increasing sequence of members, of $J$ then
$\bigcup_{\varepsilon<\cf(\delta)} A_\varepsilon\in J$).
\endproclaim

$\underline{\rm Then}$:
\item{(A)} the set
$$\{\delta<\lambda: \hbox{ if }\delta\in S^*\hbox{ and }\otp(a_\delta)=\theta
\hbox{ then }f'_\delta=f_\delta\mod J\}$$
contains a club of $\lambda$.
\item{(B)} the set
$$\eqalign{\{\delta<\lambda\colon &\hbox{if }\alpha\in S\hbox{ and }\delta=
\sup(\delta\cap a_\alpha)\hbox{ and }\otp(\alpha\cap a_\delta)=\theta\cr
\ &\hbox{then }f'_{\alpha\cap a_\delta}=f_{\alpha\cap a_\delta}\mod
J\}\cr}$$
contains a club of $\lambda$.

\demo{Proof}
We concentrate on proving (A).
Suppose $\delta\in S^*$, and $f_\delta\not= f'_\delta\mod J$.
Let
$$
\eqalign{
A_{1,\delta} & =\{i<\kappa\colon f_\delta(i)<f'_\delta(i)\}\cr
A_{2,\delta} & =\{i<\kappa\colon f_\delta(i)>f'_\delta(i)\},\cr}
$$
So $A_{1,\delta}\cup A_{2,\delta}\in J^+$,
suppose first $A_{1,\delta}\in J^+$.
By Definition 2.3(2A), for every $i\in A_{1,\delta}$ for every
large enough $\alpha\in a_\delta$, $f_\delta(i)<f'_\alpha(i)$,
say for $\alpha\in a_\delta\setminus\beta_i$.
As $J$ is $\cf(\delta)$-indecomposable for some $\beta<\alpha$ we have
$\{i<\kappa\colon\beta_i<\beta\}\in J^+$ so $f_\delta \restriction A_{1,
\delta}< f'_\beta \restriction A_{1, \delta}$ (and $\beta<\delta$). Now
by clause (c),
$E =:\{\delta<\lambda\colon$ for every $\beta<\delta$ we have
$f'_\beta<f_\delta\mod J\}$ is a club of $\lambda$, and so we have
proved
$$
\delta\in E\Rightarrow A_{1,\delta}\in J.
$$

If $\bigwedge_{\alpha< \lambda} f_\alpha\le f'_\alpha$
(first possibility in clause (e)) also $A_{2,\delta}\in J$ hence
for no $\delta\in S^*\cap
E$ do we have $f_\delta\not= f'_\delta\mod J$.
If the second possibility of clause (e) holds, we can
interchange $\bar f$, $\bar f'$ hence $[\delta\in E\Rightarrow
A_{2,\delta}\in J]$ and we are done.
\sqed{2.7}

We now return to investigating the $J_{<\lambda}[\bar \lambda]$, first
without using continuity conditions.

\proclaim{Lemma 2.8}
Suppose $(*)$ of 1.5 and
$\lambda=\cf(\lambda)\in\pcf(\bar\lambda)$;
$\underline{\rm Then}$ $\lambda$ is semi normal for
$\bar\lambda$.
\endproclaim

\demo{Proof}
We assume $\lambda$ is not semi normal for $\bar\lambda$ and
eventually get a contradiction.
Note that by our assumption $(\prod \bar \lambda, <_I)$ is
$\theta^+$-directed hence $\min\pcf_I(\bar \lambda)\geq \theta^+$ (by
1.3(4)(v)) hence
let us define by induction on $\xi\le\theta$, $\bar
f^\xi=\langle f_\alpha^\xi\colon\alpha<\lambda\rangle$, $B_\xi$
and $D_\xi$ such that:

\item{(I)}
\itemitem{(i)} $f^\xi_\alpha\in\Prod\bar\lambda$

\itemitem{(ii)} $\alpha<\beta<\lambda\Rightarrow f^\xi_\alpha\le
f^\xi_\beta\mod J_{<\lambda}[\bar\lambda]$

\itemitem{(iii)}
$\alpha<\lambda\& \xi<\theta\Rightarrow
f_\alpha^\xi\le f_\alpha^\theta\mod J_{<\lambda}[\bar\lambda]$

\itemitem{(iv)} for $\zeta<\xi<\theta$ and
$\alpha<\lambda\colon f^\zeta_\alpha\upharpoonright A^*_\xi\le
f^\xi_\alpha\upharpoonright A^*_\xi$

\item{(II)}
\itemitem{(i)} $D_\xi$ is an ultrafilter on $\kappa$ such that:
$\cf(\Prod\bar\lambda /D_\xi)=\lambda$

\itemitem{(ii)} $\langle f^\xi_\alpha
/D_\xi\colon\alpha<\lambda\rangle$ is not cofinal in $\Prod\bar\lambda
/D_\xi$

\itemitem{(iii)} $\langle f^{\xi+1}_\alpha
/D_\xi\colon\alpha<\lambda\rangle$ is increasing and cofinal in
$\Prod\bar\lambda /D_\xi$; moreover

\itemitem{(iii)$^+$} $B_\xi\in D_\xi$ and $\langle
f_\alpha^{\xi+1}\colon\alpha<\lambda\rangle$ is increasing and
cofinal in $\prod\bar\lambda/(J_{<\lambda}[\bar\lambda]+(\kappa\setminus
B_\xi))$

\itemitem{(iv)} $f_0^{\xi+1} /D_\xi$ is above $\{ f_\alpha^\xi
/D_\xi\colon\alpha<\lambda\}$.

\subdemo{For $\xi=0$}
no problem.
[Use 1.8(1)+(4)].

\subdemo{For $\xi$ limit $<\theta$}
Let $g_\alpha^\xi\in\Prod\bar\lambda$ be defined by
$g_\alpha^\xi(i)=\sup\{ f_\alpha^{\zeta}(i)\colon \zeta<\xi\}$
for $i\in A^*_\xi$ and $f_\alpha^\xi(i)=0$ else,
(remember that $\kappa\setminus A^*_\xi \in I^*$).
Then choose by induction on $\alpha<\lambda$, $f^\xi_\alpha\in \prod
\bar \lambda$ such that $g^\xi_\alpha\leq f^\xi_\alpha$ and
$\beta<\alpha \Rightarrow f_\beta< f_\alpha\ \mod\ J_{<\lambda}[\bar
\lambda]$. This is possible by 1.5 and clearly the requirements
(I)(i),(ii),(iv) are satisfied. 
Use 2.2(1) to find an appropriate $D_\xi$ (i.e.\ satisfying
II(i)+(ii)).
Now $\langle f_\alpha^\xi\colon\alpha<\lambda\rangle$ and
$D_\xi$ are as required.

\subdemo{For $\xi=\theta$}
Choose $f_\alpha^\theta$ by induction of $\alpha$ satisfying
I(i), (ii), (iii) (possible by 1.5).

\subdemo{For $\xi=\zeta+1$}
Use 1.6 to choose $B_\zeta\in D_\zeta\cap J_{\leq\lambda}[\bar
\lambda]\setminus J_{<\lambda}[\bar \lambda]$.
Let $\langle g_\alpha^\xi\colon\alpha<\lambda\rangle$ be cofinal
in $(\Prod\bar\lambda, <_{D_\xi})$ and even in $(\prod\bar\lambda,
<_{J_<[\bar\lambda]+(\kappa\setminus B_\xi)})$ and without loss of
generality $\bigwedge_{\alpha<\lambda} f_\alpha^\zeta / D_\zeta <
g_0^\xi / D_\zeta$ and $\bigwedge_{\alpha<\lambda} f_\alpha^\zeta
\upharpoonright {A^*_\xi}\le
g_\alpha^\xi\upharpoonright {A^*_\xi}$.
We get $\langle f_\alpha^\xi\colon\alpha<\lambda\rangle$
increasing and cofinal $\mod (J_{<\lambda}[\bar\lambda] +
(\kappa\setminus B_\xi))$ such that
$g^\xi_\alpha\leq f^\xi_\alpha$
by 1.5 from
$\langle g_\alpha^\xi\colon\alpha<\lambda\rangle$.
Then get $D_\xi$ as in the case ``$\xi$ limit''.

So we have defined the $f_\alpha^\xi$'s and $D_\xi$'s.
Now for each $\xi<\theta$ we apply (II) (iii)$^+$ for $\langle
f_\alpha^{\xi+1}\colon\alpha<\lambda\rangle$, $\langle
f_\alpha^\theta\colon\alpha<\lambda\rangle$.
We get a club $C_\xi$ of $\lambda$ such that:
$$
\alpha<\beta\in C_\xi \Rightarrow f^\theta_\alpha\upharpoonright
B_\xi<f_\beta^{\xi+1}\upharpoonright B_\xi\mod
J_{<\lambda}[\bar\lambda]
\leqno(*)
$$
So $C=:\bigcap_{\xi<\theta} C_\xi$ is a club of $\lambda$.
By 2.2(1) applied to $\langle
f_\alpha^\theta\colon\alpha<\lambda\rangle$ (and the assumption
``$\lambda$ is not semi-normal for $\bar\lambda$'') there is
$g\in\prod\bar\lambda$ such that
$$
\neg g\le f_\alpha^\theta \mod J_{<\lambda}[\bar\lambda]
\hbox{ for } \alpha<\lambda
\leqno{(*)_1}
$$
by 1.5 $\wlog$
$$
f^\xi_0<g\mod J_{<\lambda}[\bar\lambda]\quad
\hbox{ for } \xi<\theta
\leqno{(*)_2}
$$
For each $\xi<\theta$, by II (iii), (iii)$^+$ for some
$\alpha_\xi<\lambda$ we have
$$
g\upharpoonright B_\xi< f^{\xi+1}_{\alpha_\xi}\upharpoonright
B_\xi\mod J_{<\lambda}
[\bar\lambda]
\leqno{(*)_3}
$$
Let $\alpha(*)=\sup_{\xi<\theta}\alpha_\xi$, so
$\alpha(*)<\lambda$ and so
$$
g\upharpoonright B_\xi< f_{\alpha(*)}^{\xi+1} \upharpoonright
B_\xi \mod J_{<\lambda}[\bar\lambda]
\leqno{(*)_4}
$$
For $\zeta<\theta$, let $B_\zeta^*=\{i\in A^*_\zeta\colon g(i)<
f^\zeta_{\alpha(*)}(i)\}$.
By $(*)_4$, $B^*_{\xi+1}\in D_\xi$; by (II)(iv)$ + (*)_2$
we know $B^*_\xi\not\in D_\xi$, hence $B^*_\xi\not=B^*_{\xi+1}\mod
D_\xi$ hence $B_\xi^*\not= B_{\xi+1}^*\mod
J_{<\lambda}[\bar\lambda]$.

\noindent
On the other hand by (I)(iv) for each $\zeta< \theta$ we have  $\langle
B^*_\xi\cap A_\zeta^*\colon\xi\le\zeta\rangle$ is
$\subseteq$-increasing and (as $A^*_\zeta=\kappa\ \mod\
J_{<\lambda}[\bar \lambda]$ for each $\zeta< \theta$) hence by I(iv)
we have $\langle B^*_\xi/I^*: \xi<\theta\rangle$ is $\subseteq$-increasing,
and by the previous sentence $B^*_\xi\not= B^*_{\xi+1}\ \mod\
J_{<\lambda}[\bar \lambda]$ hence
$\langle B^*_\xi /I^*\colon
\xi<\theta\rangle$ is strictly $\subseteq$-increasing.
Together clearly $\langle B^*_{\xi+1}\cap A^*_{\xi+1}\setminus
B_\xi^*\colon \xi<\theta\rangle$ is a sequence of $\theta$
pairwise disjoint members of $(J_{<\lambda}[\bar\lambda])^+$, hence
of $(I^*)^+$;
contradiction to $\theta\ge\wsat (I^*)$.
\sqed{2.8}

\rem{Definition 2.9}
\item{(1)} We say
$\langle B_\lambda\colon\lambda\in\frc\rangle$ is a
$\underline{\rm generating~sequence}$ for $\bar\lambda$ if:
\itemitem{(i)} $ B_\lambda\subseteq\kappa$ and
$\frc\subseteq\pcf(\bar\lambda)$

\itemitem{(ii)} $J_{\le\lambda}[\bar\lambda]=J_{<\lambda}[\bar\lambda]+
 B_\lambda$  for each $\lambda\in \frc$

\item{(2)} We call $\bar B=\langle
 B_\lambda\colon\lambda\in\frc\rangle$ $\underline{\rm smooth}$ if:
$$
i\in B_\lambda\&\lambda_i\in \frc\Rightarrow  B_{\lambda_i}\subseteq
B_\lambda.
$$

\item{(3)} We call
$\bar B=\langle B_\lambda\colon\lambda\in\Rang(\bar\lambda)\rangle$
$\underline{\rm closed}$ if for each $\lambda$
$$
 B_\lambda\supseteq\{ i<\kappa\colon\lambda_i\in
\pcf(\bar\lambda\upharpoonright B_\lambda)\}
$$

\rem{Fact 2.10}
Assume $(*)$ of 1.5.

\item{(1)} Suppose $\frc\subseteq\pcf(\bar\lambda)$,
$\bar B=\langle B_\lambda\colon\lambda\in\frc\rangle$ is
a generating sequence for $\bar\lambda$, and $ B\subseteq\kappa$,
$\pcf(\bar\lambda\upharpoonright B)\subseteq\frc$
$\underline{\rm then}$
for some
finite $\frd\subseteq\frc$,
$ B\subseteq\bigcup_{\mu\in\frd} B_\mu\mod I^*$.

\item{(2)} $\cf(\Prod\bar\lambda/ I^*)=\max\pcf(\bar\lambda)$

\rem{Remark 2.10A}
For another proof of 2.10(2) see 2.12(2)+ 2.12(4) and for another use
of the proof of 2.10(2) see 2.14(1).

\demo{Proof}
(1) If not, then
$I=I^*+\{ B\cap\bigcup_{\mu\in\frd} B_\mu\colon\frd\subseteq
\frc$, $\frd$ finite$\}$ is a family of subsets of $\kappa$,
closed under union, $ B\not\in I$, hence there is an
ultrafilter $D$ on $\kappa$ disjoint from I to which $B$ belongs.
Let $\mu=:\cf(\Prod_{i<\kappa}\lambda_i/D)$;
necessarily $\mu\in\pcf(\bar\lambda\upharpoonright B)$, hence
by the last assumption of 2.10(1) we have $\mu\in\frc$.
By 1.8(2) we know $ B_\mu\in D$ hence $ B\cap B_\mu\in
D$, contradicting the choice of $D$.

(2)
The case $\theta=\aleph_0$ is trivial (as
$\wsat(I^*)\leq \aleph_0$ implies $\calP(\kappa)/I^*$ is a Boolean
algebra satisfying the $\aleph_0$-c.c. (as here we can substract) hence
this Boolean algebra is finite hence also
$\pcf(\bar\lambda)$ is
finite) so we assume $\theta>\aleph_0$.
For $B\in(I^*)^+$ let
$\lambda(B)=\max\pcf_{I^*\upharpoonright
B}(\bar\lambda\upharpoonright B)$.

We prove by induction on $\lambda$ that for every $B\in(I^*)^+$,
$\cf(\prod\bar\lambda,
<_{I^*+(\kappa\setminus B)})=\lambda(B)$ when
$\lambda(B)\leq\lambda$;
this will suffice (use $B=\kappa$ and
$\lambda=|\prod\limits_{i<\kappa}\lambda_i|^+$).
Given $B$ let $\lambda=\lambda(B)$, by notational change
$\wlog B=\kappa$. By 1.9, $\pcf(\prod\bar\lambda)$ has a last element,
necessarily it is $\lambda=: \lambda(B)$.
Let $\langle f_\alpha\colon\alpha<\lambda\rangle$ be
$<_{J_{<\lambda}[\bar\lambda]}$ increasing cofinal in
$\prod\bar\lambda/J_{<\lambda}[\bar\lambda]$, it clearly
exemplifies $\max\pcf(\bar\lambda)\le\cf(\prod\bar\lambda/I^*)$.
Let us prove the other inequality.
For $A\in J_{<\lambda}[\bar\lambda]\setminus I^*$ choose
$F_A\subseteq\prod\bar\lambda$ which is cofinal in
$\prod\bar\lambda/(I^*+(\kappa\setminus A))$, $|F_A|=\lambda(A)<\lambda$
(exists by the induction hypothesis).
Let $\chi$ be a large enough regular, and we now choose by
induction on $\varepsilon<\theta$, $N_\varepsilon$, $g_\varepsilon$ such
that:
\item{(A)}
\itemitem{(i)} $N_\varepsilon\prec(H(\chi), \in,<^*_\chi)$

\itemitem{(ii)} $\Vert N_\varepsilon\Vert=\lambda$

\itemitem{(iii)} $\langle
N_\varepsilon\colon\xi\le\varepsilon\rangle\in N_{\varepsilon+1}$

\itemitem{(iv)} $\langle N_\varepsilon\colon\varepsilon<\theta\rangle$
is increasing continuous

\itemitem{(v)} $\{\varepsilon\colon\varepsilon\le\lambda+1\}\subseteq
N_0$, $\{\bar\lambda, I^*\}\in N_0$, $\langle
f_\alpha\colon\alpha<\lambda\rangle\in N_0$ and the function
$A\mapsto F_A$ belongs to $N_0$.

\item{(B)}
\itemitem{(i)} $g_\varepsilon\in\prod\bar\lambda$ and $g_\varepsilon\in
N_{\varepsilon+1}$

\itemitem{(ii)} for no $f\in N_\varepsilon\cap\prod\bar\lambda$ does
$g_\varepsilon<_{I^*}f$

\itemitem{(iii)} $\zeta<\varepsilon\& \lambda_i>|\varepsilon|
\Rightarrow g_\zeta(i)<g_\varepsilon(i)$.

There is no problem to define $N_\varepsilon$, and if we cannot
choose $g_\varepsilon$ this means that
$N_\varepsilon\cap\prod\bar\lambda$ exemplifies $\cf(\prod\bar\lambda,
<)\le\lambda$ as required.
So assume $\langle N_\varepsilon, g_\varepsilon\colon\varepsilon<\theta
\rangle$ is defined.
For each $\varepsilon<\theta$ for some $\alpha(\varepsilon)<\lambda$,
$g_\varepsilon< f_{\alpha(\varepsilon)}\mod J_{<\lambda}[\bar\lambda]$
hence $\alpha(\varepsilon)\le\alpha<\lambda\Rightarrow
g_\varepsilon<_{J_{<\lambda}[\bar\lambda]}f_\alpha$.
As $\lambda=\cf(\lambda)> \theta$, we can
choose $\alpha<\lambda$ such that
$\alpha>\bigcup_{\varepsilon<\theta}\alpha(\varepsilon)$.
Let $B_\varepsilon=\{ i<\kappa\colon g_\varepsilon(i)\ge
f_\alpha(i)\}$;
so for each $\xi<\theta$ we have $\langle B_\varepsilon\cap
A^*_\xi\colon\varepsilon<\xi\rangle$ is
increasing with $\varepsilon$, (by clause (B)(iii)), hence as
usual as $\theta\ge\wsat(I^*)$ (and $\theta>\aleph_0$) we can find
$\varepsilon(*)<\theta$ such that $\bigwedge_nB_{\varepsilon(*)+n}=
B_{\varepsilon(*)}\mod I^*$ [why do we not demand $\varepsilon\in
(\varepsilon(*), \theta) \Rightarrow B_\varepsilon = B_{\varepsilon(*)}
\ \mod\ I^*$? as $\theta$ may be singular].
Now as $g_{\varepsilon(*)}\in N_{\varepsilon(*)+1}$ and $f_\alpha\in
N_0\prec N_{\varepsilon(*)+1}$ clearly, by its definition,
$B_{\varepsilon(*)}\in
N_{\varepsilon(*)+1}$ hence $F_{B_{\varepsilon(*)}}\in
N_{\varepsilon(*)+1}$.
Now:
$$
\eqalign{
g_{\varepsilon(*)+1}\upharpoonright (\kappa\setminus
B_{\varepsilon(*)}) =_{I^*} g_{\varepsilon(*)+1}\upharpoonright
(\kappa\setminus
B_{\varepsilon(*)+1}) &<  f_\alpha\upharpoonright(\kappa\setminus
B_{\varepsilon(*)+1}) \cr
& =_{I^*} f_\alpha\upharpoonright (\kappa\setminus
B_{\varepsilon(*)}) \cr}
$$
[why first equality and last equality? as
$B_{\varepsilon(*)+1}=B_{\varepsilon(*)}\ \mod\ I^*$, why the $<$ in the
middle? by the definition of $B_{\varepsilon(*)+1}$].

But $g_{\varepsilon(*)+1}\upharpoonright B_{\varepsilon(*)}\in
\prod\limits_{i\in B_{\varepsilon(*)}}\lambda_i$, and
$B_{\varepsilon(*)}\in J_{<\lambda}[\bar \lambda]$ as $g_\varepsilon <
f_{\alpha(\varepsilon)}\leq f_\alpha\ \mod\ J_{<\lambda}[\bar
\lambda]$ so for some $f\in
F_{B_{\varepsilon(*)}}\subseteq \prod\bar\lambda$ we have
$g_{\varepsilon(*)+1}
\upharpoonright B_{\varepsilon(*)}< f\upharpoonright
B_{\varepsilon(*)}\mod I^*$.
By the last two sentences
$$
g_{\varepsilon(*)+1}<\max \{ f, f_\alpha\} \mod I^*
\leqno(*)
$$
Now $f_\alpha\in N_{\varepsilon(*)+1}$ and $f\in
N_{\varepsilon(*)+1}$ (as $f\in
F_{B_{\varepsilon(*)}}$, $|F_{B_{\varepsilon(*)}}|\leq \lambda$,
$\lambda+1\subseteq N_{\varepsilon(*)+1}$ the function $B \mapsto F_B$
belongs to $N_0 \prec N_{\varepsilon(*)+1}$ and
$B_{\varepsilon(*)}\in N_{\varepsilon(*)+1}$ as
$\{g_{\varepsilon(*)}, f_\alpha\}\in N_{\varepsilon(*)+1}$) so together
$$
\max\{ f, f_\alpha\}\in N_{\varepsilon(*)+1};
\leqno(**)
$$
But $(*)$, $(**)$ together contradict the choice of
$g_{\varepsilon(*)+1}$ (i.e.\ clause (B)(ii)).
\sqed{2.10}

\rem{Definition 2.11}
\item{(1)} We say that $I^*$ satisfies the $\pcf$-th for (the
regular) $(\bar\lambda, \theta)$ if 
$\prod \bar \lambda/ I^*$ is $\theta$-directed and 
$(\prod \bar \lambda, <_{J_{<\lambda}[\bar\lambda]})$ is 
$\lambda$-directed and
we can find $\langle
B_\lambda\colon\lambda\in\pcf_{I^*}(\bar\lambda)\rangle$, such
that:
\itemitem{} $B_{\lambda}\subseteq\kappa$,
$J_{<\lambda}[\bar\lambda, I^*]= I^*+\{ B_\mu\colon\mu\in\lambda\cap
\pcf_{I^*}(\bar\lambda)\}$,
$B_\lambda\not\in J_{<\lambda}[\bar\lambda, I^*]$ and
$\Prod(\bar\lambda\upharpoonright B_\lambda)/J_{<\lambda}[\bar\lambda,
I^*]$ has true cofinality $\lambda$ (so $B_\lambda \in J_{\leq
\lambda}[\bar \lambda]\setminus J_{<\lambda}[\bar \lambda]$ and $J_{\leq
\lambda}[\bar \lambda] = J_{<\lambda}[\bar \lambda]+B_\lambda$).

\item{(1A)} We say that $I^*$ satisfies the weak $\pcf$-th for
$(\bar\lambda, \theta)$ if
\itemitem{} $(\prod\bar \lambda, <_{I^*})$ is $\theta$-directed

\itemitem{} each $(\prod \bar \lambda, <_{J_{<\lambda}[\bar \lambda]})$
is $\lambda$-directed and

\itemitem{} there are $B_{\lambda, \alpha}\subseteq \kappa$ for
$\alpha<\lambda\in \pcf_{I^*}(\bar \lambda)$ such that
$$
\alpha<\beta<\mu\in \pcf_{I^*}(\bar \lambda) \Rightarrow B_{\mu,
\alpha} \subseteq B_{\mu, \beta}\ \mod\ J_{<\mu}[\bar \lambda,
I^*]
$$
$$
J_{<\lambda}[\bar \lambda]= I^*+\{B_{\mu, \alpha}\colon \alpha<\mu<\lambda,
\mu\in \pcf_{I^*}(\bar \lambda)\}
$$
$$\hbox{and } 
(\prod\bar \lambda, <_{J_{<\lambda}[\bar \lambda]})\hbox{ is }
\lambda\hbox{-directed and }
$$ 
$$(\prod(\bar \lambda \restriction B_{\mu, \alpha}),
<_{J_{<\lambda}[\bar \lambda]})\hbox{ has true cofinality }\lambda
$$

\item{(1B)} We say that $I^*$ satisfies the weaker $\pcf$-th for $(\bar
\lambda, \theta)$ $\underline{\rm if}$ 
$(\prod\bar \lambda, <_{I^*})$ is $\theta$-directed and each $(\prod
\bar \lambda, <_{J_{<\lambda}[\bar \lambda})$ 
is $\lambda$-directed and for any
ultrafilter $D$ on $\kappa$ disjoint to $J_{<\theta}[\bar \lambda]$
letting $\lambda=\tcf(\prod\bar \lambda, <_D)$ we have: $\lambda\geq
\theta$ and for some $B\in D\cap J_{\leq \lambda}[\bar \lambda]
\setminus J_{<\lambda}[\bar \lambda]$, the partial order $(\prod (\bar
\lambda \restriction B), <_{J_{<\lambda}[\bar \lambda]})$ has true
cofinality $\lambda$.

\item{(1C)} We say that $I^*$ satisfies the weakest $\pcf$-th for $(\bar
\lambda, \theta)$ $\underline{\rm if}$ $(\prod\bar \lambda, <_{I^*})$ is
$\theta$-directed and $(\prod\bar \lambda, <_{J_{<\lambda}[\bar
\lambda]})$ is $\lambda$-directed for any $\lambda\geq \theta$ 

\item{(1D)} Above we write $\bar \lambda$ instead $(\bar \lambda,
\theta)$ when we mean
$$
\theta= \max\{\theta: (\prod\bar \lambda, <_{I^*})\hbox{ is }
\theta^+\hbox{-directed}\}.
$$

\item{(2)} We say that $I^*$ satisfies the $\pcf$-th for
$\theta$ if for any regular $\bar\lambda$ such that
$\lim\inf_{I^*}(\bar\lambda)\geq\theta$, we have: $I^*$ satisfies
the $\pcf$-th for $\bar\lambda$.
We say that $I^*$ satisfies the $\pcf$-th above $\mu$ if it
satisfies the $\pcf$-th for $\bar\lambda$ with
$\lim\inf_{I^*}(\bar\lambda)>\mu$. Similarly (in both cases) for the
weak $\pcf$-th and the weaker $\pcf$-th.

\item{(3)} Given $I^*$, $\theta$ let $J_\theta^\pcf=\{ A\subseteq\kappa
\colon A\in I^*$ or $A\not\in I^*$ and $I^*+(\kappa\setminus A)$
satisfies the $\pcf$-theorem for $\theta\}$.\newline
$J_\theta^\wsat=:\{ A\subseteq\kappa\colon \wsat(I^*\upharpoonright
A)\le\theta\hbox{ or }A\in I^*\}$;\newline
similarly $J_\theta^{{\rm w}\pcf}$; we may write $J_\theta^x[I^*]$.

\item{(4)} We say that $I^*$ satisfies the pseudo $\pcf$-th for $\bar
\lambda$ $\underline{\rm if}$ for every ideal $I$ on $\kappa$ extending
$I^*$, for some $A\in I^+$ we have $(\prod(\bar \lambda\restriction A),
<_I)$ has a true cofinality.

\proclaim{Claim 2.12}
\item{(1)} If $(*)$ of 1.5 then $I^*$ satisfies the weak $\pcf$-th for
$(\bar\lambda, \theta^+)$.

\item{(2)} If $(*)$ of 1.5 holds, and $\prod\bar \lambda/I^*$ is
$\theta^{++}$-directed (i.e. $\theta^+<\min\bar\lambda$) or just there is
a continuity condition for $(\theta^+, \theta)$)
$\underline{\rm then}$
$I^*$ satisfies the $\pcf$-th for $(\bar \lambda, \theta^+)$.

\item{(3)} If $I^*$ satisfy the $\pcf$-th for $(\bar \lambda, \theta)$
$\underline{\rm then}$
$I^*$ satisfy the weak $\pcf$-th for $(\bar \lambda, \theta)$ which implies
that $I^*$ satisfies the weaker $\pcf$-th for $(\bar \lambda, \theta)$,
which implies that $I^*$ satisfies the weakest $\pcf$-th for $(\bar
\lambda, \theta)$.

\endproclaim

\demo{Proof}
(1) Let appropriate $\bar\lambda$ be given. By 1.5, 1.8 most demands
holds, but we are left with normality. 
By 2.8, if $\lambda\in\pcf(\bar\lambda)$, then $\bar\lambda$ is
semi normal for $\lambda$. This finishing the proof of (1).

(2) Let $\lambda\in \pcf(\bar \lambda)$ and let $\bar f$, $\bar B$ be
as in 2.2(4). 
By 2.4(1)+(2) there is $\bar a$, a $(\lambda,\theta)$-continuity condition;
by 2.5(1) $\wlog\bar f$ obeys $\bar a$, by 2.6(1) the relevant
$B_\alpha/ I^*$ are eventually constant which suffices by 2.2(2).

(3) Should be clear.
\sqed{2.12}

\proclaim{Claim 2.13}
Assume $(\prod \bar \lambda, <_{I^*})$ is given (but possibly
$(*)$ of 1.5 fails).
\item{(1)} If $I^*$, $\bar \lambda$ satisfies (the conclusions of) 1.6,
$\underline{\rm then}$ $I^*$, $\bar \lambda$ satisfy (the conclusion of)
1.8(1), 1.8(2), 1.8(3), 1.8(4), 1.9.

\item{(1A)} If $I^*$ satisfies the weaker $\pcf$-th for $\bar \lambda$
$\underline{\rm then}$ they satisfy the conclusion of 1.6 (and 1.5).

\item{(2)} If $I^*$, $\bar \lambda$ satisfies (the conclusion of) 1.5
$\underline{\rm then}$ $I^*$, $\bar \lambda$
satisfies (the conclusion of) 1.10.

\item{(2A)} If $I^*$ satisfies the weakest $\pcf$-th for $\bar
\lambda$ $\underline{\rm then}$ $I^*$, $\bar \lambda$ satisfy the
conclusion of 1.5. 
\item{(3)} If $I^*$, $\bar \lambda$ satisfies 1.5, 1.6
$\underline{\rm then}$ $I^*$, $\bar \lambda$ satisfies
2.2(1) (for 2.2(2) - no
assumptions).

\item{(4)} If $I^*$, $\bar \lambda$ satisfies 1.8(1), 1.8(2)
$\underline{\rm then}$ $I^*$, $\bar \lambda$ satisfies 2.2(3)

\item{(5)} If $I^*$, $\bar \lambda$ satisfies 1.8(2)
$\underline{\rm then}$ $I^*$, $\bar \lambda$ satisfies 2.10(1).

\item{(6)} If $I^*$ $\bar \lambda$ satisfy 1.8(1) + 1.8(3)(i)
$\underline{\rm then}$ $I^*$, $\bar \lambda$ satisfies 1.8(2)

\item{(7)} If $I^*$, $\bar \lambda$ satisfies 1.8(1) + 1.8(2) and is
semi normal 
$\underline{\rm then}$ 2.10(2) holds i.e. 
$$\cf(\prod\bar \lambda,
<_{I^*})\leq \sup\pcf_{I^*}(\lambda).
$$

\endproclaim

\demo{Proof}
(1) We prove by parts.

\subdemo{proof of 1.8(2)} Let $\lambda= \tcf(\prod \bar \lambda /D)$; by
the definition of $\pcf$, $D\cap J_{<\lambda}[\bar \lambda]=\emptyset$.
Also by 1.6 for some $B\in D$ we have $\lambda=\tcf(\prod(\bar
\lambda\restriction  B), <_{J_{<\lambda}[\bar \lambda]})$, so by the
previous sentence $B\notin J_{<\lambda} [\bar \lambda]$, and by 1.4(5)
we have $B\in J_{\leq \lambda} [\bar \lambda]$, together we finish.

\subdemo{proof of 1.8(1)} Repeat the proof of 1.8(1) replacing the use
of 1.5 by 1.8(2).

\subdemo{proof of 1.8(3)$(i)$} Let $J=: \bigcup\limits_{\mu<\lambda}
J_{<\mu} [\bar \lambda]$, so $J\subseteq J_{<\lambda}[\bar \lambda]$ is
an ideal because $\langle J_{<\mu}[\bar \lambda]: \mu<\lambda\rangle$ is
$\subseteq$-increasing (by
1.4(2)), if equality fail choose $B\in J_{<\lambda} [\bar
\lambda]\setminus J$ and choose $D$ an ultrafilter on $\kappa$ disjoint to
$J$ to which $B$ belongs.
Now if $\mu=\cf(\mu)<\lambda$ then $\mu^+< \lambda$
(as $\lambda$ is a limit cardinal) and
$\mu=\cf(\mu) \
\&\ \mu^+<\lambda \Rightarrow D\cap J_{\leq \mu}[\bar \lambda]= D\cap
J_{<\mu^+}[\bar \lambda]=\emptyset$ hence by 1.8(2) we have $\mu\neq
\cf(\prod\bar \lambda /D)$. Also if $\mu=\cf(\mu)\geq \lambda$ then
$D\cap J_{<\mu} [\bar \lambda] \subseteq D\cap J_{<\lambda}[\bar
\lambda] = \emptyset$ hence by 1.8(2) we have $\mu\neq \cf(\prod\bar
\lambda/ D)$. Together contradiction by 1.3(7).

\subdemo{proof of 1.8(3)$(ii)$} Follows.

\subdemo{proof of 1.8(4)} Follows.

\subdemo{proof of 1.9} As in 1.9.

\noindent
(1A) Check.

\noindent
(2) Read the proof of 1.10.

\noindent
(2A) Check.

\noindent
(3) The direction $\Rightarrow$ is proved directly as in the proof of
2.2(1) (where the juse of 1.8(1) is justified by 2.13(1)).

So let us deal with the direction $\Leftarrow$. So assume $\bar
f=\langle f_\alpha: \alpha<\lambda\rangle$ is a sequence of members of
$\prod \bar \lambda$ which is $<_{J_{<\lambda}[\bar \lambda]}$-increasing such
that for every ultrafilter $D$ on $\kappa$ disjoint to
$J_{<\lambda}[\bar \lambda]$ we have: $\lambda=\tcf(\prod \bar
\lambda, <_D)$ iff 
$\bar f$ is unbounded (equivalently cofinal) in $(\prod \bar \lambda, <_D)$.
By (the conclusion of) 1.5 wlog $\bar f$ is $<_{J_{<\lambda}[\bar
\lambda]}$-increasing, and let
$$
J=:\{A\subseteq \kappa: A\in J_{<\lambda}[\bar \lambda]\hbox{ or }\bar
f\hbox{ is cofinal in }(\prod\bar \lambda, <_{J_{<\lambda}[\bar
\lambda]+(\kappa\setminus A)}\}.
$$
Clearly $J$ is an ideal on $\kappa$ (by 1.3(2)(v)), and $J_{<\lambda}[\bar
\lambda]\subseteq J\subseteq J_{\leq \lambda}[\bar \lambda]$. If $J\neq
J_{<\lambda}[\bar \lambda]$ choose $A\in J_{\leq \lambda}[\bar
\lambda]\setminus J$ and an ultrafilter $D$ on $\kappa$ disjoint to $J$
to which $A$ belongs.

By (the conclusion of) 1.6, there is $A\in J\cap D$; contradiction, so
actually $J=J_{\leq \lambda}[\bar \lambda]$. By 1.5 there is $g\in
\prod\bar \lambda$ such that $f_\alpha<g\ \mod\ J_{\leq \lambda}[\bar
\lambda]$ for each $\alpha<\lambda$, and let $B_\alpha=:\{i<\kappa:
g(i)\leq f_\alpha(i)\}$. Hence $B_\alpha\in J_{\leq \lambda}[\bar
\lambda]$ (by the previous sentence) and $\langle
B_\alpha/J_{<\lambda}[\bar \lambda]: \alpha<\lambda\rangle$ is
$\subseteq$-increasing (as $\langle f_\alpha: \alpha<\lambda\rangle$ is
$<_{J_{<\lambda}[\bar \lambda]}$-increasing). Lastly if $B\in J_{\leq
\lambda}[\bar \lambda]$, but $B\setminus B_\alpha\notin
J_{<\lambda}[\bar \lambda]$ for each $\alpha<\lambda$, let $D$ be an
ultrafilter on $\kappa$ disjoint to $J_{<\lambda}[\bar \lambda]+
\{B_\alpha: \alpha<\lambda\}$ but to which $B$ belongs, so $\tcf(\prod
\bar \lambda, <_D)=\lambda$ (by 1.8(3) which holds by 2.12(1)) but
$\{f_\alpha/D: \alpha<\lambda\}$ is bounded by $g/D$ (as $f_\alpha/D\leq
g/D$ by the definition of $B_\alpha$), contradiction. So the sequence
$\langle B_\alpha: \alpha<\lambda\rangle$ is as required.

\noindent
4) -- 6) Left to the reader.

\noindent
7) Let for $\lambda\in\pcf(\bar \lambda)$, $\langle B^\lambda_i:
i<\lambda\rangle$ be such that $J_{\leq \lambda} [\bar \lambda] =
J_{<\lambda}[\bar \lambda] + \{B^\lambda_i: i<\lambda\}$ (exists by
seminormality; we use only this equality). Let $\langle f^{\lambda, i}_\alpha:
\alpha<\lambda\rangle$ be cofinal in $(\prod(\bar \lambda \restriction
B^\lambda_i), <_{J_{,\lambda}[\bar \lambda]})$, it exists by 1.8(1). 
Let $F$ be the closure of
$\{f^{\lambda, i}_\alpha: \alpha<\lambda, i<\lambda, \lambda\in
\pcf(\bar \lambda)\}$, under the operation $\max\{g, h\}$. Clearly
$|F|\leq \sup\pcf(\bar \lambda)$, so it suffice to prove that $F$ is a
cover of $(\prod\bar\lambda, <_{I^*})$.
Let $g\in
\prod\bar \lambda$, if $(\exists f\in F)(g\leq f)$ we are done, if not
$$
I=\{ A\cup \{i<\kappa: f(i)>g(i)\}: f\in F, A\in I^*\}
$$
is $\aleph_0$-directed, $\kappa\notin I$, so there is an ultrafilter $D$
on $\kappa$ disjoint to $I$, (so $f\in F \Rightarrow g <_D f$) and let
$\lambda=\tcf(\prod \bar \lambda /D)$, so by 1.8(2) we have $D\cap J_{\leq
\lambda} [\bar \lambda] \setminus J_{<\lambda}[\bar \lambda] \neq
\emptyset$, hence for some $i<\lambda$, $B^\lambda_i \in D$, and we get
contradiction to the choice of the $\{f^{\lambda, \alpha}_\alpha:
\alpha<\lambda\}$ ($\subseteq F$).
\sqed{2.13}

\proclaim{Claim 2.14}
If $I^*$ satisfies pseodo $\pcf$-th $\underline{\rm then}$

\item{(1)} $\cf(\prod\bar \lambda, <_{I^*})= \sup \pcf_{I^*}(\bar
\lambda)$

\item{(2)} We can find $\langle (J_{\zeta}, \theta_\zeta): \zeta<
\zeta^*\rangle$, $\zeta^*$ a successor ordinal such that $J_0=I^*$,
$J_{\zeta+1}= \{A\subseteq \kappa: \hbox{ if }A\notin J_\zeta\hbox{ then }
\tcf(\prod(\bar \lambda\restriction A), <_{J_\zeta})=\theta_\zeta\}$ and for
no $A\in (J_\zeta)^+$ does $(\prod (\bar \lambda\restriction A),
<_{J_\zeta})$ has true cofinality which is $<\theta_\zeta$.

\item{(3)} If $I^*$ satisfies the weaker $\pcf$-th for $\bar \lambda$
$\underline{\rm then}$ $I^*$ satisfies the pseudo $\pcf$-th for $\bar
\lambda$.

\endproclaim

\demo{Proof}
1) Similar to the proof of 2.10(2).

2) Check (we can also present those ideals in other ways).

3. Check.
\sqed{2.14}


\head \S3 Reduced products of cardinals\endhead

We characterize here the cardinalities
$\prod\limits_{i<\kappa}\lambda_i/D$ and $T_D(\langle\lambda_i\colon
i<\kappa\rangle)$ using $\pcf$'s and the amount of regularity of
$D$ (in 3.1-3.4).
Later we give sufficient conditions for the existence of
$<_D$-lub or $<_D$-eub.
Remember the old result of Kanamori [Kn] and Ketonen [Kt]: for
$D$ an ultrafilter the sequence
$\langle\alpha/D\colon\alpha<\kappa\rangle$ (i.e.\ the constant
functions) has a $<_D$-lub if $\reg(D)<\kappa$; and see [Sh-g,
III 3.3] (for filters).
Then we turn to depth of ultraproducts of Boolean algebras.

The questions we would like to answer are (restricting ourselves to
``$\lambda_i\ge 2^\kappa$'' or ``$\lambda_i\ge 2^{2^\kappa}$'' and
$D$ an ultrafilter on $\kappa$ will be good enough).
\proclaim{Question A}
What can be $\Car_D =:\{\prod\limits_{i<\kappa}\lambda_i /D\colon
\lambda_i$ a cardinal for $i<\kappa\}$ i.e.\ characterize it by
properties of $D$; (or at least $\Card_D\setminus 2^\kappa$)
(for $D$ a filter also $T_D(\prod\lambda_i)$ is natural).
\endproclaim

\proclaim{Question B}
What can be $\DEPTH_D^+=\{\Depth^+
(\prod\limits_{i<\kappa}\lambda_i/D)\colon\lambda_i$ a regular cardinal$\}$
(at least $\DEPTH_D^+\setminus 2^\kappa$, see Definition 3.18).
\endproclaim

If $D$ is an $\aleph_1$-complete ultrafilter, the answer is clear.
For $D$ a regular ultrafilter on $\kappa$,
$\lambda_i\ge\aleph_0$ the answer to question A is known
([CK]) in fact it was the reason for defining ``regularity of
filters'' (for $\lambda_i<\aleph_0$ see [Sh7], [Sh-a, VI-\S 3 Th 3.12
and pp 357-370]
better [Sh-c VI\S 3] and Koppleberg [Ko].)
For $D$ a regular ultrafilter on $\kappa$, the answer to the question
is essentially completed in 3.22(1), the remaining problem can be
answered by $pp$ (see [Sh-g]) except the restriction $(\forall
\alpha<\lambda)(|\alpha|^{\aleph_0}<\lambda)$, which can be removed
if the $\cov = pp$ problem is
completed (see [Sh-g, AG]).
So the problem is for the other ultrafilters $D$, on which we give
a reasonable amount on information translating to a $\pcf$ problem,
sometimes depending on the $\pcf$
theorem.

\rem{Definition 3.1}
\item{(1)} For a filter $D$ let $\reg(D)=\Min\{\theta\colon D$
is not $\theta$-regular$\}$ (see below).

\item{(2)} A filter $D$ is $\theta$-regular if there are $A_\varepsilon
\in D$ for $\varepsilon<\theta$ such that the intersection of any
infinitely many $A_\varepsilon$-s' is empty.

\item{(3)} For a filter $D$ let
$$\eqalign{ \reg_*(D)=\Min\{\theta\colon &
\hbox{ there are no }A_\varepsilon\in D^+\hbox{ for }\varepsilon<\theta
\hbox{ such that }\cr
\ & \hbox{ no }i<\kappa
\hbox{ belongs to infinitely many }A_\varepsilon\hbox{'s}\}\cr}
$$
and
$$
\eqalign{\reg_\otimes(D)=:\{\theta: & \hbox{ there are no }
A_\varepsilon\in D^+\hbox{ for }\varepsilon<\theta\hbox{ such that}:\cr
\ & \ \varepsilon< \zeta \Rightarrow A_\zeta\subseteq A_\varepsilon\ \mod\
D\hbox{ and no }i<\kappa\cr
\ & \hbox{ belongs to infinitely many
}A_\varepsilon\hbox{'s}\}.\cr}
$$

\item{(4)} $\reg^\sigma(D)= \min\{ \theta: D$ is not $(\theta,
\sigma)$-regular$\}$ where ``$D$ is $(\theta, \sigma)$-regular'' means
that there are $A_\varepsilon \in D$ for 
$\alpha<\theta$ such that the intersection of any $\sigma$ of them is
empty. Lastly $\reg^\sigma_* (D)$, $\reg_\otimes^\sigma(D)$ are
defined similarly using $A_\varepsilon \in D^+$.
Of course $\reg(I)$ etc.\ means $\reg(D)$ where $D$ is the dual filter.

\rem{Definition 3.2}
\item{(1)} Let
$$\eqalign{\htcf_{D,\mu}(\Prod\gamma_i)=\sup\{\tcf
(\Prod_{i<\kappa}\lambda_i/ D)\colon & \mu\le\lambda_i=\cf\lambda_i\le
\gamma_i\hbox{ for }i<\kappa\hbox{ and }\cr
\ &\ \tcf(\Prod\lambda_i/D)\hbox{ is well defined}\}\hbox{\ \ and}\cr}
$$
$$\hcf_{D,\mu}(\Prod_{i<\kappa}\gamma_i)=
\sup\{\cf(\Prod\lambda_i /D)\colon\mu\le\lambda_i=
\cf\lambda_i\le\gamma_i\};
$$
if $\mu=\aleph_0$ we may omit it.

\item{(2)} For $E$ a family of filters on $\kappa$
let $\hcf_{E,\mu}(\Prod_{i<\kappa}\alpha_i)$ be
$$\eqalign{\sup\{\tcf(\Prod_{i<\kappa}\lambda_i/D)\colon & D\in E
\hbox{ and }
\mu\le\lambda_i=\cf\lambda_i\le\alpha_i\hbox{ for }
i<\kappa\hbox{ and }\cr
\ & \ \tcf(\Prod_{i<\kappa}\lambda_i/D)\hbox{ is well defined}\}.\cr}
$$
Similarly for $\hcf_{E,\mu}$ (using $\cf$ instead $\tcf$).

\item{(3)} $\hcf^*_{D,\mu}(\Prod_{i<\kappa}\alpha_i)$ is
$\hcf_{E, \mu}(\Prod_{i<\kappa}\alpha_i)$ for $E=\{ D'\colon D'$ a filter
on $\kappa$ extending $D\}$.
Similarly for $\htcf^*_{D,\mu}$.

\item{(4)} When we write $I$ e.g.\ in $\hcf_{I,\mu}$ we mean
$\hcf_{D,\mu}$ where $D$ is the dual filter.

\proclaim{Claim 3.3}
\item{(1)} $\reg (D)$ is always regular

\item{(2)} If $\theta<\reg_*(D)$ $\underline{\rm then}$ some
filter extending $D$
is $\theta$-regular.

\item{(3)} $\wsat(D)\le\reg_*(D)$

\item{(4)} $\reg(D)\le \reg_\otimes(D)\le \reg_*(D)$

\item{(5)} $\reg_*(D)= \min\{\theta:$ no ultrafilter $D_1$ on
$\kappa$ extending $D$ is $\theta$-regular$\}$

\item{(6)} If $D\subseteq E$ are filters on $\kappa$ then:
\itemitem{(a)} $\reg(D)\leq \reg(E)$
\itemitem{(b)} $\reg_*(D)\geq \reg_*(E)$

\endproclaim

\demo{Proof}
Should be clear. E.g (2) let $\langle u_\varepsilon\colon
\varepsilon<\theta\rangle$ list 
the finite subsets of $\theta$, and let $\{ A_\varepsilon\colon
\varepsilon<\theta\}\subseteq D^+$ exemplify ``$\theta<\reg_*(D)$''.
Now let $D^*=:\{ A\subseteq\kappa\colon$ for some finite
$u\subseteq\theta$, for every $\varepsilon<\theta$ we have:
$u\subseteq u_\varepsilon\Rightarrow A_\varepsilon\subseteq A\mod D\}$,
and let $A^*_\varepsilon=\bigcup\{ A_\zeta\colon \varepsilon\in
u_\zeta\}$.
Now $D^*$ is a filter on $\kappa$ extending $D$ and for $\varepsilon
<\theta$ we have $A^*_\varepsilon\in D$.
Finally the intersection of
$A^*_{\varepsilon_0}\cap
A^*_{\varepsilon_1}\cap\ldots$ for distinct $\varepsilon_n<\theta$
is empty, because for any memeber $j$
of it we can find $\zeta_n<\theta$ such that $j\in A_{\zeta_n}$ and
$\varepsilon_n\in u_{\zeta_n}$. Now $\underline{\rm if}$
$\{\zeta_n\colon n<\omega\}$ is
infinite then there is no such $j$ by the choice of
$\langle A_\varepsilon\colon \varepsilon<\theta\rangle$, and
$\underline{\rm if}$ $\{\zeta_n\colon n<\omega\}$ is finite then
$\wlog\ \bigwedge\limits_{n,\omega}\zeta_n=\zeta_0$ contradicting ``$u_{\zeta_0}$ is
finite'' as $\bigwedge\limits_{n<\omega} \varepsilon_n\in u_{\zeta_n}$.
\sqed{3.3}

\rem{Observation 3.4}
$|\Prod_{i<\kappa}\lambda_i/I|\ge\vert\aleph_0^\kappa/ I\vert$
holds when $\bigwedge_{i<\kappa}\lambda_i\ge
\aleph_0$.

\rem{Observation 3.5}
\item{(1)} $|\Prod_{i<\kappa}\lambda_i/I|\ge\htcf^*_{I}(\Prod_{i<\kappa}
\lambda_i)$.

\item{(2)} If $I^*$ satisfies the $\pcf$-th for $\bar\lambda$ or even
the weaker $\pcf$-th or even the pseudo $\pcf$-th for $\bar \lambda$
(see Definition 2.11) 
$\underline{\rm then}$: $\cf(\prod\bar\lambda/ I^*)=\max\pcf_{I^*}
(\bar\lambda)$. 

\item{(3)} If $I^*$ satisfies the $\pcf$-th for $\mu$ for and
$\min(\bar\lambda)\geq \mu$ 
$\underline{\rm then}$
$$\hcf_{D,\mu}(\prod\bar\lambda)=\hcf^*_{D,\mu}(\prod\bar\lambda)=
\htcf^*_{D,\mu}(\prod\bar\lambda)$$
whenever $D$ is disjoint to $I^*$.

\item{(4)} $\hcf_{E, \mu} (\prod\limits_{i<\kappa}\lambda_i)= \hcf^*_{E,
\mu}(\prod\limits_{i<\kappa} \lambda_i)$.

\item{(5)} $\prod\limits_{i<\kappa}\lambda_i/I\geq \hcf_{I,
\mu}(\prod\limits_{i<\kappa}\lambda_i)=
\hcf^*_{I,\mu}(\prod\limits_{i<\kappa}\lambda_i) \geq\htcf^*_{I,\mu}
(\prod\limits_{i<\kappa}\lambda_i)$ and 
$\hcf_{I,\mu}(\prod\limits_{i<\kappa}\lambda_i)\geq \htcf_{I,
\mu}(\prod\limits_{i<\kappa}\lambda_i)$.

\rem{Remark 3.5A} In 3.5(3) concerning $\htcf_{D, \mu}$ see 3.10. 

\demo{Proof}
1) By the definition of $\htcf^*_I$ it suffices to show
$\vert\Prod_{i<\kappa}\lambda_i/I\vert\ge\tcf(\Prod\lambda'_i/I')$,
when $I'$ is an ideal on $\kappa$ extending $I$,
$\lambda'_i=\cf\lambda'_i\le\lambda_i$ for $i<\kappa$ and
$\tcf(\Prod_{i<\kappa}\lambda'_i/I')$ is well defined.
Now
$|\Prod_{i<\kappa}\lambda_i/I|\ge|\Prod_{i<\kappa}\lambda'_i/I|\ge
|\Prod_{i<\kappa}\lambda'_i/I'|\ge\cf(\Prod\lambda'_i/I')$,
so we have finished.

2) By 2.13(1) and 1.9 and 2.14.

3) Left to the reader (see Definition 2.11(2)).

4), 5) Check.
\sqed{3.5}.

\proclaim{Claim 3.6}
If $\lambda=|\Prod_{i<\kappa}\lambda_i/I|$ (and $\lambda_i\ge\aleph_0$
and, of course,
$I$ an ideal on $\kappa$) and $\theta<\reg(I)$ $\underline{\rm then}$
$\lambda=\lambda^\theta$. 
\endproclaim

\demo{Proof}
For each $i<\kappa$, let
$\langle\eta^i_\alpha\colon\alpha<\lambda_i\rangle$ list the
finite sequences from $\lambda_i$.
Let $M_i=(\lambda_i, F_i, G_i)$ where
$F_i(\alpha)=\lg(\eta^i_\alpha)$, $G_i(\alpha, \beta)$ is
$\eta^i_\alpha(\beta)$ if $\beta<\lg(\eta^i_\alpha)$
($=F_i(\alpha)$), and $F(\alpha, \beta)=0$ otherwise;
let $M=\Prod_{i<\kappa} M_i/I$ so $\Vert M\Vert= |\Prod\lambda_i/I|$ and
let $M=(\Prod\lambda_i/I, F, G)$.
Let $\langle A_i\colon i<\theta\rangle$ exemplifies $I$ is
$\theta$-regular.
Now

\item{$(*)_1$} We can find $f\trans{\in}{\omega}$ and
$f_\varepsilon\in\Prod_{i<\kappa}f(i)$ for $\varepsilon<\theta$ such that:
$\varepsilon<\zeta<\theta\Rightarrow f_\varepsilon<_I f_\zeta$ [just for
$i<\kappa$ let $w_i=\{\varepsilon<\theta\colon i\in A_\varepsilon\}$,
it is finite and let $f(i)=|w_i|$ and $f_\varepsilon(i)=|\varepsilon\cap
w_i|<f(i)$, and note $\varepsilon<\zeta\& i\in A_\varepsilon\cap A_\zeta
\Rightarrow f_\varepsilon(i)<f_\zeta(i)$].

\item{$(*)_2$} For every sequence $\bar g= \langle g_\varepsilon\colon
\varepsilon<\theta\rangle$ of members of $\Prod_{i<\kappa}\lambda_i$,
there is $h\in \Prod_{i<\kappa}\lambda_i$ such that $\varepsilon<\theta
\Rightarrow M\vDash F(h/I, f_\varepsilon/I)= g_\varepsilon/I$
[why? let, in the notation of $(*)_1$, $h(i)$ be such that
$\eta^i_{h(i)}=\langle g_\varepsilon(i)\colon \varepsilon\in
w_i\rangle$ (in the natural order)].

\noindent
So in $M$, every $\theta$-sequence of members is coded by at least
one member so $\Vert M\Vert^\theta=\Vert M\Vert$, but $\Vert M\Vert
=|\Prod_{i<\kappa}\lambda_i/I|$ hence we have proved 3.6.
\sqed{3.6}

\rem{Fact 3.7}
\item{(1)} For $D$ a filter on $\kappa$, $\langle A_1,
A_2\rangle$ a partition of $\kappa$ and (non zero) cardinals $\lambda_i$
for $i<\kappa$ we have
$$|\Prod_{i<\kappa}\lambda_i/D|=|\Prod_{i<\kappa}\lambda_i/(D\shvor
+A_1)|
\times |\Prod_{i<\kappa}\lambda_i/(D+A_2)|
$$
(note:
$|\Prod_{i<\kappa}\lambda_i/\calP(\kappa)|=1$).

\item{(2)} $D^{[\mu]}=:\{
A\subseteq\kappa\colon|\Prod_{i<\kappa}\lambda_i /(D+(\kappa\setminus
A))|<\mu\}$ is a filter on $\kappa$ ($\mu$ an infinite cardinal of
course) and if
$\aleph_0\le\mu\le\Prod_{i<\kappa}\lambda_i/D$ then $D^{[\mu]}$
is a proper filter.

\item{(3)} If $\lambda\leq |\prod\limits_{i<\kappa} \lambda_i/ I|$,
($\lambda_i$ infinite, of course, $I$ an ideal on $\kappa$)
and $A\in I^+ \Rightarrow |\prod\limits_{i\in A}\lambda_i / I|\geq
\lambda$ and $\sigma <\reg_\otimes (I)$ $\underline{\rm then}$
$|\prod\lambda_i/I|\geq \lambda^\sigma$

\demo{Proof} Check (part (3) is like 3.6).

\proclaim{Claim 3.8}
If $D\subseteq E$ are filters on $\kappa$ then
$$
|\Prod_{i<\kappa}\lambda_i/D|\le |\Prod_{i<\kappa}\lambda_i/E|+
\sup_{A\in E\setminus D}|\Prod_{i<\kappa}\shvor\lambda_i/
(D+(\kappa\setminus A))|+
(2^\kappa/D)+\aleph_0.
$$
We can replace $2^\kappa/D$ by $|\calP|$ if $\calP$ is a maximal
subset of $E$ such that $A\not= B\in\calP\Rightarrow
(A\setminus B)\cup(B\setminus A)\not= \emptyset \mod D$.
\endproclaim

\demo{Proof} Think.

\proclaim{Lemma 3.9}
$|\Prod_{i<\kappa}\lambda_i/D|\le
(\theta^\kappa/D+\hcf_{D,\theta}(\Prod_{i<\kappa}\lambda_i))^{<\theta}$
(see Definition 3.2(1)) provided that:
$$
\theta\ge\reg_\otimes(D)
\leqno(*)
$$
\endproclaim

\rem{Remark 3.9A}
1) If $\theta=\theta^+_1$, we can replace $\theta^\kappa/D$ by
$\theta_1^\kappa/D$. In general we can replace $\theta^\kappa/D$ by
$\sup\{\prod\limits_{i<\kappa} f(i)/D: f\in \theta^\kappa\}$.

\noindent
2) If $D$ satisfies the $\pcf$-th above $\theta$
(see 2.11(1A), 2.12(2)) then by 3.5(3) we can use $\htcf^*$ (sometime
even $\htcf$, see 3.10). But by 3.7(1)
we can ignore the $\lambda_i\leq \theta$, and when $i<2 \Rightarrow
\lambda_i> \theta$ we know that 1.5(*)$(\alpha)$ holds by 3.3(3).

\demo{Proof}
Let
$\lambda=\theta^\kappa/D+ \hcf_{D,\theta}(\prod\limits_{i<\kappa}\lambda_i)$.
Let for
$\zeta<\theta$, $\mu_\zeta =: \lambda^{\vert\zeta\vert}$
i.e.\ $\mu_\zeta =:
(\theta^\kappa/D+\hcf_{D,\theta}\prod\limits_{i
<\kappa} \lambda_i)^{|\zeta|}$,
clearly $\mu_\zeta=\mu_\zeta^{|\zeta|}$.
Let $\chi=\beth_8(\sup_{i<\kappa}\lambda_i)^+$ and
$N_\zeta\prec(H(\chi), {\in}, <^*_\chi)$ be such that
$\Vert N_\zeta\Vert= \mu_\zeta$,
$N^{\le|\zeta|}\subseteq N_\zeta$, $\lambda+1\subseteq N_\zeta$
and $\{ D, \langle\lambda_i\colon
i<\kappa\rangle\}\in N_\zeta$ and $[\varepsilon<\zeta\Rightarrow
N_\varepsilon\prec N_\zeta]$.
Let $N=\cup\{ N_\zeta\colon\zeta<\theta\}$.
Let $g^*\in \Prod_{i<\kappa}\lambda_i$ and we shall find $f\in N$
such that $g^*=f\mod D$, this will suffice.
We shall choose by induction on $\zeta<\theta$, $f_\zeta^e(e<3)$
and ${\bar A}^\zeta$ such that:
\item{(a)} $f^e_\zeta\in\Prod_{i<\kappa}(\lambda_i+1)$

\item{(b)} $f^1_\zeta\in N_\zeta$ and $f^2_\zeta\in N_\zeta$.

\item{(c)} $\bar A^\zeta=\langle A^\zeta_i\colon i<\kappa
\rangle\in N_\zeta$.

\item{(d)} $\lambda_i\in A^\zeta_i\subseteq\lambda_i+1$,
$|A^\zeta_i|\leq |\zeta|+1$, and
$\langle A^\zeta_i\colon \zeta<\theta\rangle$ is increasing
continuous (in $\zeta$).

\item{(e)} $f_\zeta^0(i)=\Min(A^\zeta_i\setminus g^*(i))$;
note: it is well defined as $g^*(i)<\lambda_i\in
A^\zeta_i$

\item{(f)} $f_\zeta^1=f_\zeta^0\mod D$

\item{(g)} $g^*<f_{\zeta}^2<f_\zeta^1\mod (D+\{i<\kappa\colon
g^*(i)\not= f_\zeta^1(i)\})$.

\item{(h)} if $g^*(i)\neq f^1_\zeta(i)$ then $f^2_\zeta(i)\in
A^{\zeta+1}_i$

So assume everything is defined for every $\varepsilon<\zeta$.
If $\zeta=0$, let $A^\zeta_i=\{\lambda_i\}$, if $\zeta$ limit
$A^\zeta_i=\bigcup_{\varepsilon<\zeta} A_i^\varepsilon$,
for $\zeta=\varepsilon+1$,
$A_i^\zeta$ will be defined in stage $\varepsilon$.
So arriving to $\zeta$, $\bar A^\zeta$ is well defined and it belongs
to $N_\zeta$: for $\zeta=0$ check, for $\zeta=\varepsilon+1$, done in
stage $\varepsilon$,
for $\zeta$ limit it belongs to $N_\zeta$ as
we have $N_\zeta^{\le|\zeta|}\subseteq N_\zeta$ and $\xi<\zeta
\Rightarrow N_\xi \prec N_\zeta$.
Now use clause (e) to define $f_\zeta^0/D$.
As $\langle A_i^\zeta\colon i<\kappa\rangle\in N_\zeta$, $|A_i^\zeta|
<\theta$ and
$\theta^\kappa/D\le\lambda<\lambda+1\subseteq N_\zeta$, clearly
$|\prod\limits_{i<\kappa} |A^\zeta_i|/ D|\leq \lambda$ hence $\{f/D: f\in
\prod\limits_{i<\kappa} A^\zeta_i\}\subseteq N_\zeta$
hence $f^0_\zeta/D\in N_\zeta$ hence there is $f_\zeta^1\in N_\zeta$ such
that $f_\zeta^1\in f^0_\zeta/D$ i.e. clause (f) holds.
As $g^*\leq f^0_\zeta$ clearly $g^*\le f_\zeta^1\mod D$,
let $y_0^\zeta=:\{i<\kappa\colon
g^*(i)\ge f_\zeta^1(i)\}$,
$y_1^\zeta=:\{i<\kappa\colon i\not\in y_0^\zeta$ and
$\cf(f_\zeta^1(i))<\theta\}$ and
$y_2^\zeta=:\kappa\setminus y_0^\zeta\setminus y_1^\zeta$.
So $\langle y^\zeta_e\colon e<3\rangle$ is a partition of $\kappa$ and
$g^*<f_\zeta^1\mod(D+y_e^\zeta)$ for $e=1,2$.

Let $y^\zeta_4=\{i<\kappa: \cf(f^1_\zeta(i))\geq \theta\}$ so
$f^1_\zeta\in N_\zeta$, and $\theta\in N_\zeta$ hence $y^\zeta_4 \in
N_\zeta$, so $(\prod\limits_{i<\kappa} f^1_\zeta (i), <_{D+
y^\zeta_4})\in N_\zeta$. Now
$$
\cf(\Prod_{i<\kappa}f_\zeta^1(i),
<_{D+y_4^\zeta})\le \hcf_{D+y_4^\zeta, \theta}(\Prod_{i<\kappa}
\lambda_i)
\le\hcf_{D,\theta}(\Prod_{i<\kappa}\lambda_i)
\subseteq\lambda<\lambda+1\subseteq N_\zeta
$$
hence there is $F\in N_\zeta$, $|F|\leq \lambda$, $F\subseteq
\prod\limits_{i\in y^\zeta_4} f^1_\zeta(i)$ such that:
$$
(\forall g) [g\in \prod\limits_{i\in y^\zeta_4} f^1_\zeta(i) \Rightarrow
(\exists f\in F)(g< f\ \mod\ (D+y^\zeta_4)))].
$$
As $\lambda+1\subseteq N$ necessarily $F\subseteq N_\zeta$. Apply the
property of $F$ to $(g\restriction y^\zeta_2)\cup 0_{(\kappa \setminus
y^\zeta_2)}$ and get $f^\zeta_4\in F\subseteq N$ such that $g^* <
f^\zeta_4\ \mod\ (D+ y^\zeta_2)$. Now use similarly
$\prod\limits_{i<\kappa} \cf(f_\zeta^1(i))/(D+y_1^\zeta)\le |\theta^\kappa /D
|\le\lambda$; by the proof of 3.7(1) there is a function
$f_\zeta^2\in N_\zeta\cap\Prod_{i<\kappa}f^1_\zeta(i)$
such that $g^*\upharpoonright(y_1^\zeta+y_2^\zeta)< f_\zeta^2\mod D$.
Let $A_i^{\zeta+1}$ be:
$A_i^\zeta$ if $i\in y_0^\zeta$ and
$A_i^\zeta\cup
\{f_\zeta^2(i)\}$ if $i\in y^\zeta_1 \cup y_2^\zeta$.

It is easy to check clauses (g), (h).
So we have carried the definition.
Let
$$
X_\zeta=:\{ i<\kappa\colon  f_{\zeta+1}^0(i)< f_\zeta^0(i)\}.
$$
Note that by the choice of $f_\zeta^1$, $f_{\zeta+1}^1$ we know 
$X_\zeta=y_1^\zeta\cup y_2^\zeta\ \mod\ D$,
if this last set is not $D$-positive then
$g^*\ge f_\zeta^1\mod D$, hence $g^*/D=f_\zeta^1/D\in N_\zeta$,
contradiction,
so $y^\zeta_1\cup y_2^\zeta\not=\emptyset\mod D$ hence $X_\zeta\in
D^+$. Also $\langle y^\zeta_1\cup y^\zeta_2: \zeta<\theta\rangle$ is
$\subseteq$-decreasing hence $\langle X_\zeta/D: \zeta< \theta\rangle$
is $\subseteq$-decreasing.

Also if $i\in X_{\zeta_1}\cap X_{\zeta_2}$ and $\zeta_1< \zeta_2$ then
$f_{\zeta_2}^0(i)\le f_{\zeta_1+1}^0(i)<f_{\zeta_1}^0(i)$
(first inequality: as $A^{\zeta_1+1}_i\subseteq A_i^{\zeta_2}$ and
clause (e) above, second inequality by the definition of $X_{\zeta_1}$),
hence for each ordinal $i$ the set
$\{ \zeta<\theta\colon i\in X_\zeta\}$ is finite.
So $\theta<\reg_\otimes(D)$, contradiction to the assumption $(*)$.
\sqed{3.9}

Note we can conclude

\proclaim{Claim 3.9B}
$$\eqalign{
\prod\limits_{i<\kappa}\lambda_i/D =
\sup\{(\prod_{i<\kappa}f(i))^{<\reg_\otimes(D_1)} + & \hcf_{D_1}(\prod\limits_{i<\kappa}\lambda_i)^{<\reg_\otimes(D_1)}\colon
D_1\hbox{ is a filter on }\kappa\cr
\ & \hbox{ extending }D\hbox{ such that }\cr
\ & \qquad\qquad A\in D^+_1 \Rightarrow \prod_{i<\kappa}
\lambda_i/ (D_1+A)= \prod_{i<\kappa}\lambda_i/D_1\cr
\ & \qquad \qquad\hbox{ and } f\in \theta^\kappa, f(i)\leq \lambda_i\}\cr}
$$
\endproclaim

\demo{Proof} The inequality $\geq$ should be clear by 3.7(3). For the
other direction let $\mu$ be the right side cardinality and let
$D_1=\{\kappa \setminus A:$ if $A\in D^+$ then
$\prod\limits_{i<\kappa} \lambda_i/D \leq \mu\}$, so we know by 3.7(2)
that $D_1$ is a filter on $\kappa$ extending $D$. Now $\mu\geq
{\aleph_0}^\kappa/D$ (by the term $(\prod\limits_i f(i)/
D_1)^{<\reg_\otimes(D_1)}$) so by 3.8 we have $\prod\limits_{i<\kappa}
\lambda_i/D_1>\mu$. By 3.9 (see 3.9A(1)) we get a contradiction.
\sqed{3.9B}

Next we deal with existence of $<_D-\eub$.

\proclaim{Claim 3.10}
1) Assume $D$ a filter on $\kappa$, $g^*_\alpha\trans \in\Ord$ for
$\alpha<\delta$, $\bar g^*=\langle g^*_\alpha\colon \alpha<\delta
\rangle$ is
$\leq_D$-increasing, and
$$
\cf(\delta)\ge\theta\ge\reg_*(D).
\leqno(*)
$$

$\underline{\rm Then}$ at least one of the following holds:
\item{(A)} $\langle g^*_\alpha\colon \alpha<\delta\rangle$ has a
$<_D$-eub $g\trans \in\Ord$; moreover $\theta\le\lim\inf_D\langle\cf
[g(i)]\colon i<\kappa\rangle$

\item{(B)} $\cf(\delta)=\reg_*(D)$

\item{(C)} for some club $C$ of $\delta$ and some $\theta_1<\theta$ and
$\gamma_i<\theta_1^+$ and $w_i \subseteq \Ord$ of order type $\gamma_i$
for $i<\kappa$, there are $f_\alpha\in \Prod_{i<\kappa}
w_i$ (for $\alpha\in C$) such that $f_\alpha(i)=\min(w_i\setminus
g^*_\alpha(i))$ and
$\alpha \in C \ \&\ \beta\in C\ \&\
\alpha<\beta\Rightarrow
f_\alpha\le_D f_\beta\ \&\ \neg f_\alpha=_D f_\beta\ \&\ \neg
f_\alpha\leq_D g^*_\beta\ \&\ g^*_\alpha \leq f_\alpha$.

2) In (C) above if for simplicity $D$ is an ultrafilter we can find
$w_i\subseteq\Ord$, $\otp(w_i)=\gamma_i$,
$\langle\alpha_\xi\colon\xi<\cf(\delta)\rangle$ increasing continuous
with limit $\delta$, and $h_\varepsilon\in\prod\limits_{i<\kappa}w_i$
such that $f_{\alpha_\varepsilon}<_D h_\varepsilon <_D
f_{\alpha_{\varepsilon+1}}$, moreover, $\bigwedge\limits_{i<\kappa}
\gamma_i<\omega$.
\endproclaim

\demo{Proof}
1) Let $\sigma=\reg_*(D)$. We try to choose by induction on
$\zeta<\sigma$, $g_\zeta$,
$f_{\alpha,\zeta}$ (for $\alpha<\delta$), $\bar A^\zeta$,
$\alpha_\zeta$ such that
\item{(a)} $\bar A^\zeta=\langle
A_i^\zeta\colon i<\kappa\rangle$.

\item{(b)} $A_i^\zeta=\{ f_{\alpha_\varepsilon, \varepsilon}(i),
g_\zeta(i)\colon \varepsilon<\zeta\}\cup\{[\sup_{\alpha<\delta}
g_\alpha^*(i)]+1\}$.

\item{(c)} $f_{\alpha, \zeta}(i)=\Min(A_i^\zeta\setminus
g^*_\alpha(i))$ (and $f_{\alpha,\zeta}\trans{\in}{\Ord}$, of
course).

\item{(d)} $\alpha_\zeta$ is the first $\alpha$,
$\bigcup_{\varepsilon<\zeta}\alpha_\varepsilon<\alpha<\delta$ such that
$[\beta\in[\alpha,
\delta)\Rightarrow f_{\beta, \zeta}=f_{\alpha, \zeta}\mod D]$ if
there is one.

\item{(e)} $g_\zeta\le f_{\alpha_\zeta, \zeta}$ moreover $g_\zeta
<\max\{ f_{\alpha_\zeta, \zeta}, 1_\kappa\}$ but for no $\alpha<\delta$
do we have $g_\zeta< \max\{g^*_\alpha, 1\}\ \mod\ D$.

\medskip
Let $\zeta^*$ be the first for which they are not defined (so
$\zeta^*\le\sigma$).
Note
$$
\varepsilon<\xi<\zeta^*\&  \alpha_\xi\le
\alpha<\delta\Rightarrow f_{\alpha_\varepsilon, \varepsilon}=_D
f_{\alpha, \varepsilon}
\& f_{\alpha, \xi}\le f_{\alpha, \varepsilon} \& f_{\alpha, \xi}\not=_D
f_{\alpha, \varepsilon}.
\leqno(*)
$$               
[Why last phrase? applying clause (e) above, second phrase with
$\alpha$, $\varepsilon$ here standing for $\alpha$, $\zeta$ there we get
$A_0 =: \{i<\kappa:
\max\{g^*_{\alpha}(i), 1\}\leq g_\varepsilon(i)\}\in D^+$ and
applying clause (e) above first phrase with $\varepsilon$ here standing
for $\zeta$ there we get $A_1=\{i<\kappa:
g_\varepsilon(i)< f_{\alpha, \varepsilon}(i)$ or
$g_\varepsilon(i)=0=f_{\alpha, \varepsilon}(i)\}\in D$,
hence $A_0 \cap A_1 \in D^+$,
and $g_\varepsilon(i)>0$ for $i\in A_0\cap A_1$ (even for $i\in A_0$).
Also by clause (c) above $g^*_\alpha(i)\leq g_\varepsilon(i) \Rightarrow
f_{\alpha, \xi}(i)\leq g_\varepsilon(i)$. Now by the last two
sentences
$i\in A_0\cap A_1
\Rightarrow g^*_{\alpha}(i)\leq g_\varepsilon(i)< f_{\alpha,
\varepsilon}(i) \Rightarrow f_{\alpha, \xi}(i)\leq g_\varepsilon(i)<
f_{\alpha, \varepsilon}(i)$, together $f_{\alpha, \xi }\not=_D f_{\alpha,
\varepsilon}$ as required]

\subdemoinfo{Case A}{$\zeta^*=\sigma$ and $\bigcup_{\zeta<\sigma}
\alpha_\zeta<\delta$}
Let $\alpha(*)=\bigcup_{\zeta<\sigma}\alpha_\zeta$, for
$\zeta<\sigma$ let $y_\zeta=\{i<\kappa\colon
f_{\alpha(*), \zeta}(i)\not= f_{\alpha(*), \zeta+1}(i)\}
\not=\emptyset\mod D$.
Now for $i<\kappa$, $\langle f_{\alpha(*), \zeta}(i)\colon
\zeta<\sigma\rangle$ is non increasing so $i$ belongs to
finitely many $y_\zeta$'s only, so $\langle y_\zeta\colon
\zeta<\sigma\rangle$ contradict $\sigma\ge\reg_*(D)$.

\subdemoinfo{Case B}{$\zeta^*=\sigma$ and
$\bigcup_{\zeta<\sigma}\alpha_\zeta=\delta$}
So possibility (B) of Claim 3.10 holds.

\subdemoinfo{Case C}{$\zeta^*<\sigma$}

Still $A_i^{\zeta^*}(i<\kappa)$, $f_{\alpha, \zeta^*}(\alpha<\delta)$
are well defined.

\subdemo{Subcase C1}
$\alpha_{\zeta^*}$ cannot be defined.

Then possibility C of 3.10 holds (use $w_i=: A^{\zeta^*}_i$,
$f_\beta=f_{\alpha_{\zeta^*}+\beta, \zeta^*}$).

\subdemo{Subcase C2}
$\alpha_{\zeta^*}$ can be defined.

Then $f_{\alpha_{\zeta^*}, \zeta^*}$ is a $<_D$-eub of $\langle
g^*_\alpha\colon\alpha<\delta\rangle$ as otherwise there is
$g_{\zeta^*}$ as required in clause (e).
Now $f_{\alpha^*_\zeta, \zeta^*}$ is almost as required in possibility
(A) of Claim 3.10 only the second phrase is missing.
If for no $\theta_1<\theta$, $\{i<\kappa\colon\cf
[f_{\alpha_{\zeta^*}, \zeta^*}(i)]\le\theta_1\}\in D^+$, then
possibility (A) holds.

So assume $\theta_1< \theta$ and $B=:\{i<\kappa: \aleph_0\leq
\cf[f_{\alpha_{\zeta^*}, \zeta^*}(i)]\leq \theta_1\}$ belongs to $D^+$,
we shall try to prove that possibility (C) holds, thus finishing. Now we
choose $w_i$ for $i<\kappa$: for $i\in \kappa$ we let $w^0_i=: \{
f_{\alpha_{\zeta^*}, \zeta^*}(i), [\sup\limits_{\alpha<\delta}
g^*_\alpha(i)]+1\}$, for $i\in B$ let $w^1_i$ be an unbounded subset of
$f_{\alpha_{\zeta^*}, \zeta^*}(i)$ of order type
$\cf[f_{\alpha_{\zeta^*}, \zeta^*}(i)]$ and for $i\in \kappa\setminus
B$ let $w^1_i=\emptyset$, lastly let $w_i=w^0_i\cup
w^1_i$, so $|w_i|\leq \theta_1$ as required in possibility (C). Define
$f_\alpha\in {}^\kappa\Ord$ by $f_\alpha(i)=\min(w_i\setminus
g^*_\alpha(i))$ (by the choice of $w^0_i$ it is well defined). So
$\langle f_\alpha: \alpha<\delta\rangle$ is $\leq_D$-increasing; if for
some $\alpha^*<\delta$, for every $\alpha\in [\alpha^*, \delta)$ we have
$f_\alpha/D= f_{\alpha^*}/D$, we could define $g_{\zeta^*}\in
{}^\kappa\Ord$ by:

$g_{\zeta^*}\restriction B= f_{\alpha^*}$ (which is $<
f_{\alpha_{\zeta^*}, \zeta^*}$),

$g_{\zeta^*}\restriction (\kappa\setminus B) = 0_{\kappa\setminus B}$.

\noindent
Now $g_{\zeta^*}$ is as required in clause (e) so we get contradiction
to the choice of $\zeta^*$. So there is no $\alpha^*<\delta$ as above
so for some club $C$ of $\delta$ we have $\alpha< \beta\in C\Rightarrow
f_\alpha \neq_D f_\beta$, so we have actually proved possibility (C).

2) Easy (for $\bigwedge_i\gamma_i<\omega$, $\wlog\theta=\reg_*(D)$ but
$\reg_*(D)=\reg(D)$ so $\theta_1<\reg(D)$).
\sqed{3.10}

\proclaim{Claim 3.11}
\item{(1)} In 3.10(1), if $\lambda=\delta=\cf(\lambda)$, $\bar g^*$ obeys
$\bar a$ ($\bar a$ as in 2.1),
$\bar a$ a $\theta$-weak $(S, \theta)$ continuily condition,
$S\subseteq\lambda$
unbounded, then clause (C) of 3.10 implies:\newline
(C)$'$ there are $\theta_1<\reg_*(D)$ and $A_{\varepsilon}\in D^+$ for
$\varepsilon<\theta$ such that the intersection of any $\theta_1^+$
of the sets $A_\varepsilon$ is empty (equivalently $i<\kappa
\Rightarrow (\exists^{\le\theta_1}\varepsilon)[i\in A_\varepsilon]$
(reminds $(\sigma, \theta_1^+)$-regularity of ultrafilters).

\item{(2)} We can in 3.10(1) weaken the assumption $(*)$ to $(*)'$
below if in the conclusion we weaken clause (A) to (A)$'$ where
\itemitem{$(*)'$} $\cf(\delta)\ge\theta\ge\reg(D)$

\itemitem{(A)$'$} there is a $ \leq_D$-upper bound $f$ of
$\{ g_\alpha^*\colon
\alpha<\delta\}$ such that \newline
no $f'<_D f$ (of course $f'\trans \in\Ord)$
is a $\leq_D$-upper bound of $\{g_\alpha^*\colon
\alpha<\delta\}$\newline
and
$\theta\le\lim\inf_D\langle\cf[f(i)]\colon i<\kappa\rangle$

\item{(3)} If $g_\alpha^*\trans \in\Ord$, $\langle g_\alpha^*\colon
\alpha<\delta\rangle$ is $<_D$-increasing and $f\trans\in\Ord$
satisfies (A)$'$ above and
\itemitem{$(*)''$} $\cf(\delta)\ge\wsat(D)$ and for some $A\in D$ for
every $i<\kappa$, $\cf(f(i))\geq \wsat(D)$

\item{} $\underline{\rm then}$ for some $B\in D^+$ we have
$\prod\limits_{i<\kappa}\cf[f(i)]/(D+B)$ has true cofinality
$\cf(\delta)$.
\endproclaim

\rem{Remark} Compare with 2.6.

\demo{Proof}
1) By the choice of $\bar a=\langle a_\alpha: \alpha<
\lambda\rangle$ as $C$ (in clause (c) of 3.11(1)) is a club of
$\lambda$, we can find $\beta<\lambda$ such that letting $\langle
\alpha_\varepsilon: \varepsilon <\theta\rangle$ list $\{\alpha\in
a_\beta: \otp (\alpha \cap a_\beta)< \theta\}$ (or just a subset of it)
we have $(\alpha_\varepsilon, \alpha_{\varepsilon+1})\cap
C\not=\emptyset$.

Let $\gamma_\varepsilon \in (\alpha_\varepsilon,
\alpha_{\varepsilon+1})\cap C$, and $\xi_\varepsilon \in
(\alpha_\varepsilon, \alpha_{\varepsilon+1})$ be such that
$\{\alpha_\zeta: \zeta\leq \varepsilon \}\subseteq a_{\xi_\varepsilon}$, and as
we can use $\langle \alpha_{2\varepsilon}: \varepsilon< \theta\rangle$,
wlog $\xi_\varepsilon < \gamma_\varepsilon$. For $\zeta< \theta$ let
$B_\zeta=\{i< \kappa: f_{\alpha_\zeta}(i)< f_{\beta_\zeta}(i)<
f_{\gamma_\zeta}(i)< f_{\alpha_{\zeta+1}}(i)$ and
$\sup\{f_{\alpha_\xi}(i)+1: \xi< \zeta\}< \sup \{f_{\alpha_\xi}(i)+1:
\xi< \zeta+1\}$.

\noindent
2) In the proof of 3.10 we replace clause (e) by
\itemitem{(e$'$)} $g_\zeta\le f_{\alpha_\zeta, \zeta}$ and for
$\alpha<\delta$ we have $f_\alpha\le g_\zeta\mod D$

\noindent
3) By 1.8(1)
\sqed{3.11}

\proclaim{Claim 3.12}
\item{(1)}
Assume $\lambda=\tcf(\prod\bar\lambda/D)$ and $\mu=\cf(\mu)<\lambda$
$\underline{\rm then}$ there is $\bar\lambda'<_D\bar\lambda$,
$\bar\lambda'$ a
sequence of regular cardinals and $\mu=\tcf(\prod\bar\lambda'/D)$
provided that
$$
\mu>\reg_*(D),\ \min(\bar \lambda)> \reg^{\sigma^+}_*(D)\hbox{
whenever }\sigma< \reg_*(D)
\leqno(*)
$$

\item{(2)} Let $I^*$ be the ideal dual to $D$, and assume $(*)$
above.
If $(*)(\alpha)$ of 1.5 holds and $\mu$ is semi-normal (for
$(\bar\lambda, I^*)$) $\underline{\sl then}$ it is normal.
\endproclaim

\demo{Proof}
\subdemo{Case\ 1} $\mu<\lim\inf_D(\bar \lambda)$

We let
$$
\lambda'= \cases{ \mu & if $\mu< \lambda_i$\cr
                   1 & if $\mu \geq \lambda_i$\cr}
$$
and we are done.

\subdemo{Case\ 2} $\lim\inf_D(\bar \lambda)\geq \theta\geq \reg_*(D)$,
$\mu>\theta$, and $(\forall \sigma< \reg_*(D))[\reg^\sigma_*(D)<
\theta]$.

Let $\theta=: \reg_*(D)$. There is an unbounded $S\subseteq \mu$ and an
$(S, \theta)$-continuity system $\bar a$ (see 2.4). As $\prod\bar
\lambda/D$ has true cofinality $\lambda$, $\lambda> \mu$ clearly there
are $g^*_\alpha\in \prod\bar \lambda$ for $\alpha< \mu$ such that $\bar
g^*=\langle g^*_\alpha: \alpha<\mu\rangle$ obeys $\bar a$ (exists as
$\theta\leq \lim\inf_D(\bar \lambda)$).

Now if in claim 3.10(1) for $\bar g^*$ possibility (A) holds, we are done. By
3.11(1) we get that for some $\sigma< \reg_*(D)$ $\reg^\sigma_*(I)\geq
\mu$, contradiction.

\subdemo{Case\ 3}
$\lim\inf_D(\bar \lambda)\geq \theta$ $\reg_*(D)$, $\mu\geq \theta$, and
$(\forall \sigma< \reg_*(D))[\reg^\sigma_*(D)<\theta]$.

Like the proof of [Sh-g, Ch II 1.5B] using the silly square.
\medskip

\centerline{$\ast$\qquad$\ast$\qquad$\ast$}
\medskip

We turn to other measures of $\prod \bar \lambda /D$.

\rem{Definition 3.13}
\item{(a)} $T_D^0(\bar\lambda)=\sup\{ |F|\colon
F\subseteq\Prod\bar\lambda$ and $f_1\not= f_2\in F\Rightarrow
f_1\not=_D f_2\}$.

\item{(b)} $T^1_D(\bar\lambda)=\Min\{ |F|$: (i)
$F\subseteq\Prod\bar\lambda$

\hskip81pt(ii) $f_1\not= f_2\in F\Rightarrow f_1\not=_D f_2$

\hskip80pt(iii) $F$ maximal under (i)+(ii)$\}$

\item{(c)} $T^2_D(\bar\lambda)=
\Min\{ |F|\colon F\subseteq \Prod\bar\lambda$
and for every $f_1\in \Prod\bar\lambda$, for some
$f_2\in F$ we have $\neg f_1\not=_D f_2\}$.

\item{(d)} If
$T^0_D(\bar\lambda)=T^1_D(\bar\lambda)=T^2_D(\bar\lambda)$ then
let $T_D(\bar\lambda)=T^l_D(\bar\lambda)$ for $l<3$.

\item{(e)} for $f\in {}^\kappa\Ord$ and $\ell<3$ let $T^l_D(f)$ means
$T^l_D(\langle f(\alpha)\colon\alpha<\kappa\rangle)$.

\proclaim{Theorem 3.14}
\item{(0)} If $D_0\subseteq D_1$ are filters on $\kappa$ then
$T^{\ell}_{D_0}(\bar \lambda)\leq T^{\ell}_{D_1}(\bar \lambda)$
for $\ell=0,2$. Also if $\kappa=A_0\cup A_1$, $A_0\in D^+$, and $A_1\in
D^+$ then $T^\ell_D(\bar \lambda)=\min\{ T^\ell_{D+A_0}(\bar \lambda),
T^\ell_{D+A_1}(\bar \lambda)\}$ for $\ell=0, 2$.

\item{(1)} $\htcf_D(\Prod\bar\lambda)\le T^2_D(\bar\lambda)\le
T^1_D(\bar\lambda)\le T^0_D(\bar\lambda)$

\item{(2)} If $T^0_D(\bar\lambda)> |\calP(\kappa)/D|$ or just
$T^0_D(\lambda)>\mu$, and $\calP(\kappa)/D$ satisfies the $\mu^+$-c.c.
$\underline{\sl then}$
$T^0_D(\bar\lambda)=T^1_D(\bar\lambda)=T^2_D(\bar\lambda)$ so
the supremum in 3.13(a) is obtained (so e.g.
$T^0_D(\bar\lambda)>2^\kappa$ suffice)

\item{(3)} $T_D^0(\bar\lambda)^{<\reg D}=
T^0_D(\bar\lambda)$ (each $\lambda_i$ infinite of course).

\item{(4)}
$[\htcf_D\Prod_{i<\kappa}f(i)]\leq T^2_D(f)\le
[\htcf_D\Prod_{i<\kappa}f(i)]^{<\reg(D)}+
\wsat(D)^\kappa/D$

\item{(5)} If $D$ is an ultrafilter
$|\Prod\bar\lambda/D|=T^e_D(\bar\lambda)$ for $e\le 2$.

\item{(6)} In (4), if $\bigwedge_{i<\kappa} f(i)\ge 2^\kappa$ (or just
$(\wsat(D)+2)^\kappa/D\leq \min\limits_{i<\kappa} f(i)$), the second and
third terms are equal.

\item{(7)} If the sup in the definition of $T_D^0(\bar\lambda)$
is not obtained then it has cofinality $\ge\reg(D)$ and even is regular.
\endproclaim

\demo{Proof}
(0) Check.

(1) First assume
$\mu=:T^2_D(\bar\lambda)<\htcf_D(\prod\bar\lambda)$; then we can
find $\mu^* =\cf(\mu^*)\in (\mu, \htcf_D(\prod\bar\lambda)]$ and
$\bar\mu=\langle \mu_i\colon i<\kappa\rangle$, a sequence of
regular cardinals, $\bigwedge_{i<\kappa}\mu_i\le\lambda_i$ such that
$\mu^*=\tcf(\prod\bar\mu /D)$ and let $\langle f_\alpha\colon
\alpha<\mu^*\rangle$ exemplify this. Now let $F$ exemplify $\mu
= T^2_D(\bar\lambda)$, for each $g\in F$ let
$$
g'\in \prod\limits_{i<\kappa}\mu_i \hbox{ be }\colon g'(i)=
\cases{ g(i) & if $g(i)<\mu_i$\cr
0 & otherwise.\cr}
$$
So there is $\alpha(g)<\mu^*$ such that $g'<_D f_{\alpha(g)}$.
Let $\alpha^*=\sup\{\alpha(g)\colon g\in F\}$, now
$\alpha^*<\mu^*$ (as $\mu^*=\cf\mu^*>\mu=|F|)$.
So $g\in F\Rightarrow g\not=_D f_{\alpha^*}$, contradiction.
So really $T^2_D(\bar\lambda)\le\htcf_D(\prod\bar\lambda)$ as
required.

If $F$ exemplifies the value of $T^1_D(\bar\lambda)$, it also
exemplifies $T^2_D(\bar\lambda)\le |F|$ hence
$T^2_D(\bar\lambda)\le T^1_D(\bar\lambda)$.

Lastly if $F$ exemplifies the value of $T^1_D(f)$ it also
exemplifies $T^0_D(\bar\lambda)\ge |F|$, so
$T^1_D(\bar\lambda)\le T^0_D(\bar\lambda)$.

(2) Let $\mu$ be $|\calP(\kappa)/D|$ or at least $\mu$ is such that
the Boolean 
algebra $\calP(\kappa)/D$ satisfies the $\mu^+$-c.c.
Assume that the desired conclusion fails so $T^2_D(\bar\lambda)<
T^0_D(\bar\lambda)$, so there is $F_0\subseteq\prod\bar\lambda$, such that
$[f_1\not= f_2\in F_0\Rightarrow f_1\not=_D f_2]$, and
$|F_0|>T^2_D(\bar\lambda)+\mu$ (by the definition of $T^0_D(\bar\lambda)$).
Also there is $F_2\subseteq\prod\bar\lambda$ exemplifying the value
of $T^2_D(\bar\lambda)$.
For every $f\in F_0$ there is $g_f\in F_2$ such that
$\neg f\not=_D g_f$ (by the choice of $F_2$).
As $|F_0|> T^2_D(\bar\lambda) + \mu$ for some $g\in
F_2$, $F^*=:\{ f\in F_0\colon g_f= g\}$ has cardinality $>T^2_D(f)+\mu$.
Now for each $f\in F^*$ let $A_f=\{ i<\kappa\colon f(i)=g(i)\}$
clearly $A_f\in D^+$. Now $f\mapsto A_f/D$ is a function from $F^*$ into
$\calP(\kappa)/D$, hence, if $\mu=|\calP(\kappa)/D|$, it is not one to one
(by cardinality consideration) so for some $f'\not= f''$ from $F^*$ (hence form
$F_0$) we have $A_{f'} /D= A_{f''}/D$; but so
$$
\{ i<\kappa\colon f'(i) = f''(i)\}\supseteq \{ i<\kappa\colon
f'(i) =g(i)\} \cap \{ i<\kappa\colon f''(i) = g(i)\}
= A_{f'}/D
$$
hence is $\not=\emptyset\mod D$, so $\neg f'\not=_D f''$,
contradition the choice of $F_0$.
If $\mu\neq |\calP(\kappa)/D|$ (as $F^*\subseteq F_0$ by the choice of $F_0$)
we have:
$$f_1\not= f_2 \in F^* \Rightarrow  A_{f_1}\cap A_{f_2}= \emptyset\ \mod\ D$$
so $\{A_f: f\in F^*\}$ contradicts ``the Boolean algebra $\calP(\kappa)/D$
satisfies the $\mu^+$-c.c.".

(3) Assume that $\theta<\reg(D)$ and\footnote\dag{$\leq^+$ means the left
side is a supremum, right bigger than the left or equal but the supremum
is obtained}
$\mu\le^+ T^0_D(\bar\lambda)$. As $\mu\le^+ T^0_D(\bar\lambda)$ we can find
$f_\alpha\in\prod\bar\lambda$ for $\alpha<\mu$ such that
$[\alpha<\beta\Rightarrow f_\alpha\not=_D f_\beta]$. Also (as $\theta<
\reg(D)$) we can find $\{ A_\varepsilon\colon\varepsilon<\theta\}\subseteq D$
such that for every $i<\kappa$ the set $w_i=: \{\varepsilon<\theta\colon i\in
A_\varepsilon\}$ is finite. Now for every function $h\colon\theta\to\mu$ we
define $g_h$, a function with domain $\kappa$:
$$
g_h(i) = \{ (\varepsilon, f_{h(\varepsilon)}(i))\colon \varepsilon\in w_i\}
$$
So $|\{ g_h(i)\colon h\transa{\in}{\mu}\} |\le (\lambda_i)^{|w_i|}=\lambda_i$,
and if $h_1\not= h_2$ are $\transa{\rm from}{\mu}$ then for some
$\varepsilon<\theta$, $h_1(\varepsilon)\not= h_2(\varepsilon)$ so
$B_{h_1, h_2}= \{ i\colon  f_{h_1(\varepsilon)} (i)\not=
f_{h_2(\varepsilon)}(i)\}\in D$ that is
$B_{h_1, h_2}\cap A_\varepsilon\in
D$ so

$\otimes_1$ if $i\in B_{h_1, h_2}\cap A_\varepsilon$ then $\varepsilon
\in w_i$, so $g_{h_1}(i)\not= g_{h_2}(i)$.

$\otimes_2$ $B_{h_1, h_2}\cap A_\varepsilon \in D$

\noindent So $\langle g_h\colon h \transa{\in}{\mu}\rangle$ exemplifies
$T^0_D(\bar\lambda)\ge \mu^\theta$. If the supremum in the definition of
$T^0_D(\bar\lambda)$ is obtained we are done. If not then $T_D^0(\bar\lambda)$
is a limit cardinal, and by the proof above:
$$[\mu<T^0_D(\bar\lambda)\quad \&\quad
\theta<\reg(D)\quad\Rightarrow\quad \mu^\theta<T^0_D(\bar\lambda)].$$
So if $T^0_D(\bar\lambda)$ has cofinality $\ge\reg(D)$ we are done;
otherwise let it be $\sum_{\varepsilon<\theta}\mu_\varepsilon$ with
$\mu_\varepsilon<T^0_D(\bar\lambda)$ and $\theta<\reg D$. Note that by the
previous sentence $T^0_D(\bar\lambda)^\theta =
T^0_D(\bar\lambda)^{<\reg(D)} = \prod\limits_{\varepsilon<\theta}
\mu_\varepsilon$, and let $\{ f_\alpha^\varepsilon\colon \alpha<
\mu_\varepsilon\}\subseteq\prod\bar\lambda$ be such that $[\alpha<
\beta\Rightarrow f_\alpha^\varepsilon\not=_D f_D^\varepsilon]$ and repeat the
previous proof with $f^\varepsilon_{h(\varepsilon)}$ replacing
$f_{h(\varepsilon)}$.

\subdemo{(4) For the first inequality} assume it fails so $\mu=: T^2_D(f)
< \htcf_D(\prod\limits_{i<\kappa} f(i))$ hence for some $g\in \prod\limits_{i<
f(i)} (f(i)+1)$, $\tcf(\prod\limits_{i<\kappa} g(i), <_D)$ is $\lambda$ with
$\lambda=\cf(\lambda)>\mu$. Let $\langle f_\alpha:
\alpha<\lambda\rangle$ exemplifies this. Let $F$ be as in the definition
of $T^2_D(f)$, now for each $h\in F$, there is $\alpha(h)<\lambda$ such that
$$\{i< \kappa:\hbox{ if }h(i)< g(i)\hbox{ then }h(i)<
f_{\alpha(g)}(i)\}\in D.$$
Let $\alpha^*=\sup\{\alpha(h)+1: h\in F\}$, now $f_{\alpha^*}\in
\prod\limits_{i<\kappa} f(i)$ and $h\in F\Rightarrow h\not=_D
f_{\alpha^*}$ contradicting the choice of $F$.

\subdemo{for the second inequality} Repeat the proof of 3.9 except that
here we prove $F=: \bigcup\limits_{\zeta<\theta} (N_\zeta \cap
\prod\limits_{i<\kappa} f(i))$ exemplifies $T^2_D(f)\leq \lambda$;
we replace clause (g) in the proof  by
\item{(g)$'$} $g^*< f^2_{\zeta+1}< f^1_\zeta\mod D$

\noindent the construction is for $\zeta<\reg(D)$ and if we find
satisfy $\neg f^1_\zeta\not=_D g^*$ we are done.

(5) Straightforward.

(6) Note that all those cardinals are $\geq 2^\kappa$ and $2^\kappa \geq
\wsat(D)^\kappa/D$. Now write successively inequalities from (2), (4),
(1) and (3):
$$
T^0_D(f)= T^2_D(f) \leq [\htcf_D \prod\limits_{i<\kappa} f(i)]^{<\reg(D)}
\leq [T^0_D(f)]^{<\reg(D)}= T^0_D(f).
$$

(7) See proof of part (3). Moreover, if
$\mu=\sum\limits_{\varepsilon<\tau} \mu_\varepsilon$, $\tau<
T^0_D(\bar \lambda)$, $\mu_\varepsilon< T^0_D(\bar \lambda)$ as
exemplified by $\{f_\varepsilon: \varepsilon<\tau\}$,
$\{f^\varepsilon_\alpha: \alpha< \mu_\varepsilon\}$ respectively. Let
$g_\alpha$ 
be: if $\sum\limits_{\varepsilon <\zeta} \mu_\varepsilon < \alpha <
\sum\limits_{\varepsilon\leq \zeta} \mu_\varepsilon$ then
$g_\alpha(i)= (f_\varepsilon (i), f^\varepsilon_\alpha(i))$. So
$\{g_\alpha: \alpha< \mu\}$ show: if $T^0_D(\bar \lambda)$ is singular
then the supremum is obtained. 
\sqed{3.14}

\proclaim{Claim 3.15}
Assume $D$ is a filter on $\kappa$, $f\trans\in\Ord$, $\mu^{\aleph_0}=\mu$
 and $2^\kappa<\mu$, $T_D(f)$,
(see Definition 3.13(d) and Theorem 3.14(2)).
If $\mu< T_D(f)$ $\underline{\rm then}$ for some sequence
$\bar\lambda\le f$ of regulars,
$\mu^+=\tcf(\prod\bar\lambda /D)$, or at least

\item{$(*)$} there are $\langle\langle\lambda_{i,n}\colon n<n_i\rangle
\colon i<\kappa\rangle$,
$\lambda_{i, n}=\cf(\lambda_{i, n})< f(i)$ and a filter $D^*$ on
$\bigcup_{i<\kappa}\{ i\}\times n_i$ such that: $\mu^+=\tcf
(\prod\limits_{(i, n)}\lambda_{i, n} /D^*)$ and $D=\{
A\subseteq\kappa\colon\bigcup_{i\in A}\{i\}\times n_i\in D^*\}$.

Also the inverse is true.
\endproclaim

\rem{Remark 3.15A}
(1) It is not clear whether the first possibility may fail. We have
explained earlier the doubtful role of $\mu^{\aleph_0}= \mu$.

(2) We can replace $\mu^+$ by any regular $\mu$ such that
$\bigwedge_{\alpha<\mu}|{\alpha}|^{\aleph_0}<\mu$ and then we use
3.14(4) to get $\mu\le^+ T_D(f)$.

(3) The assumption $2^\kappa< \mu$ can be omitted.

\demo{Proof}
The inverse should be clear (as in the proof of 3.6, by 3.14(3)).

${\rm Wlog}\, f(i)>2^\kappa$ for $i<\kappa$, and trivially
$(\wsat(D))^\kappa/D\le 2^\kappa$, so by 3.14(4)
$$
T_D(f)\le
[\htcf_D(\prod\limits_{i<\kappa}f(i)]^{<\reg(D)}.
$$
If $\mu<\htcf_D(\prod\limits_{i<\kappa} f(i))$ we are done (by 3.12(1)),
so assume
$\htcf_D(\prod\limits_{i<\kappa} f(i))\le\mu$, but we have assumed
$\mu< T_D(f)$ so we can conclude $\mu^{<\reg(D)}\ge\mu^+$.
Let $\chi\le\mu$ be minimal such that
$\bigvee_{\theta<\reg(D)}\chi^\theta
\ge\mu$, and let $\theta=:\cf(\chi)$ so, as $\mu>2^\kappa$
we know $\chi^{\cf\chi}=
\chi^{<\reg(D)}=\mu^{<\reg(D)}\ge \mu^+$, $\chi>2^\kappa$,
$\bigwedge_{\alpha<\chi}|\alpha|^{<\reg(D)}<\chi$.
By the assumption $\mu=\mu^{\aleph_0}$ we know $\theta>\aleph_0$
(of course $\theta$ is regular).
By [Sh-g, VIII 1.6(2), IX 3.5] and [Sh513, 6.12] there is
a strictly increasing sequence
$\langle\mu_\varepsilon\colon\varepsilon<\theta\rangle$ of regular
cardinals with limit $\chi$ such that
$\mu^+=\tcf(\prod\limits_{\varepsilon<\theta}\mu_\varepsilon
/J_\theta^{\bd})$.

As clearly $\chi\le\htcf_D(\prod\limits_{i<\kappa}f(i))$, we can find
for each $\varepsilon<\theta$,
a sequence $\bar\lambda^\varepsilon=\langle
\lambda_i^\varepsilon\colon i<\kappa\rangle$ such that
$\lambda_i^\varepsilon=\cf(\lambda_i^\varepsilon)\leq f(i)$, and
$\tcf(\prod\limits_{i<\kappa}\lambda_i^\varepsilon/D)=\mu_\varepsilon$,
also $\wlog\lambda_i^\varepsilon>2^\kappa$.
Let $\langle A_\varepsilon\colon\varepsilon<\theta\rangle$ exemplify
$\theta<\reg(D)$ and $n_i=|\{\varepsilon<\theta\colon i\in A_\varepsilon
\}|$ and $\{\lambda_{i,n}\colon n<\omega\}$ enumerate
$\{\lambda_i^\varepsilon\colon \varepsilon$ satisfies $i\in A_\varepsilon\}$,
so we have gotten $(*)$.
\sqed{3.15}

\rem{Conclusion 3.16}
Suppose $D$ is an $\aleph_1$-complete filter on $\kappa$.
If $\lambda_i\ge 2^\kappa$ for $i<\kappa$ and $\sup_{A\in D^+}
T_{D+A}(\bar\lambda)>\mu^{\aleph_0}$ $\underline{\rm then}$ for
some $\lambda'_i=\cf(\lambda'_i)\le \lambda_i$ we have
$$
\sup_{A\in D^+}\htcf_{D+A}(\prod\limits_{i<\kappa}\lambda'_i)>\mu.
$$

\rem{Conclusion 3.17}
Let $D$ be an $\aleph_1$-complete filter on $\kappa$.
If for $i<\kappa$, $B_i$ is a Boolean algebra and $\lambda_i<
\Depth^+(B_i)$ (see below) and
$$
2^\kappa<\mu^{\aleph_0}<\sup_{A\in D^+} T_{D+A}(\bar\lambda)
$$
$\underline{\rm then}$ $\mu^+<\Depth^+(\prod\limits_{i<\kappa} B_i/D)$.

\demo{Proof}
Use 3.25 below and 3.16 above.

\rem{Definition 3.18}
For a partial order $P$ (e.g.\ a Boolean algebra) let
$\Depth^+(P)=\Min\{\lambda:$ we cannot find $a_\alpha\in P$ for
$\alpha<\lambda$ such that $\alpha<\beta\Rightarrow a_\alpha<_P
a_\beta\}$.

\rem{Discussion 3.19}
\item{(1)} We conjecture that in 3.16 (and 3.17) the assumption
``$D$ is $\aleph_1$-complete'' can be omitted.

\item{(2)} Note that our results are for $\mu=\mu^{\aleph_0}$
only; to remove this we need to improve the theorem on $pp=cov$
(i.e.\ to prove $\cf(\lambda)=\aleph_0<\lambda\Rightarrow
pp(\lambda)=cov(\lambda, \lambda, \aleph_1, 2)$ (or $\sup\{
pp(\mu)\colon\cf\mu=\aleph_0<\mu<\lambda\}=\cf(S_{\le\aleph_0}(\lambda),
\subseteq)$ (see [Sh-g], [Sh430, \S1]), which seems to me a very
serious open problem (see [Sh-g, Analitic guide, 14]).

\item{(3)} In 3.17, if we can find
$f_\alpha\in\prod\limits_{i<\kappa}\lambda_i$ for
$\alpha<\lambda\colon[\alpha<\beta<\lambda\Rightarrow
f_\alpha\le f_\beta\mod D]$ and $\neg f_\alpha =_D f_{\alpha+1}$
then $\lambda<\Depth^+(\prod\limits_{i<\kappa} B_i/D)$.
But this does not help for $\lambda$ regular $> 2^\kappa$.

\item{(4)}
We can approach 3.15 differently, by 3.20--3.23 below.

\proclaim{Claim 3.20}
If $2^{2^\kappa}\le \mu<T_D({\bar\lambda})$, (or at least
$2^{|D|+\kappa}\le \mu<T_D(\bar\lambda)$) and
$\mu^{<\theta}=\mu$, $\underline{\rm then}$ for some
$\theta$-complete filter $E\subseteq D$ we have $T_E(\bar\lambda)>\mu$.
\endproclaim

\demo{Proof}
Wlog $\theta$ is regular (as $\mu^{<\theta}=\mu\ \&$,
$\cf(\theta)<\theta\Rightarrow \mu^{<\theta^+}=\mu$).
Let $\{ f_\alpha\colon\alpha<\mu^+\}\subseteq\prod{\bar\lambda}$,
be such that $[\alpha<\beta\Rightarrow f_\alpha \not=_D f_\beta]$.
We choose by induction on $\zeta$, $\alpha_\zeta<\mu^+$ as follows:
$\alpha_\zeta$ is the minimal ordinal $\alpha<\mu^+$ such that
$E_{\zeta,\alpha}\subseteq D$ where $E_{\zeta,\alpha}=$ the
$\theta$-complete filter generated by
$$
\big\{\{i<\kappa\colon f_{\alpha_\varepsilon}(i)\not=
f_\alpha(i)\}\colon \varepsilon<\zeta\big\}
$$
(note: each generator of $E_{\zeta,\alpha}$ is in $D$ but not
necessarly $E_{\zeta,\alpha}\subseteq D$!).

\noindent
Let $\alpha_\zeta$ be well defined if $\zeta<\zeta^*$, clearly
$\varepsilon<\zeta\Rightarrow \alpha_\varepsilon<\alpha_\zeta$.
Now if $\zeta^*<\mu^+$, then clearly
$\alpha^*=\bigcup_{\zeta<\zeta^*}\alpha_\zeta<\mu^+$ and for every
$\alpha\in(\alpha^*, \mu^+)$, $E_{\zeta^*, \alpha}
\not\subseteq D$, so
for every such $\alpha$ there are $A_\alpha\in D^+$ and
$a_\alpha\in [\zeta^*]^{<\theta}$ such
that $A_\alpha=\bigcup_{\varepsilon\in a_\alpha}\{ i<\kappa\colon
f_{\alpha_\varepsilon}(i)=f_\alpha(i)\}$.
But for every $A\in D^+$, $a\in [\zeta^*]^{<\theta}$ we have
$$
\{\alpha\colon\alpha\in(\alpha^*,\mu^+), A_\alpha=A,
a_\alpha=a\}\subseteq
\{\alpha\colon f_\alpha\upharpoonright A\in
\Prod_{i<\kappa}\{ f_{\alpha_\varepsilon}(i)
\colon \varepsilon\in a_\alpha\}\},
$$
hence has cardinality $\le\theta^\kappa\le 2^{\kappa}<\mu$.
Also $|[\zeta^*]^{<\theta}|\le\mu^{<\theta}<\mu^+$,
$\vert D^+\vert\le 2^\kappa <\mu^\kappa$
so we get easy contradiction.

So $\zeta^*=\mu^+$, but the number of possible $E$'s is
$\le 2^{2^\kappa}$, hence for some $E$ we have $|\{ \varepsilon<\mu^+\colon
E_{\varepsilon, \alpha_\varepsilon}=E\}|=\mu^+$.
Necessarily $E\subseteq D$ and $E$ is $\theta$-complete,
and $\{ f_{\alpha_\varepsilon}\colon\varepsilon<\mu^+$, and
$E_{\alpha_\varepsilon}=E\}$ exemplifies $T_E(\bar\lambda)>\mu$,
so $E$ is as required.
\sqed{3.20}

\rem{Fact 3.21}
1. In 3.20 we can replace $\mu^+$ by $\mu^*$ if $2^{2^\kappa}<\cf
(\mu^*)\le\mu^*\le T_D^0({\bar\lambda})$ and
$\bigwedge_{\alpha<\mu^*}|\alpha|^{<\theta}<\mu^*$.

2.  We can, in 3.20, [and 3.21(1)] replace
``$T_D(\bar\lambda)>\mu$'' by
``$\prod\bar\lambda /D$ has an increasing
sequence of lengths $>\mu[\ge\mu]$'', we can deduce this also otherwise.

\proclaim{Claim 3.22}
\item{(1)} If $2^\kappa<|\prod{\bar\lambda}/D|$, $D$ an ultrafilter on
$\kappa$, $\mu=\cf(\mu)\le|\prod\bar\lambda/ D|$,
$\bigwedge_{i<\kappa} |i|^{\aleph_0}<\mu$, and $D$ is regular
$\underline{\rm then}$ $\mu<\Depth^+(\prod\limits_{i<\kappa}\lambda_i/ D)$

\item{(2)} Similarly for $D$ just a filter.
\endproclaim

\demo{Proof}
${\rm Wlog}\,\lambda=\lim_D{\bar\lambda}=\sup({\bar\lambda})$,
so $|\prod{\bar\lambda}/D|=\lambda^\kappa$ (by [CK]).
If $\mu\le\lambda$ we are done;
otherwise let $\chi=\Min\{\chi\colon\chi^\kappa=
\lambda^\kappa\}$, so $\chi^{\cf(\chi)}=\lambda^\kappa$,
$\cf(\chi)\le\kappa$ but $\lambda<\mu\le\lambda^\kappa$ hence
$\lambda^{\aleph_0}<\mu$ hence $\cf(\chi)>\aleph_0$, also by $\chi's$
minimality $\bigwedge_{i<\chi}|i|^{\cf\chi}\leq |i|^\kappa<\chi$, and
remember $\chi<\mu=\cf\mu\le \chi^{\cf\chi}$ so by [Sh-g.\ VIII 1.6(2)]
there is $\langle\mu_\varepsilon\colon \varepsilon<\cf(\chi)\rangle$
strictly increasing sequence of regular cardinals with limit $\chi$,
$\prod\limits_{\varepsilon< \cf(\chi)}
\mu_\varepsilon/J_{\cf\chi}^{bd}$ has true
confinality $\mu$.
Let $\chi_\varepsilon=\sup\{
\mu_\zeta\colon\zeta<\varepsilon\}+2^\kappa$, let
$\newagei\colon\kappa\to\cf(\chi)$ be
$\newagei(i)=\sup\{\varepsilon+1\colon
\lambda_i\ge\chi_\varepsilon\}$.
If there is a function $h\in\prod\limits_{i<\kappa}\newagei(i)$
such that $\bigwedge_{j<\cf(\chi)}
\{i<\kappa\colon h(i)<j\}=\emptyset\mod D$
then $\prod\limits_{i<\kappa}\mu_{h(i)}/D$ has true cofinality
$\mu$ as required; if not $(D, \newagei)$ is weakly normal (i.e. there
is no such $h$ - see [Sh420]).
But for $D$ regular, $D$ is $\cf(\chi)$-regular, some $\langle
A_\varepsilon\colon\varepsilon<\cf(\kappa)\rangle$ exemplifies it and
$h(i)=\max\{ \varepsilon\colon\varepsilon< \newagei(i)$ and $i\in
A_\varepsilon\}$ (maximum over a finite set) is as required.
\sqed{3.22}

\rem{Discussion 3.23}
\item{1.} In 3.20 (or 3.21) we can apply [Sh 410, \S6] so
$\mu=\tcf(\prod \bigcup_{i<\mu}\fra_i/D^*$, where
$D=\{ A\subseteq\kappa\colon\bigcup_{i\in A}\fra_i\in D^*\}$ and each
$\fra_i$ is finite.

In 3.15 we have gotten this also for $\mu\in(2^\kappa, 2^{2^\kappa})$.

\proclaim{Claim 3.24}
If $D$ is a filter on $\kappa$, $B_i$ is the interval Boolean
algebra on the ordinal $\alpha_i$, and $|\prod\limits_{i<\kappa}
\alpha_i/D|>2^\kappa$ $\underline{\rm then}$ for regular $\mu$
we have $\mu<\Depth^+(\prod\limits_{i<\kappa} B_i/D)$ iff for some
$\mu_i\le\alpha_i$ (for $i<\kappa$) and $A\in D^+$, the true
cofinality of $\prod\limits_{i<\kappa}\mu_i /(D+A)$) is well defined and
equal to $\mu$.
\endproclaim

\demo{Proof}
The $\Rightarrow$ (i.e. only if direction) is clear. For the $\Leftarrow$
direction assume $\mu$ is regular $< \Depth^+(\prod\limits_{i<\kappa}
B_i/D)$ so there are $f_\alpha\in \prod\limits_{i<\kappa} B_i$ such that
$\prod\limits_{i<\kappa} B_i/D \vDash f_\alpha/D < f_\beta /D$ for
$\alpha<\beta$.

Wlog $\mu>2^\kappa$. Let $f_\alpha(i)= \bigcup\limits_{\ell< n(\alpha,
i)}[j_{\alpha, i, 2\ell}, j_{\alpha, i. 2\ell+1})$ where $j_{\alpha, i,
\ell}< j_{\alpha, i, \ell+1}< \alpha_i$ for $\ell< 2n(\alpha, i)$.
As $\mu=\cf(\mu)> 2^\kappa$
wlog $n_{\alpha, i}= n_{i}$. By [Sh430, 6.6D] (see more [Sh513, 6.1]) we
can find $A\subseteq A^* =: \{(i, \ell): i<\kappa, \ell< 2n_\alpha\}$
and $\langle \gamma^*_{i, \ell} : i< \kappa, \ell< 2n_i\rangle$ such
that $(i, \ell)\in A \Rightarrow \gamma^*_{i, \ell}$ is a limit ordinal
and

\item{$(*)$} for every $f\in \prod\limits_{(i, \ell)\in A} \gamma^*_{i, \ell}$
and $\alpha< \mu$ there is $\beta\in (\alpha, \mu)$ such that

\itemitem{} $(i, \ell)\in A^* \setminus A \Rightarrow j_{\alpha, i,
\ell} = \gamma^*_{i, \ell}$

\itemitem{} $(i, \ell) \in A \Rightarrow f(i, \ell)< j_{\alpha, i, \ell}
< \gamma^*_{i, \ell}$

\itemitem{} $(i, \ell) \in A \Rightarrow \cf(\gamma^*_{i, \ell})>
2^\kappa$

Let $\ell(i)= \max\{\ell< 2n(i): (i, \ell)\in A\}$ and let  $B= \{i: \ell(i)$
well defined$\}$. Clearly $B\in D^+$ (otherwise we can find $\alpha<
\beta< \mu$ such that $f_\alpha/D= f_\beta/D$, contradiction). For $(i,
\ell)\in A$ define $\beta^*_{i, \ell}$ by
$\beta^*_{i, \ell} = \sup \{\gamma^*_{j, m}+1: (j, m)\in A^*$ and
$\gamma^*_{j, m}< \gamma^*_{i, \ell}\}$. Now $\beta^*_{i, \ell}<
\gamma^*_{i, \ell}$ as $\cf(\gamma^*_{i, \ell})> 2^\kappa$. Let
$$
\eqalign{
Y=\{\alpha<\mu: & \hbox{ if } (i, \ell)\in A^*\setminus A \hbox{ then
}j_{\alpha, i, \ell}= \gamma^*_{i, \ell} \cr
 & \hbox{ and if } (i, \ell)\in A \hbox{ then } \beta^*_{i, \ell} <
j_{\alpha, \ell, i} < \gamma^*_{\ell, i}\}\cr}
$$
Let $B_1=\{ i\in B: \ell(i) \hbox{ is odd}\}$. Clearly $B_1 \subseteq
B$ and $B\setminus B_1 =\emptyset\ \mod\ D$ (otherwise as in $(*)_1$,
$(*)_2$ below get contradiction)
hence $B_1\in D^+$. Now
\item{$(*)_1$} for $\alpha<\beta$ from $Y$ we have
$$
\langle j_{\alpha, i, \ell(i)}: i\in B_1\rangle \leq \langle j_{\beta, i,
\ell(i)}: i\in B_1\rangle\ \mod\ (D\restriction B_1)
$$
[Why? as $f_\alpha/D$ was non decreasing in $\prod\limits_{i<\kappa} B_i/D$]

\item{$(*)_2$} for every $\alpha\in Y$ for some $\beta$, $\alpha<
\beta\in Y$ we have
$$
\langle j_{\alpha, i, \ell(i)}: i\in B_1\rangle < \langle j_{\beta, i,
\ell(i)}: i\in B_1\rangle\ \mod\ (D\restriction B_1)
$$
[Why? by $(*)$ above]

Together for some unbounded $Z\subseteq Y$, $\big\langle \langle
j_{\alpha, \ell, \ell(i)}: i\in B_1\rangle / (D\restriction B_1):
\alpha\in Z\big\rangle$ is $<_{D\restriction B_1}$-increasing, so it has
a $<_{(D\restriction B_1)}-\eub$ (as $\mu> 2^\kappa$), say $\langle j^*_i:
i\in B_1\rangle$ hence $\prod\limits_{i\in B_1} j^*_i / (D\restriction
B_1)$ has true cofinality $\mu$, and clearly $j^*_i \leq \gamma^*_{i,
\ell(i)} \leq \alpha_i$, so we have finished.
\sqed{3.24}

\proclaim{Claim 3.25}
If $D$ is a filter on $\kappa$, $B_i$ a Boolean algebra,
$\lambda_i<\Depth^+(B_i)$ $\underline{\rm then}$
\item{(a)} $\Depth(\prod\limits_{i<\kappa} B_i/D)\ge\sup_{A\in D^+}\tcf(\prod
_{i<\kappa}\lambda_i /(D+A))$ (i.e.\ on the cases $\tcf$ is well defined).

\item{(b)} $\Depth^+(\prod\limits_{i<\kappa} B_i/D)$ is
$\ge\Depth^+(\calP(\kappa)/D)$ and is at least
$$
\sup\{[\tcf(\prod\limits_{i<\kappa}\lambda'_i/(D+A))]^+\colon\lambda'_i
<\Depth^+(B_i), A\in D^+\}.
$$
\endproclaim

\demo{Proof} Check.

\proclaim{Claim 3.26}
Let $D$ be a filter on $\kappa$, $\langle \lambda_i: i<\kappa\rangle$ a
sequence of cardinals and $2^\kappa< \mu=\cf(\mu)$. Then $(\alpha)
\Leftrightarrow (\beta) \Rightarrow (\gamma) \Rightarrow (\delta)$, and
if $(\forall \sigma< \mu) (\sigma^{\aleph_0}<\mu)$ we also have
$(\gamma) \Leftrightarrow (\delta)$ where

\item{$(\alpha)$} if $B_i$ is a Boolean algebra, $\lambda_i<
\Depth^+(B_i)$ $\underline{\rm then}$ $\mu<
\Depth^+(\prod\limits_{i<\kappa} B_i/D)$

\item{$(\beta)$} there are $\mu_i=\cf(\mu_i)\leq \lambda_i$ for
$i<\kappa$ and $A\in D^+$ such that $\mu= \tcf(\prod\mu_i / (D+A))$

\item{$(\gamma)$}
there are $\langle\langle\lambda_{i,n}\colon n<n_i\rangle
\colon i<\kappa\rangle$,
$\lambda_{i, n}=\cf(\lambda_{i, n})< \lambda_i$ and a filter $D^*$ on
$\bigcup_{i<\kappa}\{ i\}\times n_i$ such that:
$$
\mu=\tcf
(\prod\limits_{(i, n)}\lambda_{i, n} /D^*)\hbox{ and }
D=\{A\subseteq\kappa\colon\hbox{ the set }\bigcup_{i\in A}\{i\}\times
n_i\hbox{ belongs to } D^*\}.
$$

\item{$(\delta)$} for some $A\in D^+$, $\mu\leq T_{D+A}(\langle
\lambda_i: i<\kappa\rangle)$

\endproclaim

\rem{Remark}
So the question whether $(\alpha) \Leftrightarrow (\delta)$ assuming
 $(\forall \sigma< \mu) (\sigma^{\aleph_0}<\mu)$
 is equivalent to $(\beta) \leftrightarrow (\gamma)$ which is a ``pure"
$\pcf$ problem.

\demo{Proof}
Note $(\gamma) \Rightarrow (\delta)$ is easy (as in 3.15, i.e. as in the
proof of 3.6, only easier). Now $(\beta) \Rightarrow (\gamma)$ is trivial and $(\beta)
\Rightarrow (\alpha)$ by 3.25. Next $(\alpha) \Rightarrow (\beta)$ holds
as we can use $(\alpha)$ for $B_i=:$ the interval Boolean algebra of the
order $\lambda_i$ and use 3.24. Lastly assume
 $(\forall \sigma< \mu) (\sigma^{\aleph_0}<\mu)$, now $(\gamma)
\Leftrightarrow (\delta)$ by 3.15.
\sqed{3.26}

\rem{Discussion}
We would like to have (letting $B_i$ denote Boolean algebra)
$$
\Depth^{(+)}(\prod\limits_{i< \kappa} B_i /D)\geq \prod\limits_{i<\kappa}
\Depth^{(+)} (B_i)/D
$$
if $D$ is just filter we should use $T_D$ and so by the problem of
attainment (serious by Magidor Shelah [MgSh433]), we ask
\item{$\otimes$} for $D$ an ultrafilter on $\kappa$, does $\lambda_i<
\Depth^+(B_i)$ for $i<\kappa$ implies
$$
\prod\limits_{i<\kappa}\lambda_i/D < \Depth^+(\prod\limits_{i<\kappa} B_i / D)
$$
at least when $\lambda_i > 2^\kappa$;

\item{$\otimes'$} for $D$ a filter on $\kappa$, does $\lambda_i<
\Depth^+(B_i)$ for $i<\kappa$ implies, assuming $\lambda_i> 2^\kappa$
for simplicity,
$$
T_D (\langle \lambda_i: i<\kappa\rangle) < \Depth^+ (\prod\limits_{i<\kappa}
B_i/D)
$$
As explained in 3.26 this is a $\pcf$ problem.

However changing the invariant (closing under homomorphisms, see [M])
we get a nice result; this will be presented in [Sh580]. 

\head \S 4 Remarks on the conditions for the $\pcf$ analysis
\endhead

We consider a generalization whose interest is not so clear.

\proclaim{Claim 4.1}
Suppose $\bar\lambda=\langle \lambda_i\colon i<\kappa\rangle$ is a
sequence of regular cardinals, and $\theta$ is a cardinal and $I^*$ is
an ideal on
$\kappa$; and $H$ is a function with domain $\kappa$. We cosnsider the
following statements:
{\itemindent30pt
\item{$(**)_H$}
$\lim\inf_{I^*}(\bar\lambda)\ge\theta\ge\wsat(I^*)$ and $H$ is a
function from $\kappa$ to $\calP(\theta)$ such that:}
$$
\eqalign{
\hbox{(a)} & \hbox{ for every } \varepsilon<\theta\hbox{ we have }\quad
\{ i<\kappa\colon\varepsilon\in H(i)\}=\kappa\mod I^*\cr
\hbox{(b)} & \hbox{ for } i<\kappa \hbox{ we have }
\otp(H(i))\le\lambda_i \hbox{ or at least } \{ i<\kappa\colon
|H(i)|\ge\lambda_i\}\in I^*\cr}
$$
\item{$(**)^+$} similarly but
\item{(b)$^+$} for $i<\kappa$ we have $\otp(H(i)) < \lambda_i$

\item{(1)} In 1.5 we can replace the assumption $(*)$ by $(**)_H$
above.

\item{(2)} Also in 1.6, 1.7, 1.8, 1.9, 1.10, 1.11 we can replace
1.5$(*)$ by $(**)_H$.

\item{(3)} Suppose in Definition 2.3(2) we say $\bar f$ obeys
$\bar a$ for $H$ (instead of for $\bar A^*$) if

\itemitem{(i)} for $\beta\in a_\alpha$ such that
$\varepsilon=:\otp(a_\alpha)<\theta$ we have
$$
\otp(a_\beta), \otp(a_\alpha)\in H(i)\Rightarrow f_\beta(i)\le
f_\alpha(i)
$$
and in 2.3(2A), $f_\alpha(i)=\sup\{ f_\beta(i)\colon\beta\in
a_\alpha$ and $\otp(a_\beta), \otp(a_\alpha)\in H(i)\}$.

Then we can replace 1.5$(*)$ by $(**)_H$ in 2.5, 2.5A, 2,6; and
replace 1.5$(*)$ by $(**)^+_H$ in 2.7 (with the natural changes).
\endproclaim

\demo{Proof}
(1) Like the proof of 1.5, but defining the $g_\varepsilon$'s by
induction on $\varepsilon$ we change requirement (ii) to
\item{(ii)$'$} if $\zeta<\varepsilon$, and $i\in H(\zeta)\cap
H(\varepsilon)$ then $g_\zeta(i)<g_\varepsilon(i)$.

We can not succeded as
$$
\langle (B^\varepsilon_{\alpha(*)}\setminus
B_{\alpha(*)}^{\varepsilon+1})\cap\{ i<\kappa\colon \varepsilon,
\varepsilon+1\in H(i)\}\colon \varepsilon<\theta\rangle
$$
is a sequence of $\theta$ pairwise disjoint member of $(I^*)^+$.

In the induction, for $\varepsilon$ limit let
$g_\varepsilon(i)<\cup\{ g_\zeta(i)\colon\zeta\in H(i)$ and
$\varepsilon\in H(i)\}$ (so this is a union at most
$\otp(H(i)\cap\varepsilon)$ but only when $\varepsilon\in H(i)$
hence is $\langle\otp(H(i))\le\lambda_i$).

(2) The proof of 1.6 is the same, in the proof of 1.7 we again
replace (ii) by (ii)$'$.
Also the proof of the rest is the same.

(3) Left to the reader.
\sqed{4.1}

We want to see how much weakening $(*)$ of 1.5 to
``$\lim\inf_{I^*}(\bar\lambda)\ge\theta\ge\wsat(I^*)$ suffices.
If $\theta$ singular or $\lim\inf_{I^*}(\bar \lambda)>\theta$ or just $(\prod \bar
\lambda, <_{I^*})$ is $\theta^+$-directed then case $(\beta)$
of 1.5 applies.
This explains $(*)$ of 4.2 below.

\proclaim{Claim 4.2}
Suppose $\bar\lambda=\langle \lambda_i\colon i<\kappa\rangle$,
$\lambda_i=\cf(\lambda_i)$, $I^*$ an ideal on $\kappa$, and
$$
\lim\inf\nolimits_I(\bar\lambda)=\theta\ge\wsat(I^*),\quad
\theta\hbox{ regular }
\leqno(*)
$$

$\underline{\rm Then}$ we can define  a sequence $\bar J=\langle
J_\zeta\colon\zeta<\zeta(*)\rangle$ and an ordinal $\zeta(*)\le \theta^+$
such that
\item{(a)} $\bar J$ is an increasing continuous sequence of ideals on
$\kappa$.

\item{(b)} $J_0=I^*$, $J_{\zeta+1}=:\{ A\colon A\subseteq\kappa,$
and: $A\in J_\zeta$ or we can find $h\colon A\to\theta$ such that
$\lambda_i>h(i)$ and $\varepsilon<\theta\Rightarrow\{ i\colon
h(i)<\varepsilon\} \in J_\zeta\}$.

\item{(c)} for $\zeta<\zeta(*)$ and $A\in J_{\zeta+1}\setminus
J_\zeta$, the pair $(\prod\bar\lambda, J_\zeta+(\kappa\setminus A))$
(equivalently $\prod\bar\lambda\upharpoonright A,
J_\zeta\upharpoonright A)$) satisfies condition 1.5$(*)$ (case
$(\beta)$) hence its consequences, (in particular it satisfies the
weak $\pcf$-th for $\theta$).

\item{(d)} if $\kappa\not\in \cup_{\zeta<\zeta(*)} J_\zeta$ then
$(\prod\bar\lambda, \cup_{\zeta<\zeta(*)} J_\zeta)$ has true cofinality
$\theta$.
\endproclaim

\demo{Proof}
Straight. (We define $J_\zeta$ for $\zeta\leq \theta^+$ by clause (b)
for $\zeta=0$, $\zeta$ successor and as
$\bigcup\limits_{\varepsilon<\zeta} J_\varepsilon$ for $\zeta$ limit.
Clause $(c)$ holds by claim 4.4 below. It should be clear that
$J_{\theta^++1} = J_{\theta^+}$, and let $\zeta(*)=\min\{\zeta\colon
J_{\zeta+1}=\bigcup\limits_{\varepsilon<\zeta}J_\varepsilon\}$ so we are
left with checking clause (d). If $A\in J^+_{\zeta(*)}$, $h\in
\prod\limits_{i\in A}\lambda_i$, choose by induction on $\zeta<\theta$,
$\varepsilon(\zeta)<\theta$ increasing with $\zeta$ such that
$\{i<\kappa\colon h(i)\in (\varepsilon(\zeta), \varepsilon(\zeta+1))\in
J^+_{\zeta(*)}$. If we succeed we contradict $\theta\geq \wsat(I^*)$ as
$\theta$ is regular. So for some $\zeta<\theta$, $\varepsilon(\zeta)$ is
well defined but not $\varepsilon(\zeta+1)$. As $J_{\zeta(*)}=
J_{\zeta(*)+1}$, clearly $\{i<\kappa\colon h(i)\leq
\varepsilon(\zeta)\}=\kappa\ \mod\ J_{\zeta(*)}$. So
$g_\varepsilon(i)=\cases{ \varepsilon & if $\varepsilon<\lambda_i$\cr
0 & if $\varepsilon \geq \lambda_i$\cr}$ exemplifies $\tcf(\prod\bar
\lambda/ J_{\zeta(*)})=0$.
\sqed{4.2}

Now:

\proclaim{Conclusion 4.3}
Under the assumptions of 4.2, $I^*$ satisfies the pseudo $\pcf$-th (see
Definition 2.11(4)), hence $\cf(\prod\bar \lambda, <_{I^*}) =
\sup\pcf_{I^*}(\bar \lambda)$ (see 2.14).
\endproclaim

\proclaim{Claim 4.4}
Under the assumption of 4.2, if $J$ is an ideal on $\kappa$ extending
$I^*$ the following conditions are equivalent
\item{(a)} for some $h\in \prod\bar \lambda$, for every
$\varepsilon<\theta$ we have $\{i\in A: h(i)<\varepsilon\}\in J$
\item{(b)} $(\prod\bar \lambda, <_{J+(\kappa\setminus A)})$ is
$\theta^+$-directed.
\endproclaim

\demo{Proof}
$\underline{(a)\Rightarrow (b)}$

Let $f_\zeta\in \prod\bar \lambda$ for $\zeta<\theta$, we define $f^*\in
\prod\bar \lambda$ by
$$
f^*(i)=\sup\{ f_\zeta(i)+1: \zeta< h(i)\}.
$$
Now $f^*(i)<\lambda_i$ as $h(i)<\lambda_i=\cf(\lambda_i)$ and
$f_\zeta\restriction A <_J f^*\restriction A$ as $\{i\in A:
h(i)<\zeta\}\in J$.

\noindent
$\underline{(b)\Rightarrow (a)}$

Let $f_\zeta$ be the following function with domain $\kappa$:
$$
f_\zeta(i)= \cases{ \zeta & if $\zeta<\lambda_i$\cr
     0 & if $\zeta\geq \lambda_i$\cr}
$$
As $\lim\inf_{I^*} \geq \theta$, clearly $\varepsilon<\zeta \Rightarrow
f_\varepsilon<_{I^*} f_\zeta$ and of course $f_\zeta \in \prod\bar
\lambda$. By ourassumption (b) there is $h\in \prod \bar \lambda$ such
that $\zeta<\theta \Rightarrow f_\zeta \restriction A< h\restriction A\
\mod\ J$. Clearly $h$ is as required.
\sqed{4.4}

\newpage


\bigskip
\references{XXSh000}

\ref [CK]
C. C. Chang and H. J. Keisler,
\sl Model Theory,
\rm North Holland Publishing Company (1973).

\ref [Kn]
A. Kanamori,
\sl Weakly normal filters and irregular ultra-filter,
\rm Trans of A.M.S., {\bf 220} (1976) 393--396.

\ref [Ko]
S. Koppelberg,
{\sl Cardinalities of ultraproducts of finite sets},
The Journal of Symbolic Logic, {\bf 45} (1980) 574--584.

\ref [Kt]
J. Ketonen,
\sl Some combinatorial properties of ultra-filters,
\rm Fund Math. VII (1980) 225--235.

\ref [M]
J. D. Monk,
\sl Cardinal Function on Boolean Algebras,
\rm Lectures in Mathmatics, ETH Z\"urich, Bikh\"auser, Verlag,
Baser, Boston, Berlin, 1990.

\ref [Sh-b]
S. Shelah,
\sl Proper forcing
\rm Springer Lecture Notes, 940 (1982) 496+xxix.

\ref [Sh-g]
S. Shelah,
{\sl Cardinal Artihmetic,}
volume 29 of Oxford Logic Guides,
General Editors: Dov M. Gabbai, Angus Macintyre and Dana Scott,
{\rm Oxford Univertsity Press,} 1994.

\ref[Sh7]
S. Shelah,
{\sl On the cardinality of ultraproduct of finite sets},
Journal of Symbolic Logic, {\bf 35} (1970) 83--84.

\ref [Sh345]
S. Shelah,
\sl Products of regular cardinals and cardinal invariants of Boolean Algebra,
\rm Israel Journal of Mathematics, {\bf 70} (1990) 129--187.

\ref [Sh400a]
S. Shelah,
{\sl Cardinal arithmetic for skeptics},
{American Mathematical Society. Bulletin. New Series}, {\bf 26} (1992)
197--210. 

\ref [Sh 420]
S. Shelah,
\sl Advances in Cardinal Arithmetic,
\rm Proceedings of the Conference in Banff, Alberta, April 1991.

\ref [Sh430]
S. Shelah,
{\sl Further cardinal arithmetic},
Israel Journal of Mathematics, {\bf accepted}.

\ref [MgSh433]
M. Magidor and S. Shelah,
{\sl $\lambda _i$ inaccessible $> \kappa , \prod_{i<\kappa } \lambda
_i/D$ of order type $\mu ^+$}, preprint.

\ref[Sh580]
S. Shelah, 
{\it Strong covering without squares, the tree numbers},
in preparation.

\endreferences


\shlhetal

\shalom
\bye